\documentclass{amsart}
\usepackage{graphicx}
\usepackage{amssymb}
\usepackage{epstopdf}
\usepackage{nicefrac}
\usepackage{enumerate}
\usepackage{float}

\usepackage{tabularx}
\usepackage{hyperref}

\newtheorem{thm}{THEOREM}[section]
\newtheorem{conj}[thm]{CONJECTURE}
\newtheorem{cor}[thm]{COROLLARY}
\newtheorem{defn}[thm]{DEFINITION}

\newtheorem{lemma}[thm]{LEMMA}

\newtheorem{prop}[thm]{PROPOSITION}

\newtheorem{remark}[thm]{REMARK}

\newcommand{\ds}{\displaystyle}

\newcommand{\F}{{\mathcal F}} 
\newcommand{\wtF}{{\widetilde{\mathcal F}}} 
\newcommand{\G}{\Gamma}

\newcommand{\e}{{\epsilon}} 
\newcommand{\eM}{{\epsilon_{\fM}}} 
\newcommand{\eU}{{\epsilon_{\cU}}} 
\newcommand{\eFU}{{\epsilon^{\F}_{\cU}}} 
\newcommand{\eTU}{{\epsilon^{\cT}_{\cU}}} 

\newcommand{\ve}{\varepsilon}

\newcommand{\adj}{\text{\bf adj}}

\newcommand{\dF}{d_{\F}} 
\newcommand{\wtdF}{d_{\wtF}} 
\newcommand{\dM}{d_{\fM}} 

\newcommand{\dFU}{\delta^{\F}_{\cU}} 
\newcommand{\dTU}{\delta^{\cT}_{\cU}} 
\newcommand{\dT}{\delta^{\cT}} 

\newcommand{\wtkF}{{\kappa_{\wtF}}} 
\newcommand{\kL}{{\kappa_{L}}} 
\newcommand{\kX}{{\kappa_{X}}} 
\newcommand{\lF}{{\lambda_{\mathcal F}}} 
\newcommand{\rp}{\rho_{\pi}} 
\newcommand{\rt}{\rho_{\tau}} 
\newcommand{\rF}{{\lambda_{\mathcal F}^*}} 

\newcommand{\cGF}{\cG_{\F}} 
\newcommand{\GF}{\Gamma_{\F}} 

\newcommand{\PF}{{\rm Pen}_{\F}} 
\newcommand{\PX}{{\rm Pen}_{X}} 

\newcommand{\wtdp}{{\widetilde{d}\hspace{.8pt}'}}
\newcommand{\wtgp}{{\widetilde{g}\hspace{.8pt}'}}

\newcommand{\wtyp}{{\widetilde{y}\hspace{.8pt}'}}

\newcommand{\whvp}{{\widehat{v}\hspace{.5pt}'}}

\newcommand{\vzpj}{{\vec{z}_j\hspace{-2.6pt}'}}
\newcommand{\vzpk}{{\vec{z}_k\hspace{-2.6pt}'}}

\newcommand{\vwpj}{{\vec{w}_j\hspace{-2.6pt}'}}

\newcommand{\vop}{{\vec{\omega}\hspace{.5pt}'}}

\newcommand{\whrho}{{\widehat{\rho}}}
\newcommand{\whrhop}{\widehat{\rho}\hspace{1.8pt}'}
\newcommand{\wtrho}{{\widetilde{\rho}}}
\newcommand{\wtrhop}{\widetilde{\rho}\hspace{1.8pt}'}

\newcommand{\wtd}{{\widetilde{d}}}

\newcommand{\wtg}{{\widetilde{g}}}
\newcommand{\wth}{{\widetilde{h}}}

\newcommand{\wtw}{{\widetilde{w}}}
\newcommand{\wtx}{{\widetilde{x}}}
\newcommand{\wtz}{{\widetilde{z}}}
\newcommand{\wty}{{\widetilde{y}}}

\newcommand{\wtD}{\widetilde{D}}
\newcommand{\wtK}{\widetilde{K}}
\newcommand{\wtL}{{\widetilde L}}

\newcommand{\wtR}{\widetilde{R}}

\newcommand{\wtU}{{\widetilde U}}
\newcommand{\wtV}{\widetilde{V}}

\newcommand{\wtomega}{{\widetilde{\omega}}}
\newcommand{\wtpi}{{\widetilde{\pi}}}

\newcommand{\wttau}{{\widetilde{\tau}}}
\newcommand{\wtvp}{{\widetilde{\varphi}}}
\newcommand{\wtgamma}{{\widetilde{\gamma}}}
\newcommand{\wtlambda}{{\widetilde{\lambda}}}

\newcommand{\wtGamma}{{\widetilde{\Gamma}}}
\newcommand{\wtPhi}{{\widetilde{\Phi}}}

\newcommand{\wtsigma}{{\widetilde{\sigma}}}

\newcommand{\wtcH}{{\widetilde{\mathcal H}}}
\newcommand{\wtcN}{{\widetilde{\mathcal N}}}
\newcommand{\wtcP}{{\widetilde{\mathcal P}}}
\newcommand{\wtcM}{\widetilde{\mathcal M}}
\newcommand{\wtfN}{\widetilde{\mathfrak{N}}}
\newcommand{\wtfR}{\widetilde{\mathfrak{R}}}

\newcommand{\wtfT}{\widetilde{\mathfrak{T}}}
\newcommand{\wtfU}{\widetilde{\mathfrak{U}}}
\newcommand{\wtfZ}{\widetilde{\mathfrak{Z}}}

\newcommand{\wtcX}{\widetilde{\mathcal X}}
\newcommand{\wtexp}{\widetilde{\exp}}

\newcommand{\wtfB}{\widetilde{\mathfrak{B}}}

\newcommand{\whe}{{\widehat e}}
\newcommand{\whf}{{\widehat f}}

\newcommand{\whu}{{\widehat u}}
\newcommand{\whv}{{\widehat v}}

\newcommand{\whz}{{\widehat z}}

\newcommand{\whvarp}{{\widehat \varphi}}

\newcommand{\whK}{{\widehat K}}

\newcommand{\whR}{{\widehat R}}
\newcommand{\whU}{{\widehat U}}

\newcommand{\whcX}{{\widehat{\mathcal X}}}

\newcommand{\whXi}{{\widehat{\Xi}}}

\newcommand{\wthu}{\widetilde{\widehat{u}}}

\newcommand{\vl}{{\vec{\lambda}}}

\newcommand{\whfU}{\widehat{\mathfrak{U}}}

\newcommand{\oU}{{\overline{U}}}

\newcommand{\bA}{{\bf A}}
\newcommand{\bB}{{\bf B}}
\newcommand{\bC}{{\bf C}}
\newcommand{\bQ}{{\bf Q}}
\newcommand{\bU}{{\bf U}}
\newcommand{\bV}{{\bf V}}
\newcommand{\bW}{{\bf W}}

\newcommand{\mG}{{\mathbb G}}
\newcommand{\mK}{{\mathbb K}}

\newcommand{\mR}{{\mathbb R}}

\newcommand{\mZ}{{\mathbb Z}}

\newcommand{\cB}{{\mathcal B}}
\newcommand{\cC}{{\mathcal C}}

\newcommand{\cG}{{\mathcal G}}

\newcommand{\cH}{{\mathcal H}}
\newcommand{\cI}{{\mathcal I}}
\newcommand{\cJ}{{\mathcal J}}

\newcommand{\cM}{{\mathcal M}}
\newcommand{\cN}{{\mathcal N}}
\newcommand{\cO}{{\mathcal O}}
\newcommand{\cP}{{\mathcal P}}

\newcommand{\cS}{{\mathcal S}}
\newcommand{\cT}{{\mathcal T}}
\newcommand{\cU}{{\mathcal U}}
\newcommand{\cV}{{\mathcal V}}

\newcommand{\cX}{{\mathcal X}}

\newcommand{\fB}{{\mathfrak{B}}}
\newcommand{\fC}{{\mathfrak{C}}}
\newcommand{\fD}{{\mathfrak{D}}}

\newcommand{\fM}{{\mathfrak{M}}}
\newcommand{\fN}{{\mathfrak{N}}}

\newcommand{\fR}{{\mathfrak{R}}}
\newcommand{\fS}{{\mathfrak{S}}}
\newcommand{\fT}{{\mathfrak{T}}}
\newcommand{\fU}{{\mathfrak{U}}}

\newcommand{\fX}{{\mathfrak{X}}}

\newcommand{\fZ}{{\mathfrak{Z}}}

\input xy
\xyoption {all}

\newcommand{\vp}{{\varphi}}

\begin{document}

\title{Voronoi tessellations for matchbox manifolds}


\author{Alex Clark}
\thanks{AC and OL supported in part by EPSRC grant EP/G006377/1}
\address{Alex Clark, Department of Mathematics, University of Leicester, University Road, Leicester LE1 7RH, United Kingdom}
\email{adc20@le.ac.uk}

\author{Steven Hurder}
\address{Steven Hurder, Department of Mathematics, University of Illinois at Chicago, 322 SEO (m/c 249), 851 S. Morgan Street, Chicago, IL 60607-7045}
\email{hurder@uic.edu}

\author{Olga Lukina}
\address{Olga Lukina, Department of Mathematics, University of Leicester, University Road, Leicester LE1 7RH, United Kingdom}
\email{ol16@le.ac.uk}

\thanks{Version date: July 11, 2011; revised August 7, 2012}

\date{}


\begin{abstract}
Matchbox manifolds $\fM$ are a special class of foliated spaces, which includes as special examples   exceptional minimal sets of foliations, weak solenoids, suspensions of odometer and Toeplitz actions,  and tiling spaces associated to aperiodic tilings with finite local complexity. 
Some of these classes of examples are endowed with    an additional structure, that of a transverse foliation, consisting of a continuous family of Cantor sets transverse to the foliated structure. The purpose of this paper is to show that this transverse structure can be defined on all minimal matchbox manifolds. This follows from the construction of uniform stable Voronoi tessellations on a dense leaf, which is the main goal of this work. From this we define a foliated Delaunay triangulation of $\fM$, adapted to the dynamics of $\F$. The result is highly technical, but underlies the study of the basic topological structure of matchbox manifolds in general. Our methods  are unique in that we give the construction of the Voronoi tessellations for a complete Riemannian manifold $L$ of arbitrary dimension, with stability estimates.
\end{abstract}

\maketitle


\section{Introduction and Main Theorems} \label{sec-intro}

 An \emph{$n$-dimensional foliated space} $\fM$   is    a continuum  locally homeomorphic to a product of a disk in $\mathbb{R}^n$ and a Hausdorff separable topological space.    The leaves of the foliation $\F$ of $\fM$ are the
maximal connected components with respect to the fine topology on
$\fM$ induced by the plaques of the local product structure. 
A \emph{matchbox manifold} is a foliated space  such that the
local transverse models are totally disconnected, and the leaves have a smooth structure. Thus, a matchbox manifold is a continuum $\fM$ whose arc-components define a   smooth foliated structure on $\fM$. A matchbox manifold is \emph{minimal} if every leaf of $\F$ is dense in $\fM$.

The main result of this paper is  that a minimal matchbox manifold $\fM$ has an additional ``regularity property'', that it  always admits a   \emph{Cantor foliation} $\cH$, which is ``transverse'' to the foliation $\F$. 
The existence of the transverse foliation $\cH$ has a variety of applications for the study of these spaces, as it implies that  a minimal matchbox manifold $\fM$ is homeomorphic to  an inverse limit of ``covering maps'' between branched manifolds. 
The method of proof   uses the construction of stable Delaunay triangulations for the leaves of such foliations, and we give a careful development of this topic as well as discuss  the subtleties that arise for the cases when $n > 2$. 
 
The notion of a matchbox manifold is essentially the same as that of a \emph{lamination}. In low-dimensional topology, a lamination is a decomposition into leaves of a closed subset of a manifold; in holomorphic dynamics, Sullivan \cite[Appendix]{Sullivan1988} introduced Riemann surface laminations as compact topological spaces locally homeomorphic to a complex disk times a Cantor set. A similar notion is used by Lyubich and Minsky in \cite{LM1997}, and Ghys in \cite{Ghys1999}.
An embedding into a manifold is not assumed in the latter contexts, and a matchbox manifold is a lamination in this sense. 

A celebrated theorem of Bing   \cite{Bing1960} showed that if $X$ is    a homogeneous, circle-like
continuum that contains an arc, then either $X$ is homeomorphic to a circle, or to   a Vietoris solenoid.
In the course of proving this result,
Bing raised the question: If $X$ is a homogeneous
continuum, and if every proper subcontinuum of $X$ is an arc, must
$X$ then be a circle or a solenoid?
An affirmative answer to this question was given by Hagopian~\cite{Hagopian1977}, and
 subsequent alternate proofs    were given by Mislove and Rogers \cite{MR1989} and by  Aarts,  Hagopian and Oversteegen~\cite{AHO1991}. 
These authors called such spaces ``matchbox manifolds'', as  $\fM$ admits a covering by local coordinate charts
$U$ which are a product of an interval with a Cantor set, and so are intuitively   a ``box of matches.''  
The    first two authors showed in   \cite{ClarkHurder2011b}  that a homogeneous continua $\fM$ whose arc-connected components define a $n$-dimensional foliation of $\fM$  is homeomorphic to a generalized solenoid, for  $n \geq 1$. That is,  
 Bing's Theorem holds in the context   of    $n$-dimensional matchbox manifolds.   The proof makes use of Theorem~\ref{thm-foliated} below.

Foliated spaces were   introduced in the book by Moore and Schochet \cite{MS1988}, as the natural generalization of foliated manifolds for  which   the leafwise index theory of Connes could be developed \cite{Connes1994}. 
Our use of the ``foliated  space''  terminology follows that of   \cite{MS1988}, and also  \cite{CandelConlon2000}, and    the term ``matchbox manifold'' is used to distinguish a class of foliated spaces which  have totally disconnected transversals. 
  The foliation index theorem for matchbox manifolds was used to formulate   Bellisard's \emph{Gap Labeling Conjecture} \cite{Bellisard1992}, which    motivated the study of the topological properties of  the tiling space associated to an aperiodic tiling of $\mR^n$ of finite type.  The conjecture relates the $K$-theory of the $C^*$-algebra associated to the aperiodic potentials formed by translation of the operator on $\mR^n$ (see \cite{BHZ2000}), with the $K$-theory of the associated tiling space \cite{KellendonkPutnam2000}. The approaches to the solution of the \emph{Gap Labeling Conjecture} in the works  \cite{BBG2006,BOO2003,KamPut2003} motivated the study of the topology of tiling spaces, including the  work of this paper.

 Given a   repetitive, aperiodic tiling of $\mR^n$ with finite local complexity, the associated tiling space $\Omega$ is defined as the   closure   in an appropriate \emph{Gromov-Hausdorff topology} on the  space of translations of the tiling. Then the space $\Omega$ is a matchbox manifold in our sense (for example, see \cite{FS2009,SW2003,Sadun2008}.)
 One of the remarkable results of the theory of tilings of $\mR^n$ is the theorem of Anderson and Putnam and its extensions, that  the  tiling space $\Omega$  admits a presentation  as an inverse limit of a tower of branched flat manifolds \cite{AP1998,Sadun2003, Sadun2008}.  The proof   uses the combinatorial structure of the tiling, along with the observation that the choice of a basepoint $x_0 \in \mR^n$ defines by translation a family of Cantor set transversals to the foliation $\F$.

   Benedetti and Gambaudo \cite{BG2003}   proved   an equivariant version of the above results, in the case of tilings associated to a free action of a connected Lie group $G$ which again defines a foliated space $\Omega$. Part of their hypotheses  is that the action of $G$ on $\Omega$ preserves a transverse Cantor structure to the foliation $\F$ defined by the action. It is not clear if the  constructions in this paper extend to the $G$-equivariant case, generalizing the results  in \cite{BG2003}.

  Williams \cite{Williams1967,Williams1974} showed that an invariant hyperbolic set $\Omega$  for certain classes of Axiom A diffeomorphisms  are homeomorphic to   inverse limit systems. The inverse system is  defined by bonding maps which are immersions of branched manifolds, obtained from the expanding map obtained by restricting the action to the invariant hyperbolic set. These are the   \emph{Williams solenoids}, which are matchbox manifolds.  In this context, the diffeomorphism is assumed to preserve a hyperbolic splitting of the tangent bundle to the ambient manifold $M$ in a neighborhood of the smoothly embedded space $\Omega \subset M$. The leaves of $\F$ correspond to the expanding directions, while the stable, contracting directions define a foliation in a neighborhood of $\Omega$ which is transverse to $\F$. The restriction of this local stable foliation to $\Omega$ then defines a transverse Cantor structure for $\F$. These examples are a model for the study of the properties  of tiling spaces and matchbox manifolds.

 More generally, given a finitely generated group $\G$ and a  minimal topological action on a Cantor set $K$, the suspension construction yields a matchbox manifold, whose leaves are covered by a simply connected manifold of which $\G$ acts freely and cocompactly. 
 Such examples arise naturally in the   study of odometers and Toeplitz flows \cite{CoPe2008,Do2005},  which includes the class of McCord solenoids \cite{McCord1965,Schori1966}.   Suspension examples admit a fibration to a base space, and the fibers define a Cantor fibration  with fibers transverse to $\F$.
 A  variety of further examples of matchbox manifolds can be found in   the works \cite{Blanc2001,Blanc2003,Ghys1999,Lukina2011a}. 
Finally, the continua which arise as minimal sets for flows   or foliations, and have totally disconnected transverse model spaces, are examples of matchbox manifolds as well.

For the special cases above where $\fM$ is   a tiling space, a generalized or Williams solenoid, or a suspension, there is a    continuous family of local transversals to the foliation $\F$ on $\fM$. The collection of these transversals define a \emph{transverse Cantor foliation $\cH$} of $\fM$, whose ``leaves'' are   Cantor sets. This concept is  related to the notion of Smale spaces, introduced by Ruelle \cite{Ruelle1988} and developed by  Putnam \cite{Putnam1996,PS1999} and  Wieler \cite{Wieler2012a,Wieler2012b}. See   section~\ref{sec-cantor} for precise statements.

For a general matchbox manifold $\fM$,   its leaves are    Riemannian manifolds with uniform control over their geometry, but there is essentially no control over their ``transverse geometry''. For example, the transverse holonomy of the foliation $\F$ on $\fM$ need not be Lipshitz continuous, so that  many  usual techniques from foliation theory which require some degree of transverse regularity do not apply. 
The  main goal of this paper  is to show the  existence  of a transverse Cantor foliation $\cH$ for a minimal matchbox manifold $\fM$, where $\cH$ is ``well-adapted'' to the dynamical properties of the foliation $\F$, so that $\fM$ admits a covering with bi-foliated charts. 
Of equal importance though, is the method of proof using leafwise Voronoi tessellation and the associated Delaunay triangulations, as discussed below, which provides   a foundation for the general study of these spaces.

Our first   result  shows the existence of transverse Cantor foliations for matchbox manifolds with equicontinuous dynamics.

\begin{thm}\label{thm-foliated}
Let $\fM$ be an equicontinuous matchbox manifold.
Then there exists a transverse Cantor foliation $\cH$ on $\fM$ such that the projection to the leaf space $\fM \to \fM/\cH \cong M$ is a Cantor bundle map over a compact manifold $M$.
\end{thm}

The most general application of our techniques is to prove an analog of the  ``Long Box Lemma''  from the study of   flows on continua (see \cite[Lemma~5.2]{FO1996}, and  also  \cite{AHO1991,AO1991,AO1995,FO2002}). This lemma   states that every \emph{connected, contractible  orbit segment} $K$ in a $1$-dimensional matchbox manifold is contained in a bi-foliated open  subset $\fN_K$ of $\fM$, which is the ``long box'' neighborhood of $K$.
We generalize this concept to $n$-dimensional matchbox manifolds, for $n \geq 1$.

\begin{defn}\label{def-LBB}
Let $\fM$ be a matchbox manifold, $x \in \fM$ a basepoint, $L_x \subset M$ the leaf through $x$, and $\wtL_x$ the holonomy covering of $L_x$. We say that $K_x \subset L_x$ is a \emph{proper base} if $K_x$ is a union of closed foliation plaques with $x \in K_x$, and there is a \emph{connected} compact subset  
  $\wtK_x \subset \wtL_x$  such that the composition
$\ds \iota_x \colon \wtK_x \subset \wtL_x \to L_x \subset \fM$
is injective with image $K_x$.
\end{defn}
The hypotheses imply that $K_x$ is path connected, and the holonomy of $\F$ along any path in $K_x$ is trivial.  
The   ``Big Box''  concept  is analogous to a strong form of the Reeb Stability theorem for foliations of compact manifolds.

\begin{thm}[Big Box]\label{thm-tessel}
Let $\fM$ be a minimal matchbox manifold, $x \in \fM$ and $K_x \subset L_x$ a proper base. Then there exists a clopen transversal $V_x \subset \fM$ containing $x$, and a foliated homeomorphic inclusion $\Phi \colon  K_x \times V_x \to \fN_{K_x} \subset \fM$ such that the images $\ds \Phi \left( \{ \{y\} \times V_x \mid y \in K_x\}\right)$ form a continuous family of Cantor transversals for $\F | \fN_{K_x}$.
\end{thm}
The main part of the conclusion is the foliated product structure on the image $\ds \fN_{K_x} \subset \fM$. This implies that $\ds \fN_{K_x}$ has a natural projection map to the submanifold $K_x \subset L_x$.

 In the case of one-dimensional matchbox manifolds, Theorem~\ref{thm-tessel} reduces to the Long Box Lemma \cite[Lemma~5.2]{FO1996}, as   a contractible line segment $K$   has no holonomy. The requirement that $K_x$ contains no loops with holonomy is essential. 

Theorem~\ref{thm-tessel} and the techniques used to prove the result    are used in  \cite{CHL2012a}  to  prove:
\begin{thm}\label{thm-cantorfol}
Let $\fM$ be a minimal matchbox manifold. Then $\fM$ admits a Cantor foliation transverse to $\F$, such that the   leaf space $\fM/\cH$ is a branched compact manifold.
\end{thm}

Given the existence of the transverse Cantor foliation $\cH$ as in  either  Theorems~\ref{thm-foliated} and \ref{thm-cantorfol}, then the fact that the leaves of $\cH$ are totally disconnected allows decomposing the foliation $\cH$ into a finite union of subfoliations, whose quotient spaces again have the structure of a (branched) manifold. Continuing this process recursively, so that the ``leaves'' of the subfoliations have diameter tending to zero, one obtains the following two applications.   
Details  are given in the cited papers below.
 
 \begin{thm}\label{thm-CHb} \cite{ClarkHurder2011b}
 Let $\fM$ be an equicontinuous matchbox manifold. Then $\fM$ is a foliated space, homeomorphic to a weak solenoid. That is, $\fM$ is homeomorphic to the inverse limit of an infinite chain of proper covering maps between compact manifolds. 
  \end{thm}
 
 \medskip
 
 \begin{thm} \label{thm-CHLa} \cite{CHL2012a}
 Let $\fM$ be a minimal matchbox manifold. Then $\fM$ is   homeomorphic  to the inverse limit of an infinite chain of  branched covering maps  between compact branched manifolds. 
  \end{thm}

The   results of Theorems~\ref{thm-foliated}, \ref{thm-CHb} and \ref{thm-CHLa}, all have   applications to the study of the topological invariants of the space $\fM$,  analogous to those results for tiling spaces and minimal actions of $\mZ^n$ on Cantor sets,  to the dynamics of the foliation $\F$ as is given in Gambaudo and Martens \cite{GamMart2006}, and other properties of foliated spaces as discussed in \cite{Hurder2011a}.  
 
We mention one other direction of research giving an application of the techniques of this paper. 
The sequence of papers by   the authors  Forrest \cite{Forrest2000},    Giordano, Matui, Putnam and  Skau \cite{GPS2004,GMPS2010,Matui2008,Putnam2010} and Phillips \cite{Phillips2005} study   \emph{topological orbit equivalence}   for minimal  $\mZ^n$-actions on Cantor sets, culminating in a proof that they are \emph{affable}.  
  The methods we develop here for the proof of Theorem~\ref{thm-tessel} also makes it possible to consider extensions of their results   for minimal $\mZ^n$-actions to much more general cases. We refer to their works for the definition of affable equivalent relations, which is a continuous analog of the hyperfinite property in measurable dynamics.

\begin{conj}\label{conj-affable}
 Let $\fM$ be a minimal matchbox manifold, which is amenable in the sense of sense of Anantharaman-Delaroche and Renault \cite{AnRe2000}. Then $\fM$ is affable. 
\end{conj}
It seems   likely that using a combination of the methods of 
  \cite{Forrest2000} with those developed  in  \cite{GMPS2010}, it should be possible to prove the special case of Conjecture~\ref{conj-affable}  for the case of 
   minimal matchbox manifolds with leaves of polynomial growth,  which would be analogous to the proof by   Series  in  \cite{Series1979} that   a  foliated manifold with leaves of   polynomial growth is  hyperfinite.

 \section{Overview of proof}\label{sec-overview}

We give some comments on the goals of this work, and then give an overview of the paper and the proofs of the above theorems. 
     
Our first remark is that in the case where $\F$ is defined by a flow, as in the seminal works \cite{AHO1991,AO1991,AO1995,AM1988}, the construction of $\cH$ can be easily accomplished by the choice of a local transversal to the flow, and then using the flow to suitably translate it to all points of the space. However, when the leaves of $\F$ have dimension greater than one, this method does not apply. Even in the case where $\F$ is defined by an action of a connected Lie group $G$ as in \cite{BG2003}, it is part of their hypotheses that translation by $G$   yields a transverse Cantor foliation.
 
In the case where $\F$ is a smooth foliation of a compact manifold, it is generally not possible to construct a transverse foliation to $\F$, unless the codimension of $\F$ is one, and then the unit normal vector field to the leaves defines such a foliation. Thus, the existence of $\cH$  fundamentally requires that the transverse geometry be totally disconnected. 
 
 In the case where $\fM$ is   embedded  as a minimal set for a smooth foliation $\wtF$ of a compact manifold $M$, one approach to the construction of $\cH$ would be to construct a transverse foliation $\wtcH$  to $\wtF$ on some open neighborhood of $\fM$ in $M$, and restrict $\wtcH$ to $\fM$. Even with these very strong assumptions, the requirement that $\wtcH$ be transverse to $\wtF$ demands a local inductive approach, based for example on the construction of a smooth transverse triangulation to $\wtF$ as given by Benameur in \cite{Benameur1997}, using the ideas of Thurston's ``Jiggling Lemma'' in the Appendix to \cite{Thurston1974}. 
In general,  a matchbox manifold $\fM$ need not have a smooth embedding into some Euclidean space $\mR^N$ for $N \gg 0$, \emph{such that the normal bundle to the smoothly embedded leaves is well-behaved}. Thus, it is problematic whether transversality methods such as used by Benameur can be applied to the problem. 
 
The approach in this paper to the construction of $\cH$ on $\fM$  is based on using arbitrarily small but discrete approximations to the geometry of the leaves of $\F$,  and deriving corresponding estimates of the uniform continuity of the ``normal direction''. As the transversal geometry to $\F$ is totally disconnected, the concept of normal to a leaf must be defined ``discretely'' as well. 

Most of the work in the construction of $\cH$  is to show that for each leaf $L_x \subset \fM$, there exists a Voronoi set $\cN_x$ which is sufficiently fine, and the points of $\cN_x$ are ``strongly stable'' as the point $x$ varies transversally over a small but fixed transverse clopen set. The strongly stable criterion, as  made precise later in the paper, implies that the leafwise Delaunay triangulations   of the leaves of $\F$ associated to the Voronoi  sets are transversely stable, and so can be used to define local  affine transverse coordinates, which are used to define the transverse foliation $\cH$. In this way, the   more geometric construction one obtains in the     case where $\fM$ is a minimal set for a smooth foliation, is extended to the general case of matchbox manifolds, without assuming any additional regularity on the holonomy maps of $\F$. 

The second remark is that while the techniques used are based     on ``elementary methods'', they are quite technical. We give effective estimates for the construction of Delaunay triangulations in the case where the leaves of $\F$ are general Riemannian manifolds, and are not assumed to be Euclidean.  This is required for our    construction of ``stable'' Delaunay triangulations, as the stability property  is fundamentally more subtle in dimensions greater than two, than for the one and two dimensional cases.   This work  gives the  construction in full detail, as it does not seem to be dealt with in the literature, yet is the foundation for a variety of other results. 

  For each foliated coordinate chart, $\oU_i \subset \fM$, there is a natural ``vertical'' foliation whose leaves are the images of the vertical transversals
defined by the coordinate charts. The problem is that on the overlap of two charts, these vertical foliations need not match up, as the only requirement on a foliation chart for $\F$ is that the horizontal plaques in each chart ``match up''. The exception is when $\fM$ is given with a fibration structure, then the coordinates can be chosen to be adapted to the fibration structure, and so the fibers of the bundle restrict to transversals in each chart which  are compatible on overlapping charts.
The   idea of our construction is   to subdivide the topological space $\fM$ by dividing it into ``arbitrarily small'' coordinate boxes, where the problems to be solved become ``almost linear''.

The structure of the paper is as follows.
First,  we present the basic concepts and dynamical properties of matchbox manifolds as required. For the preliminary  results in  sections~\ref{sec-concepts} through section~\ref{sec-microbundles}, the proofs are often  omitted, as they can be found in  \cite{ClarkHurder2011b}. The exception is when details of the proof or the notations introduced are necessary for later development, then they are briefly outlined.

Section~\ref{sec-concepts} gives definitions and notations. Then section~\ref{sec-holonomy} introduces the holonomy pseudogroup, and gives some basic technical properties of holonomy maps.
Section~\ref{sec-mme} recalls important classical definitions from topological dynamics, adapted to the case of matchbox manifolds, and gives several results concerning the dynamical properties of matchbox manifolds. Section~\ref{sec-microbundles} develops the analog of foliated microbundles and the \textit{Reeb Stability Theorem} for matchbox manifolds. Then in section~\ref{sec-cantor}, the important notion of a \emph{transverse Cantor foliation $\cH$} for $\F$  is introduced.

Sections~\ref{sec-VT} and \ref{sec-DSC} give the classical constructions of Voronoi and Delaunay triangulations in the context of Riemannian manifolds. In section~\ref{sec-FVS}, we extend these concepts from a single leaf, to a ``parametrized version'' which applies uniformly to the leaves of a matchbox manifold $\fM$, and introduce the notion of a nice stable transversal.  In section~\ref{sec-proofstableappr} we show how to obtain a transverse Cantor foliation using the Delaunay triangulation derived from    a nice stable transversal.

Sections~\ref{sec-euclidean} and \ref{sec-inequalities} recall the classical properties  of  Voronoi tessellations and Delaunay triangulations in Euclidean space. The goal is to  establish the framework and  some estimates  that are required in the subsequent constructions in the Riemannian context. 

Sections~\ref{sec-approx}  and \ref{sec-distortions} develop the ``micro-local'' Riemannian geometry of the leaves. The point is to    give estimates on the distortion from Euclidean geometry in the local adapted Gauss coordinate systems, which are used to chose the radius of charts sufficiently small so that the leafwise Delaunay triangulations will be defined and stable.  This development  requires standard results  of local Riemannian geometry, as in Bishop and Crittenden \cite{BC1964} or Helgason \cite{Helgason1978}. We establish various   \emph{a priori}  estimates,  as in section~\ref{sec-constants}, which seem to be   unavoidable in order to prove the stability   of the leafwise Delaunay triangulations.

Finally,   sections~\ref{sec-existence} and \ref{sec-nicestable} give the inductive construction   for nice stable transversals. The procedure followed invokes all the previous preparations, including the constants defined  in section~\ref{sec-constants}. The resulting arguments are the most technical of this work. 
The paper concludes with   section~\ref{sec-proofs} where we apply these results to prove  Theorems~\ref{thm-foliated} and  \ref{thm-tessel}.

We  comment on the various notations used in this paper. In general, roman letters indicate properties of a manifold, such as sets in the leaves. Greek letters are used to denote the more conventional notations, such as the Lebesgue number $\eU$ or the equicontinuous constant $\ds \dTU$, 
 or leafwise estimates such at $\lF$.
Capital script letters generally indicate some set associated with the leafwise   sets, such as the plaque chains $\ds \cP_{\cI}(\xi)$, or the  simplicial cones $\ds  \cC_{\Delta}(z)$. 
Fraktur letters denote sets  associated with the matchbox manifold $\fM$, such as the transversals $\fT_i$ or the Voronoi cylinders $\ds  \fC^{\ell}_{\wtfR}$.

\section{Foliated spaces} \label{sec-concepts}

We introduce the basic concepts of foliated spaces and matchbox manifolds. Further discussion with examples can be found in \cite[Chapter 11]{CandelConlon2000}, \cite[Chapter 2]{MS1988} and the  papers \cite{ClarkHurder2011a,ClarkHurder2011b}.  Recall that a \emph{continuum} is a compact connected metrizable space. 

\begin{defn} \label{def-fs}
A \emph{foliated space of dimension $n$} is a   continuum $\fM$, such that  there exists a compact separable metric space $\fX$, and
for each $x \in \fM$ there is a compact subset $\fT_x \subset \fX$, an open subset $U_x \subset \fM$, and a homeomorphism defined on the closure
$\vp_x \colon \oU_x \to [-1,1]^n \times \fT_x$ such that $\vp_x(x) = (0, w_x)$ where $w_x \in int(\fT_x)$. 
Moreover, it is assumed that each $\vp_x$  admits an extension to a foliated homeomorphism  
$\whvarp_x \colon \whU_x \to (-2,2)^n \times \fT_x$ where $\oU_x \subset \whU_x$.
\end{defn}
The subspace  $\fT_x$ of $\fX$ is called   the \emph{local transverse model} at $x$.

Let $\pi_x \colon \oU_x \to \fT_x$ denote the composition of $\vp_x$ with projection onto the second factor.

For $w \in \fT_x$ the set $\cP_x(w) = \pi_x^{-1}(w) \subset \oU_x$ is called a \emph{plaque} for the coordinate chart $\vp_x$. We adopt the notation, for $z \in \oU_x$, that $\cP_x(z) = \cP_x(\pi_x(z))$, so that $z \in \cP_x(z)$. Note that each plaque $\cP_x(w)$ for $w \in \fT_x$ is given the topology so that the restriction $\vp_x \colon \cP_x(w) \to [-1,1]^n \times \{w\}$ is a homeomorphism. Then $int (\cP_x(w)) = \vp_x^{-1}((-1,1)^n \times \{w\})$.

Let $U_x = int (\oU_x) = \vp_x^{-1}((-1,1)^n \times int(\fT_x))$.
Note that if $z \in U_x \cap U_y$, then $int(\cP_x(z)) \cap int( \cP_y(z))$ is an open subset of both
$\cP_x(z) $ and $\cP_y(z)$.
The collection of sets
$$\cV = \{ \vp_x^{-1}(V \times \{w\}) \mid x \in \fM ~, ~ w \in \fT_x ~, ~ V \subset (-1,1)^n ~ {\rm open}\}$$
forms the basis for the \emph{fine topology} of $\fM$. The connected components of the fine topology are called \emph{leaves}, and define the foliation $\F$ of $\fM$.
For $x \in \fM$, let $L_x \subset \fM$ denote the leaf of $\F$ containing $x$.

Note that in Definition~\ref{def-fs}, the collection of transverse models
$\{\fT_x \mid x \in \fM\}$ need not have union equal to $\fX$. This is similar to the situation for a smooth foliation of codimension $q$, where each foliation chart projects to an open subset of $\mR^q$, but the collection of images need not cover $\mR^q$.
\begin{defn} \label{def-sfs}
A \emph{smooth foliated space} is a foliated space $\fM$ as above, such that there exists a choice of local charts $\vp_x \colon \oU_x \to [-1,1]^n \times \fT_x$ such that for all $x,y \in \fM$ with $z \in U_x \cap U_y$, there exists an open set $z \in V_z \subset U_x \cap U_y$ such that $\cP_x(z) \cap V_z$ and $\cP_y(z) \cap V_z$ are connected open sets, and the composition
$$\psi_{x,y;z} \equiv \vp_y \circ \vp_x ^{-1}\colon \vp_x(\cP_x (z) \cap V_z) \to \vp_y(\cP_y (z) \cap V_z)$$
is a smooth map, where $\vp_x(\cP_x (z) \cap V_z) \subset \mR^n \times \{w\} \cong \mR^n$ and $\vp_y(\cP_y (z) \cap V_z) \subset \mR^n \times \{w'\} \cong \mR^n$. The leafwise transition maps $\psi_{x,y;z}$ are assumed to depend continuously on $z$ in the $C^{\infty}$-topology on maps between subsets of $\mR^n$.
\end{defn}

A map $f \colon \fM \to \mR$ is said to be \emph{smooth} if for each flow box
$\vp_x \colon \oU_x \to [-1,1]^n \times \fT_x$ and $w \in \fT_x$ the composition
$y \mapsto f \circ \vp_x^{-1}(y, w)$ is a smooth function of $y \in (-1,1)^n$, and depends continuously on $w$ in the $C^{\infty}$-topology on maps of the plaque coordinates $y$. As noted in \cite{MS1988} and \cite[Chapter 11]{CandelConlon2000}, this allows one to define smooth partitions of unity, vector bundles, and tensors for smooth foliated spaces. In particular, one can define leafwise Riemannian metrics. We recall a standard result, whose proof for foliated spaces can be found in \cite[Theorem~11.4.3]{CandelConlon2000}.
\begin{thm}\label{thm-riemannian}
Let $\fM$ be a smooth foliated space. Then there exists a leafwise Riemannian metric for $\F$, such that for each $x \in \fM$, $L_x$ inherits the structure of a complete Riemannian manifold with bounded geometry, and the Riemannian geometry of $L_x$ depends continuously on $x$. In particular, each leaf $L_x$ has the structure of a complete Riemannian manifold with bounded geometry. \hfill $\Box$
\end{thm}
Bounded geometry implies, for example, that for each $x \in \fM$, there is a leafwise exponential map
$\exp^{\F}_x \colon T_x\F \to L_x$ which is a surjection, and the composition $\exp^{\F}_x \colon T_x\F \to L_x \subset \fM$ depends continuously on $x$ in the compact-open topology on maps.

\begin{defn} \label{def-mm}
A \emph{matchbox manifold} is a continuum with the structure of a
smooth foliated space $\fM$, such that  the
transverse model space $\fX$ is totally disconnected, and for each $x \in \fM$, $\fT_x \subset \fX$ is a clopen subset.
\end{defn}
All matchbox manifolds  are assumed to be smooth with a given leafwise Riemannian metric.

\subsection{Metric properties and regular covers}\label{sec-metricregular}
We first establish the local properties of a matchbox manifold, which are codified by the definition of a \emph{regular covering} of $\fM$. 
One particular property to note is the  \emph{strong local convexity} for the plaques in the leaves.

Another nuance about the study of matchbox manifolds, is that for given $x \in \fM$,
the neighborhood $\oU_x$ in Definition~\ref{def-fs} need not be ``local''. As the transversal model $\fT_x$ is totally disconnected, the set $\oU_x$ is  not   connected, and \emph{a priori} its connected components need not be contained in a suitably small metric ball around $x$.
The technical procedures described in detail in \cite[\S 2.1 - 2.2]{ClarkHurder2011b} ensure that we can always choose local charts for $\fM$ to have a \emph{uniform locality property}, as well as other metric regularity properties as discussed below.

Let $\dM \colon \fM \times \fM \to [0,\infty)$ denote the metric on $\fM$, and $d_{\fX} \colon \fX \times \fX \to [0,\infty)$ the metric on $\fX$.

For $x \in \fM$ and $\e > 0$, let $D_{\fM}(x, \e) = \{ y \in \fM \mid \dM(x, y) \leq \e\}$ be the closed $\e$-ball about $x$ in $\fM$, and $B_{\fM}(x, \e) = \{ y \in \fM \mid \dM(x, y) < \e\}$ the open $\e$-ball about $x$.

Similarly, for $w \in \fX$ and $\e > 0$, let $D_{\fX}(w, \e) = \{ w' \in \fX \mid d_{\fX}(w, w') \leq \e\}$ be the closed $\e$-ball about $w$ in $\fX$, and $B_{\fX}(w, \e) = \{ w' \in \fX \mid d_{\fX}(w, w') < \e\}$ the open $\e$-ball about $w$.

Each leaf $L \subset \fM$ has a complete path-length metric, induced from the leafwise Riemannian metric:
$$\dF(x,y) = \inf \left\{\| \gamma\| \mid \gamma \colon [0,1] \to L ~{\rm is ~ piecewise ~~ C^1}~, ~ \gamma(0) = x ~, ~ \gamma(1) = y ~, ~ \gamma(t) \in L \quad \forall ~ 0 \leq t \leq 1\right\}$$
  where $\| \gamma \|$ denotes the path-length of the piecewise $C^1$-curve $\gamma(t)$. If $x,y \in \fM$   are not on the same leaf, then set $\dF(x,y) = \infty$. For each $x \in \fM$ and $r > 0$, let $D_{\F}(x, r) = \{y \in L_x \mid \dF(x,y) \leq r\}$.

Note that the metric $\dM$ on $\fM$ and the leafwise metric $\dF$ have no relation, beyond their relative continuity properties. The metric $\dM$ is essentially   just used to define the metric topology on $\fM$, while the metric $\dF$ depends on an independent  choice of the Riemannian metric on leaves.

For each $x \in \fM$, the  {Gauss Lemma} implies that there exists $\lambda_x > 0$ such that $D_{\F}(x, \lambda_x)$ is a \emph{strongly convex} subset for the metric $\dF$. That is, for any pair of points $y,y' \in D_{\F}(x, \lambda_x)$ there is a unique shortest geodesic segment in $L_x$ joining $y$ and $y'$ and  contained in $D_{\F}(x, \lambda_x)$ (cf. \cite[Chapter 3, Proposition 4.2]{doCarmo1992}, or \cite[Theorem 9.9]{Helgason1978}). Then for all $0 < \lambda < \lambda_x$ the disk $D_{\F}(x, \lambda)$ is also strongly convex. As $\fM$ is compact and the leafwise metrics have uniformly bounded geometry, we obtain:

\begin{lemma}\label{lem-stronglyconvex}
There exists $\lF > 0$ such that for all $x \in \fM$, $D_{\F}(x, \lF)$ is strongly convex.
\end{lemma}

If $\F$ is defined by a flow without periodic points, so
that every leaf is diffeomorphic to $\mR$, then the entire leaf is
strongly convex, and $\lF > 0$ can be chosen arbitrarily. For
foliations with leaves of dimension $n > 1$, the constant $\lF$ must
be less than one half the injectivity radius     of each leaf.

The following proposition summarizes results in \cite[\S 2.1 - 2.2]{ClarkHurder2011b}.

\begin{prop}\label{prop-regular}\cite{ClarkHurder2011b}
For a smooth foliated space $\fM$, given $\eM > 0$, there exist   $\lF>0$ and a choice of local charts $\vp_x \colon \oU_x \to [-1,1]^n \times \fT_x$ with the following properties:
\begin{enumerate}
\item For each $x \in \fM$, $U_x \equiv int(\oU_x) = \vp_x^{-1}\left( (-1,1)^n \times B_{\fX}(w_x, \e_x)\right)$, where $\e_x>0$.
\item Locality: for all $x \in \fM$,  each $\oU_x \subset B_{\fM}(x, \eM)$.
\item Local convexity: for all $x \in \fM$ the plaques of $\vp_x$ are leafwise strongly convex subsets with diameter less than $\lF/2$. That is,   there is a unique shortest geodesic segment joining any two points in a plaque, and the entire geodesic segment is contained in the plaque.
\end{enumerate}
A \emph{regular covering} of $\fM$ is one that satisfies   these conditions. 
\end{prop}

By a standard argument, there exists a finite collection  $\{x_1, \ldots , x_{\nu}\} \subset \fM$ where $\vp_{x_i}(x_i) = (0 , w_{x_i})$ for $w_{x_i} \in \fX$,  and regular foliation charts $\vp_{x_i} \colon \oU_{x_i} \to [-1,1]^n \times \fT_{x_i}$
satisfying the conditions of Proposition~\ref{prop-regular},   which form an open covering of $\fM$. Moreover, without loss of generality, we can    impose a uniform size restriction on the plaques of each chart. Without loss of generality, we can assume  there exists $0 < \dFU < \lF/4$ so that for all $1 \leq i \leq \nu$ 
and $\omega \in \fT_i$ with plaque ``center point'' $x_{\omega} = \vp_{x_i}^{-1}(0 , \omega)$, then 
the   plaque  for $\vp_{x_i}$ through $x_{\omega}$ satisfies the uniform estimate of diameters:
\begin{equation}\label{eq-Fdelta}
D_{\F}(x_{\omega} , \dFU/2) ~ \subset ~ \cP_i(\omega) ~ \subset ~ D_{\F}(x_{\omega} , \dFU) 
\end{equation}
For each $1 \leq i \leq \nu$ the set $\cT_{x_i} =  \vp_i^{-1}(0 , \fT_i)$ is a compact transversal to $\F$. Again, without loss of generality, we can assume   that the transversals 
$\ds \{ \cT_{x_1} , \ldots , \cT_{x_{\nu}} \}$ are pairwise disjoint, so there exists a constant $0 < e_1 <  \dFU$ such that 
\begin{equation}\label{eq-e1}
\dF(x,y) \geq e_1 \quad {\rm for} ~ x \ne y ~,   x \in \cT_{x_i} ~ , ~ y \in \cT_{x_j} ~ , ~ 1 \leq i, j \leq \nu
\end{equation}
In particular, this implies that 
 the centers of disjoint plaques on the same leaf are separated by distance at least $e_1$.
 
Given a fixed choice of foliation covering as above, 
we simplify the notation as follows.
For $1 \leq i \leq \nu$, set $\oU_i = \oU_{x_i}$, $U_i = U_{x_i}$, and $\e_i = \e_{x_i}$.
Let $\cU = \{U_{1}, \ldots , U_{\nu}\}$ denote the corresponding open covering of $\fM$, with  coordinate maps
$$
\vp_i = \vp_{x_i} \colon \oU_i \to [-1,1]^n \times \fT_i \quad , \quad
\pi_i = \pi_{x_i} \colon \oU_i \to \fT_i \quad , \quad \lambda_i \colon \oU_i \to [-1,1]^n.
$$
For $z \in \oU_i$, the plaque of the chart $\vp_i$ through $z$ is denoted by $\cP_i(z) = \cP_i(\pi_i(z)) \subset \oU_i$. Note that the restriction $\lambda_i \colon \cP_i(z) \to [-1,1]^n$ is a homeomorphism onto.
Also, define sections
\begin{equation}\label{eq-taui}
\tau_i \colon \fT_i \to \oU_i ~ , ~ {\rm defined ~ by} ~ \tau_i(\xi) = \vp_i^{-1}(0 , \xi) ~ , ~ {\rm so ~ that} ~ \pi_i(\tau_i(\xi)) = \xi.
\end{equation}
Then $\cT_i = \cT_{x_i}$ is the image of $\tau_i$ and we let $\cT = \cT_1 \cup \cdots \cup \cT_{\nu} \subset \fM$ denote their disjoint union.

Let $\ds \fT_* = \fT_1 \cup \cdots \cup \fT_{\nu} \subset ~ \fX$; note that $\fT_*$ is compact, and if each $\fT_i$ is totally disconnected, then $\fT_*$ will also be totally disconnected.

We assume in the following that a finite regular covering $\cU$ of $\fM$ as above has been chosen.

\subsection{Foliated maps}
A map $f \colon \fM \to \fM'$ between foliated spaces is said to be a \emph{foliated map} if the image of each leaf of $\F$ is contained in a leaf of $\F'$. If $\fM'$ is a matchbox manifold, then each leaf of $\F$ is path connected, so its image is path connected, hence must be contained in a leaf of $\F'$. Thus, 
\begin{lemma} \label{cor-foliated1}
Let $\fM$ and $\fM'$ be matchbox manifolds, and $h \colon \fM' \to \fM$ a continuous map. Then $h$ maps the leaves of $\F'$ to leaves of $\F$. In  particular, any homeomorphism $h \colon \fM \to \fM$ of a matchbox manifold is a foliated map. \hfill $\Box$
\end{lemma}

A \emph{leafwise path}, or more precisely an \emph{$\F$-path}, is a continuous map $\gamma \colon [0,1] \to \fM$ such that there is a leaf $L$ of $\F$ for which $\gamma(t) \in L$ for all $0 \leq t \leq 1$. 
If $\fM$ is a matchbox manifold, and $\gamma \colon [0,1] \to \fM$ is continuous, then   $\gamma$ is a leafwise path.

\subsection{Local estimates}\label{subsec-locestimates}
We next introduce a number of constants based on the above choices, which will be used throughout the paper when making metric estimates.

Let $\eU > 0$ be a Lebesgue number for the covering $\cU$. That is, given any $z \in \fM$ there exists some index $1 \leq i_z \leq \nu$ such that the open metric ball $B_{\fM}(z, \eU) \subset U_{i_z}$.

The local projections $\pi_i \colon \oU_i \to \fT_i$ and sections
$\tau_i \colon \fT_i \to \oU_i$ are continuous maps of compact spaces, so admit uniform metric estimates as follows.
\begin{lemma}\cite{ClarkHurder2011b}\label{lem-modpi}
There exists a continuous increasing function $\rp$ (the \emph{modulus of continuity} for the projections $\pi_i$) such that:
\begin{equation}\label{eq-modpi}
\forall ~ 1 \leq i \leq \nu ~ \text{and}~ x,y \in \oU_i \quad , \quad \dM(x,y) <\rp(\e) ~ \Longrightarrow ~ d_{\fX}(\pi_i(x), \pi_i(y)) < \e ~.
\end{equation}
\end{lemma}
\proof
Set $\ds\rp(\e) = \min \left\{\e, \min \left\{ \dM(x,y) \mid 1 \leq i \leq \nu ~ , ~ x,y \in \oU_i ~ , ~ d_{\fX}(\pi_i(x), \pi_i(y)) \geq \e\right\}\right\}$.
\endproof

\begin{lemma}\cite{ClarkHurder2011b}\label{lem-modtau}
There exists a continuous increasing function $\rt$ (the \emph{modulus of continuity} for the sections $\tau_i$) such that:
\begin{equation}\label{eq-modtau}
\forall ~ 1 \leq i \leq \nu ~ \text{and}~ w, w' \in \fT_i \quad , \quad d_{\fX}(w,w') <\rt(\e) ~ \Longrightarrow ~ \dM(\tau_i(w), \tau_i(w')) < \e ~ .
\end{equation}
\end{lemma}
\proof
Set $\ds\rt(\e) = \min \left\{\e, \min \left\{ d_{\fX}(w,w') \mid 1 \leq i \leq \nu ~ , ~ w,w' \in \fT_i ~ , ~ \dM(\tau_i(w), \tau_i(w')) \geq \e\right\}\right\}$.
\endproof

Introduce two additional constants, derived from the Lebesgue number $\eU$ chosen above.
The first is derived from a ``converse'' to the modulus function $\rp$. For each $1 \leq i \leq \nu$, consider the projection map 
$\pi_i \colon \oU_i \to \fT_i$, then introduce the constant
$$ \epsilon^{\cT}_{i} = \max \left\{\e \mid  \forall ~ x \in \oU_i ~{\rm such ~ that}  ~ D_{\fM}(x, \eU/2) \subset \oU_i ~ , ~ D_{\fX}(\pi_i(x),\e) \subset \pi_i\left( D_{\fM}(x, \eU/2)\right)\right\}.
$$ 
which measures the distance from $\eU/2$-interior points of $\oU_i$ to the exterior of their transverse projection to $\fT_i$. 
Then let 
\begin{equation}\label{eq-transdiam}
\eTU = \min \left\{\epsilon^{\cT}_{i} \mid \forall ~ 1 \leq i \leq \nu \right\}.
\end{equation}
 Note that by \eqref{eq-modtau} we have the estimate $\eTU \geq \rt(\eU/2 )$.

For $y \in \fM$ recall that $D_{\F}(y, \e)$ is the  closed ball of radius $\e$ for the leafwise metric.
Introduce a form of ``leafwise Lebesgue number'', defined by
\begin{equation}\label{eq-leafdiam}
\eFU = \min \left\{ \eFU(y) \mid ~ \forall ~ y \in \fM \right\} ~ , ~ \eFU(y) = \max \left\{ \e \mid  ~ D_{\F}(y, \e) \subset D_{\fM}(y, \eU/4)\right\}.
\end{equation}
Thus, for all $y \in \fM$, $D_{\F}(y, \eFU) \subset D_{\fM}(y, \eU/4)$.
Note that for all $r > 0$ and $z' \in D_{\F}(z, \eFU)$, the triangle inequality implies that
$D_{\fM}(z', r) \subset D_{\fM}(z, r + \eU/4)$.

\section{Holonomy of foliated spaces} \label{sec-holonomy}

The holonomy pseudogroup of a smooth foliated manifold $(M, \F)$ generalizes the induced dynamical systems    associated to a section of a flow. The holonomy pseudogroup for a matchbox manifold $(\fM, \F)$ is defined analogously, although there are delicate issues of domains which must be considered. See the articles by Haefliger \cite{Haefliger2002}, and also by Hurder \cite{Hurder2011a},  for a discussion of these and related topics. 

The properties of the holonomy pseudogroup of a matchbox manifold are fundamental to its study, as they reflect the degree to which the leaves of $\F$ are intertwined, and so the lack of a global product structure for $\fM$. This is a factor in all of the constructions later.  
 
A pair of indices $(i,j)$, $1 \leq i,j \leq \nu$, is said to be \emph{admissible} if the \emph{open} coordinate charts satisfy $U_i \cap U_j \ne \emptyset$.
For $(i,j)$ admissible, define $\fD_{i,j} = \pi_i(U_i \cap U_j) \subset \fT_i \subset \fX$.  The regularity of foliation charts imply that plaques are either disjoint, or have connected intersection. This implies that there is a well-defined homeomorphism $h_{j,i} \colon \fD_{i,j} \to \fD_{j,i}$ with domain $D(h_{j,i}) = \fD_{i,j}$ and range $R(h_{j,i}) = \fD_{j,i}$. 
Note that the map $h_{j,i}$ admits a continuous extension to a local homeomorphism on the closure of its domain, 
$\overline{h}_{j,i} \colon \overline{\fD_{i,j}} \to \overline{\fD_{j,i}}$.

The maps $\cGF^{(1)} = \{h_{j,i} \mid (i,j) ~{\rm admissible}\}$ are the transverse change of coordinates defined by the foliation charts. By definition they satisfy $h_{i,i} = Id$, $h_{i,j}^{-1} = h_{j,i}$, and if $U_i \cap U_j\cap U_k \ne \emptyset$ then $h_{k,j} \circ h_{j,i} = h_{k,i}$ on their common domain of definition. The \emph{holonomy pseudogroup} $\cGF$ of $\F$ is the topological pseudogroup modeled on $\fX$ generated by   the elements of $\cGF^{(1)}$.

A sequence $\cI = (i_0, i_1, \ldots , i_{\alpha})$ is \emph{admissible}, if each pair $(i_{\ell -1}, i_{\ell})$ is admissible for $1 \leq \ell \leq \alpha$, and the composition
\begin{equation}\label{eq-defholo}
h_{\cI} = h_{i_{\alpha}, i_{\alpha-1}} \circ \cdots \circ h_{i_1, i_0}
\end{equation}
has non-empty domain.
The domain $D(h_{\cI})$ is the \emph{maximal open subset} of $\fD_{i_0 , i_1} \subset \fT_{i_0}$ for which the compositions are defined.

Given any open subset $U \subset D(h_{\cI})$ we obtain a new element $h_{\cI} | U \in \cGF$ by restriction. Introduce
\begin{equation}\label{eq-restrictedgroupoid}
\cGF^* = \left\{ h_{\cI} | U \mid \cI ~ {\rm admissible} ~ \& ~ U \subset D(h_{\cI}) \right\} \subset \cGF ~ .
\end{equation}
The range of $g = h_{\cI} | U$ is the open set $R(g) = h_{\cI}(U) \subset \fT_{i_{\alpha}} \subset \fX$. Note that each map $g \in \cGF^*$ admits a
continuous extension to a local homeomorphism $\overline{g} \colon \overline{D(g)}  \to \overline{R(g)} \subset \fT_{i_{\alpha}}$.

The orbit of a point  $w \in \fX$ by  the action of the pseudogroup $\cGF$ is denoted by 
\begin{equation}\label{eq-orbits}
\cO(w) = \{g(w) \mid g \in \cGF^* ~, ~ w \in D(g) \} \subset \fT_* ~ .
\end{equation}

Given an admissible sequence $\cI = (i_0, i_1, \ldots , i_{\alpha})$ and any $0 \leq \ell \leq \alpha$,     the truncated sequence $\cI_{\ell} = (i_0, i_1, \ldots , i_{\ell})$ is again admissible, and we introduce the holonomy map defined by the composition of the first $\ell$ generators appearing in $h_{\cI}$, 
\begin{equation}\label{eq-pcmaps}
h_{\cI_{\ell}} = h_{i_{\ell} , i_{\ell -1}} \circ \cdots \circ h_{i_{1} , i_{0}}~.
\end{equation}
Given $\xi \in D(h_{\cI})$ we adopt the notation $\xi_{\ell} = h_{\cI_{\ell}}(\xi) \in \fT_{i_{\ell}}$. So $\xi_0 = \xi$ and
$h_{\cI}(\xi) = \xi_{\alpha}$.

\medskip

Given $\xi \in D(h_{\cI})$, let $x = x_0 = \tau_{i_0}(\xi_0) \in L_x$. Introduce the \emph{plaque chain}
\begin{equation}\label{eq-plaquechain}
\cP_{\cI}(\xi) = \{\cP_{i_0}(\xi_0), \cP_{i_1}(\xi_1), \ldots , \cP_{i_{\alpha}}(\xi_{\alpha}) \} ~ .
\end{equation}
Intuitively, a plaque chain $\cP_{\cI}(\xi)$ is a sequence of successively overlapping convex ``tiles'' in $L_0$ starting at $x_0 = \tau_{i_0}(\xi_0)$, ending at
$x_{\alpha} = \tau_{i_{\alpha}}(\xi_{\alpha})$, and with each $\cP_{i_{\ell}}(\xi_{\ell})$ ``centered'' on the point $x_{\ell} = \tau_{i_{\ell}}(\xi_{\ell})$.

Recall that $\cP_{i_{\ell}}(x_{\ell}) = \cP_{i_{\ell}}(\xi_{\ell})$, so we also adopt the notation $\cP_{\cI}(x) \equiv \cP_{\cI}(\xi)$.

 \subsection{Leafwise path holonomy}\label{subsec-lph} 
 A standard construction in foliation theory, introduced by Poincar\'{e} for sections to flows, and   developed for foliations by Reeb \cite{Reeb1952} (see also \cite{Haefliger1984}, \cite{CN1985}, \cite[Chapter 2]{CandelConlon2000}) associates to a leafwise path $\gamma$ a holonomy map $h_{\gamma}$. We describe this   construction, paying particular attention to domains and metric estimates.

Let $\cI$ be an admissible sequence. For $w \in D(h_{\cI})$, we say that $(\cI , w)$ \emph{covers} $\gamma$,
if the domain of $\gamma$ admits   a partition $0 = s_0 < s_1 < \cdots < s_{\alpha} = 1$ such that  the plaque chain $\cP_{\cI}(w) = \{\cP_{i_0}(w_0), \cP_{i_1}(w_1), \ldots , \cP_{i_{\alpha}}(w_{\alpha}) \}$ satisfies
\begin{equation}\label{eq-cover}
\gamma([s_{\ell} , s_{\ell + 1}]) \subset int (\cP_{i_{\ell}}(w_{\ell}) )~ , ~ 0 \leq \ell < \alpha, \quad \& \quad \gamma(1) \in int( \cP_{i_{\alpha}}(w_{\alpha})).
\end{equation}
It follows that $w_0 = \pi_{i_0}(\gamma(0)) \in D(h_{\cI})$.

Given two admissible sequences, $\cI = (i_0, i_1, \ldots, i_{\alpha})$ and
$\cJ = (j_0, j_1, \ldots, j_{\beta})$, such that both $(\cI, w)$ and $(\cJ, v)$ cover the leafwise path $\gamma \colon [0,1] \to \fM$, then 
$$\gamma(0) \in int( \cP_{i_0}(w_0)) \cap int( \cP_{j_0}(v_0)) \quad , \quad \gamma(1) \in int(\cP_{i_{\alpha}}(w_{\alpha})) \cap int( \cP_{j_{\beta}}(v_{\beta}) )$$ Thus both $(i_0 , j_0)$ and $(i_{\alpha} , j_{\beta})$ are admissible, and
$v_0 = h_{j_{0} , i_{0}}(w_0)$, $w_{\alpha} = h_{i_{\alpha} , j_{\beta}}(v_{\beta})$.
The proof of the following can be found in \cite{ClarkHurder2011b}.
\begin{prop}\label{prop-copc}
The maps $h_{\cI}$ and
$\ds h_{i_{\alpha} , j_{\beta}} \circ h_{\cJ} \circ h_{j_{0} , i_{0}}$
agree on their common domains. \hfill $\Box$ 
\end{prop}

\subsection{Admissible sequences}\label{subsec-admissible}
 Given a leafwise path $\gamma \colon [0,1] \to \fM$, we next construct an admissible sequence
$\cI = (i_0, i_1, \ldots, i_{\alpha})$ with $w \in D(h_{\cI})$ so that $(\cI , w)$ covers $\gamma$, and has ``uniform domains''.

Inductively, choose a partition of the interval $[0,1]$, $0 = s_0 < s_1 < \cdots < s_{\alpha} = 1$ such that for each $0 \leq \ell \leq \alpha$,
$\gamma([s_{\ell}, s_{\ell + 1}]) \subset D_{\F}(x_{\ell}, \eFU)$ where $x_{\ell} = \gamma(s_{\ell})$.
As a notational convenience, we have let
$s_{\alpha+1} = s_{\alpha}$, so that $\gamma([s_{\alpha}, s_{\alpha + 1}]) = x_{\alpha}$.
Note that we can choose $s_{\ell + 1}$ to be the largest value such that $\dF(\gamma(s_{\ell}), \gamma(t)) \leq \eFU$ for all  $s_{\ell} \leq t \leq s_{\ell +1}$. Thus, we can assume  $\alpha \leq 1 + \| \gamma \|/\eFU$.

For each $0 \leq \ell \leq \alpha$, choose an index $1 \leq i_{\ell} \leq \nu$ so that $ B_{\fM}(x_{\ell}, \eU) \subset U_{i_{\ell}}$.
Note that, for all $s_{\ell} \leq t \leq s_{\ell +1}$, $B_{\fM}(\gamma(t), \eU/2) \subset U_{i_{\ell}}$, so that
$x_{\ell+1} \in U_{i_{\ell}} \cap U_{i_{\ell +1}}$. It follows that $\cI_{\gamma} = (i_0, i_1, \ldots, i_{\alpha})$ is an admissible sequence.
Set $h_{\gamma} = h_{\cI_{\gamma}}$. Then $h_{\gamma}(w) = w'$, where $w = \pi_{i_0}(x_0)$ and
$w' = \pi_{i_{\alpha}}(x_{\alpha})$.

The construction of the admissible sequence $\cI_{\gamma}$ above has the important   property that 
  $h_{\cI_{\gamma}}$ is the composition of generators of $\cGF^*$ which have a uniform lower bound estimate $ \eTU$ on the radii of the metric balls centered at the orbit, which are contained in their domains, with $ \eTU$  independent of $\gamma$.
To see this, let  $0 \leq \ell < \alpha$, and note that $x_{\ell+1} \in D_{\F}(x_{\ell +1}, \eFU)$ implies that for some $s_{\ell} < s_{\ell + 1}' < s_{\ell + 1}$, we have that
$\gamma([s_{\ell + 1}', s_{\ell + 1}]) \subset D_{\F}(x_{\ell +1}, \eFU)$. Hence,
\begin{equation} \label{eq-unifest}
B_{\fM}(\gamma(t), \eU/2) \subset U_{i_{\ell}} \cap U_{i_{\ell +1}} ~ , ~ {\rm for ~ all} ~ s_{\ell + 1}' \leq t \leq s_{\ell + 1} ~.
\end{equation}
Then for all $s_{\ell + 1}' \leq t \leq s_{\ell + 1}$, the uniform estimate defining $\eTU > 0$ in \eqref{eq-transdiam} implies that
\begin{equation} \label{eq-domains}
B_{\fX}(\pi_{i_{\ell}}(\gamma(t)), \eTU ) \subset \fD_{i_{\ell} , i_{\ell +1}} \quad \& \quad
B_{\fX}(\pi_{i_{\ell +1}}(\gamma(t)), \eTU ) \subset \fD_{i_{\ell +1} , i_{\ell}} ~ .
\end{equation}
For the admissible sequence $\cI_{\gamma} = (i_0, i_1, \ldots, i_{\alpha})$,
recall that $x_{\ell} = \gamma(s_{\ell})$ and we set $w_{\ell} = \pi_{i_{\ell}}(x_{\ell})$.
Then by the definition \eqref{eq-defholo} of $\ds h_{\cI_{\gamma}}$ the condition \eqref{eq-domains} implies that
$D_{\fX}(w_{\ell} , \eTU) \subset D(h_{\ell})$.

\medskip

There is a converse to the above construction, which associates to an admissible sequence a leafwise path.
Let $\cI = (i_0, i_1, \ldots, i_{\alpha})$ be admissible, with corresponding holonomy map $h_{\cI}$,
and choose $w \in D(h_{\cI})$ with $x = \tau_{i_0}(w)$.

For each $1 \leq \ell \leq \alpha$, recall that
$\cI_{\ell} = (i_0, i_1, \ldots, i_{\ell})$, and let $h_{\cI_{\ell}}$ denote the corresponding holonomy map. For $\ell = 0$, let $\cI_0 = (i_0 , i_0)$.
Note that $h_{\cI_{\alpha}} = h_{\cI}$ and $h_{\cI_{0}} = Id \colon \fT_0 \to \fT_0$.

For each $0 \leq \ell \leq \alpha$, set $w_{\ell} = h_{\cI_{\ell}}(w)$ and
$x_{\ell}= \tau_{i_{\ell}}(w_{\ell})$. By assumption, for $\ell > 0$, there exists $z_{\ell} \in \cP_{\ell -1}(w_{\ell -1}) \cap \cP_{\ell}(w_{\ell})$.

Let
$\gamma_{\ell} \colon [(\ell -1)/\alpha , \ell / \alpha] \to L_{x_0}$ be the leafwise piecewise geodesic segment from $x_{\ell -1}$ to $z_{\ell}$ to $x_{\ell}$. Define the leafwise path $\gamma^x_{\cI} \colon [0,1] \to L_{x_0}$ from $x_0$ to $x_{\alpha}$ to be the concatenation of these paths.
If we then cover $\gamma^x_{\cI}$ by the charts determined by the given admissible sequence $\cI$, it follows that $h_{\cI} = h_{\gamma^x_{\cI}}$.

Thus, given an admissible sequence $\cI = (i_0, i_1, \ldots, i_{\alpha})$ and $w \in D(h_{\cI})$ with $w' = h_{\cI}(w)$, the choices above determine an initial chart $\vp_{i_0}$ with ``starting point''
$x = \tau_{i_0}(w) \in U_{i_0} \subset \fM$. Similarly, there is a terminal chart $\vp_{i_{\alpha}}$ with
``terminal point'' $x' = \tau_{i_{\alpha}}(w') \in U_{i_{\alpha}} \subset \fM$. The leafwise path $\gamma^x_{\cI}$ constructed above starts at $x$, ends at $x'$, and has image contained in the plaque chain $\cP_{\cI}(x)$.

On the other hand, if we start with a leafwise path $\gamma \colon [0,1] \to \fM$, then the initial point $x = \gamma(a)$ and the terminal point $x' = \gamma(b)$ are both well-defined. However, there need not be a unique index $j_0$ such that $x \in U_{j_0}$ and similarly for the index $j_{\beta}$ such that $x' \in U_{j_{\beta}}$. Thus, when one constructs an admissible sequence $\cJ = (j_0, \ldots , j_{\beta})$ from $\gamma$, the initial and terminal charts need not be well-defined. That is, in fact, the essence of Proposition~\ref{prop-copc}, which proved that
$$h_{\cI} | U = h_{i_{\alpha} , j_{\beta}} \circ h_{\cJ} \circ h_{j_{0} , i_{0}} | U \quad {\rm for} \quad U = D(h_{\cI}) \cap D(h_{i_{\alpha} , j_{\beta}} \circ h_{\cJ} \circ h_{j_{0} , i_{0}}) ~ .$$

We conclude this discussion with a trivial observation, and an application which yields a key technical point, that the holonomy along a path is independent of ``small deformations'' of the path.

The observation is this. Let $\cI = (i_0, i_1, \ldots, i_{\alpha})$ be admissible, with associated holonomy map $h_{\cI}$. Given $w, u \in D(h_{\cI})$, then the germs of $h_{\cI}$ at $w$ and $u$ admit a common extension, namely $h_{\cI}$. Thus, if $\gamma$, $\gamma'$ are leafwise paths defined as above from the plaque chains associated to $(\cI, w)$ and $(\cI,u)$ then the germinal holonomy maps along $\gamma$ and $\gamma'$ admit a common extension by Proposition~\ref{prop-copc}. This is the basic idea behind the following technically useful result.

\begin{lemma} \cite{ClarkHurder2011b} \label{lem-domainconst}
Let $\gamma, \gamma' \colon [0,1] \to \fM$ be leafwise paths. Suppose that $x = \gamma(0), x' = \gamma'(0) \in U_i$ and
$y = \gamma(1), y' = \gamma'(1) \in U_j$. If $\dM(\gamma(t) , \gamma'(t)) \leq \eU/4$ for all $0 \leq t \leq 1$, then the induced holonomy maps $h_{\gamma}, h_{\gamma'}$ agree on their common domain $D(h_{\gamma}) \cap D(h_{\gamma'}) \subset \fT_i$.
In particular, if curves $\gamma,\gamma'$ are sufficiently close, then they define holonomy maps which have a common extension.
\end{lemma}

\subsection{Homotopy independence}
Two leafwise paths $\gamma , \gamma' \colon [0,1] \to \fM$ are \emph{homotopic} if there exists a family of leafwise paths $\gamma_s \colon [0,1] \to \fM$ with $\gamma_0 = \gamma$ and $\gamma_1 = \gamma'$. We are most interested in the special case when $\gamma(0) = \gamma'(0) = x$ and $\gamma(1) = \gamma'(1) = y$. Then $\gamma$ and $\gamma'$ are \emph{endpoint-homotopic}
if they are homotopic with $\gamma_s(0) = x$ for all $0 \leq s \leq 1$, and similarly
$\gamma_s(1) = y$ for all $0 \leq s \leq 1$. Thus, the family of curves $\{ \gamma_s(t) \mid 0 \leq s \leq 1\}$ are all contained in a common leaf $L_{x}$. The following property then follows from an inductive application of Lemma~\ref{lem-domainconst}:

\begin{lemma}\cite{ClarkHurder2011b}\label{lem-homotopic}
Let $\gamma, \gamma' \colon [0,1] \to \fM$ be endpoint-homotopic leafwise paths. Then their holonomy maps $h_{\gamma}$ and $h_{\gamma'}$ agree on some open subset $U \subset D(h_{\gamma}) \cap D(h_{\gamma'}) \subset \fT_*$. In particular, they determine the same germinal holonomy maps. \hfill $\Box$
\end{lemma}

The following is another consequence of the strongly convex property of the plaques:
\begin{lemma}\cite{ClarkHurder2011b} \label{lem-close}
Suppose that $\gamma, \gamma' \colon [0,1] \to \fM$ are leafwise paths for which $\gamma(0) = \gamma'(0) = x$ and $\gamma(1) = \gamma'(1) = x'$, and suppose that $\dM(\gamma(t), \gamma'(t)) < \eU/2$ for all $a \leq t \leq b$.
Then $\gamma, \gamma' \colon [0,1] \to \fM$ are endpoint-homotopic. \hfill $\Box$
\end{lemma}

Given $g \in \cGF^*$ and $w \in D(g)$, let $[g]_w$ denote the germ of the map $g$ at $w \in \fT_*$. Set
\begin{equation}\label{eq-holodef}
\G_{\F}^w = \{ [g]_w \mid g \in \cGF^* ~ , ~ w \in D(g) ~ , ~ g(w) = w\} ~ .
\end{equation}
Given $x \in U_i$ with $w = \pi_i(x) \in \fT_*$, the elements of $\G_{\F}^w$ form a group, and by Lemma~\ref{lem-homotopic} there is a well-defined homomorphism $h_{\F,x} \colon \pi_1(L_x , x) \to \G_{\F}^w$ which is called the \emph{holonomy group} of $\F$ at $x$.

\subsection{Non-trivial holonomy} 
Note that if $y \in L_x$ then the homomorphism
$h_{\F , y}$ is conjugate (by an element of $\cGF^*$) to the homomorphism $h_{\F , x}$.
A leaf $L$ is said to have \emph{non-trivial germinal holonomy} if for some $x \in L$, the homomorphism $h_{\F , x}$ is non-trivial. If the homomorphism $h_{\F , x}$ is trivial, then we say that $L_x$ is a \emph{leaf without holonomy}. This property depends only on $L$, and not the basepoint $x \in L$.
The foliated space $\fM$ is said to be \emph{without holonomy} if for every $x \in M$, the leaf $L_x$ is without germinal holonomy.

\begin{lemma}\cite{ClarkHurder2011b} \label{lem-noholo}
Let $\fM$ be a foliated space, and $L_x$ a leaf without holonomy. Fix a regular covering for $\fM$ as above, 
and let $w \in \fT_*$ be the local projection of a point in $L_x$. Given  plaques chains $\cI$, $\cJ$ such that $w \in Dom(h_{\cI}) \cap Dom(h_{\cJ})$ with $h_{\cI}(w) = w' = h_{\cJ}(w)$, then $h_{\cI}$ and $h_{\cJ}$ have the same germinal holonomy at $w$. Thus, for each $w' \in \cO(w)$ in the $\cGF^*$ orbit of $w$,
there is a well-defined holonomy germ $h_{w,w'}$.
\end{lemma}
\proof
The composition $g = h_{\cJ}^{-1} \circ h_{\cI}$ satisfies $g(w) = w$, so by assumption there is some open neighborhood $w \in U$ for which $g | U$ is the trivial map. That is, $ h_{\cI} | U = h_{\cJ} | U$.
\endproof

Finally, we recall a basic result of Epstein, Millet and Tischler \cite{EMT1977} for foliated manifolds, whose proof applies verbatim in the case of foliated spaces.
\begin{thm} \label{thm-emt}
The union of all leaves without holonomy in a foliated space $\fM$ is a dense $G_{\delta}$ subset of $\fM$. In particular, there exists at least one leaf without germinal holonomy. \hfill $\Box$
\end{thm}

\section{Dynamics of matchbox manifolds}\label{sec-mme}

Many of the   concepts of dynamical systems for flows (and more generally group actions) on a compact manifold admit generalizations to the foliation dynamics associated to a matchbox manifold, by considering the leaves of $\F$ in place of the orbits of the action. See \cite{Hurder2011a} for a discussion of this topic. 
We  first recall several important classical definitions from topological dynamics, adapted to the case of matchbox manifolds, and  several results concerning their dynamical properties from \cite{ClarkHurder2011b}. 

\begin{defn} \label{def-equicontinuous}
The holonomy pseudogroup $\cGF$ of $\F$ is \emph{equicontinuous} if for all $\epsilon > 0$, there exists $\delta > 0$ such that for all $g \in \cGF^*$, if $w, w' \in D(g)$ and $d_{\fX}(w,w') < \delta$, then $d_{\fX}(g(w), g(w')) < \epsilon$.
\end{defn}

\begin{defn} \label{def-expansive}
The holonomy pseudogroup $\cGF$ of $\F$ is \emph{expansive}, or more properly $\e$-expansive, if there exists $\e > 0$ such that for all $w, w' \in \fT_*$, there exists $g \in \cGF$ with $w, w' \in D(g)$ such that $d_{\fX}(g(w), g(w')) \geq \e$.
\end{defn}

Equicontinuity for $\cGF$ gives \emph{uniform} control over the domains of arbitrary compositions of generators.

\begin{prop}\cite{ClarkHurder2011b}\label{prop-uniformdom}
Assume the holonomy pseudogroup $\cGF$ of $\F$ is equicontinuous. Then there exists $\dTU > 0$ such that for every leafwise path
$\gamma \colon [0,1] \to \fM$, there is a corresponding admissible sequence $\cI_{\gamma} = (i_0 , i_1 , \ldots , i_{\alpha})$ so that
$B_{\fX}(w_0, \dTU) \subset D(h_{\cI_{\gamma}})$, where $x = \gamma(0)$ and $w_0 = \pi_{i_0}(x)$.

Moreover, for all $0 < \e_1 \leq \eTU$ there exists $0 < \delta_1 \leq \dTU$ independent of the path $\gamma$, such that
$h_{\cI_{\gamma}}(D_{\fX}(w_0, \delta_1)) \subset D_{\fX}(w', \e_1)$ where $w' = \pi_{i_{\alpha}}(\gamma(1))$.

Thus, $\cGF^*$ is equicontinuous as a family of local group actions.
\end{prop}

  Recall that a foliated space $\fM$ is \emph{minimal} if each leaf $L \subset \fM$ is dense.
The following is an immediate consequence of the definitions. 
\begin{lemma} \label{lem-minimal}
The  foliated space $\fM$ is minimal if and only if  $\cO(w)$   is dense in $\fT_*$ for all   $w \in \fT_*$. 
\end{lemma}

The following result is, at first glance, very surprising. It has been previously shown for flows \cite{AHO1991} and $\mR^n$-actions \cite{Clark2002}. The proof of it is given in detail in \cite[\S 4.1]{ClarkHurder2011b}, and fundamentally uses the conclusions of Proposition~\ref{prop-uniformdom}.

\begin{thm}\cite{ClarkHurder2011b}\label{thm-minimal}
If $\fM$ is an equicontinuous matchbox manifold, then $\fM$ is minimal.
\end{thm}

It is well-known that an equicontinuous action of a countable group $\G$ on a Cantor set $\mK$ admits a finite $\G$-invariant decomposition of $\mK$ into clopen subsets with arbitrarily small diameter. It is shown  in \cite[section 6]{ClarkHurder2011b} that a corresponding result holds for an equicontinuous pseudogroup action.

\begin{thm}\cite{ClarkHurder2011b}\label{thm-invariants}
Let $\fM$ is an equicontinuous matchbox manifold, and $w_0 \in \fT_*$ a basepoint.
Then there exists a descending chain of clopen subsets
$$\cdots \subset V_{\ell +1} \subset V_{\ell} \subset \cdots V_0 \subset \fT_*$$
such that for all $\ell \geq 0$, 
$w_0 \in V_{\ell}$ and ${\rm diam}_{\fX}(V_{\ell}) < \dTU/2^{\ell}$.

Moreover, each $V_{\ell}$ is $\cGF$-invariant in the following sense:
  if $\gamma$ is a path with initial point $\gamma(0) \in V_{\ell}$, then the holonomy map $h_{\gamma}$ satisfies $V_{\ell} \subset Dom(h_{\gamma})$, and if
$h_{\gamma}(V_{\ell}) \cap V_{\ell} \ne \emptyset$, then $h_{\gamma}(V_{\ell}) = V_{\ell}$.
\hfill $\Box$
\end{thm}

It follows that the collection $\ds \left\{ h_{\gamma}(V_{\ell}) \mid \gamma(0) \in V_{\ell} \right\}$ of subsets of the transverse space $\fT_*$ forms a finite clopen partition, and these sets   are permuted by the   action of the holonomy pseudogroup.

\medskip

If the action of the pseudogroup $\cGF$ is   expansive, then the domains of arbitrary compositions of generators for its holonomy typically do not admit uniform estimates as in Proposition~\ref{prop-uniformdom}. However, the compactness of $\fM$ implies there is a uniform  estimate on the size of the domain of a holonomy map  formed from a bounded    number of compositions used to define it.
Recall that the path length in the $\dF$ metric of the piecewise $C^1$-curve $\gamma(t)$ is denoted by   $\| \gamma \|$.

\begin{prop}\label{prop-domest} For each $\epsilon > 0$ and $r > 0$, there exists
$0 < \delta(\epsilon, r) \leq \epsilon$
so that for any piecewise smooth leafwise path $\gamma \colon [0,1] \to \fM$ with $\| \gamma \| \leq r$, then there exists
an admissible sequence $\cI = (i_0, i_1, \ldots, i_{\alpha})$ such that $(\cI , w)$ covers $\gamma$ with:

\begin{enumerate}
\item $w_0 = \pi_{i_{0}}(\gamma(0)) \in D(h_{\cI})$ and $D_{\fX}(w_0 , \delta(\epsilon, r)) \subset D(h_{\cI})$;
\item $h_{\cI}(D_{\fX}(w_0, \delta(\epsilon, r))) \subset D_{\fX}(w', \e)$ where $w' = \pi_{i_{\alpha}}(\gamma(b))$.
\end{enumerate}
\end{prop}
\proof
By the arguments of section~\ref{sec-holonomy}, there exists an admissible sequence $\cI = (i_0, i_1, \ldots, i_{\alpha})$
with $w \in D(h_{\cI})$ and $\alpha \leq 1 + \|\gamma\|/\eFU \leq 1 + r/\eFU$ where $\eFU> 0$ is defined by \eqref{eq-leafdiam}.

We estimate the size of the domain $D(h_{\cI})$.
For each $0 \leq \ell \leq \alpha$, we have
$\gamma([s_{\ell}, s_{\ell + 1}]) \subset D_{\F}(x_{\ell}, \eFU)$ where $x_{\ell} = \gamma(s_{\ell})$. Moreover, for the associated admissible sequence $\cI_{\gamma} = (i_0, i_1, \ldots, i_{\alpha})$, we have that
for all $0 \leq t \leq 1$, $B_{\fM}(\gamma(t), \frac{1}{2} \eU) \subset U_{i_{\ell}}$.

The proof will be by downward induction. Let $w_{\ell} = \pi_{i_{\ell}}(x_{\ell})$ and
set $\cI_{\ell} = (i_0, i_1, \ldots, i_{\ell})$ with corresponding holonomy map $h_{\cI_{\ell}}$.
Then $h_{\cI_{\ell}}(w_0) = w_{\ell}$. Let $h_{\ell} = h_{i_{\ell +1} , i_{\ell}}$ so that $h_{\ell} \circ h_{\cI_{\ell}} = h_{\cI_{\ell +1}}$.

For every admissible pair $(i,j)$ the holonomy homeomorphism $h_{j,i}$ is uniformly continuous as it has compact domain. Since there is only a finite number of distinct non-empty intersections $U_i \cap U_j$, for every $\epsilon > 0$ there exists a $0< \delta_\epsilon \leq \epsilon$ such that for every admissible pair $(i,j)$ if $w,w' \in D(h_{j,i})$ and $d_\fX(w,w') < \delta_{\epsilon}$ then $d_\fX \left( h_{j,i}(w),h_{j,i}(w') \right) < \epsilon$.

Recall $\eTU> 0$ as defined by \eqref{eq-transdiam}. Given $\epsilon > 0$, set $\epsilon_\alpha = \min \{\eTU/2 , \epsilon\}$. Now proceed by downward induction.
For $0<\ell \leq \alpha$ assume that $\epsilon_{\ell}$ has been defined. Then
denote $\delta_{\ell} = \delta_{\epsilon_{\ell}}$ and $\epsilon_{\ell-1} = \delta_{\ell}$ as defined using equicontinuity as above.

By the choice of the covering of the admissible sequence $\cI$, we have $D_{\fX}(w_{\ell}, \eTU) \subset D(h_{\ell})$ and so $D_{\fX}(w_{\ell}, \delta_{\ell}) \subset D(h_{\ell})$, and by the choice of $\delta_{\ell}$ we have $D_{\fX}(w_{\ell}, \delta_{\ell}) \subset D(h_{\ell} \circ h_{\ell+1} \circ \cdots \circ h_{\alpha})$. Then $\delta(\epsilon, r) = \delta_1$ satisfies the required conditions.
\endproof
 
\section{Foliated microbundles and Reeb Structure Theorem}\label{sec-microbundles}

The strategy of our construction of the transverse foliation $\cH$ begins with the choice of a compact, path-connected region $K \subset L$ of a leaf $L$, such that the holonomy along closed paths in $K$ is trivial. In this section, we discuss how this assumption implies there is a ``thickening'' of $K$ to a foliated neighborhood of $K$ in $\fM$. In later sections, we establish that for such a thickening, if it is sufficiently ``thin'', then there exists a transverse Cantor foliation defined on it. 

The existence of a foliated open neighborhood  of a compact subset $K$ follows from the analog of the  \emph{Reeb Stability Theorem}, which is one of the fundamental results of foliation theory for smooth manifolds \cite{CN1985,CandelConlon2000,Reeb1952,Tamura1992}.  The most general version of these ideas is formulated  in terms of the ``foliated microbundle'' associated to the holonomy covering of a leaf in a foliated space. (See   Milnor \cite{Milnor2009} for a discussion of the concept of foliated microbundles for manifolds). This general formulation admits a generalization to matchbox manifolds, as described below.

\subsection{Non-trivial holonomy}
Recall that we assume there is a fixed regular covering $\cU$ for $\fM$, as in Proposition~\ref{prop-regular} which we can assume to be finite.
Assume there is given a fixed transversal base-point  $w_0 \in int(\fT_1)$ and let $L_0$ be the leaf through $x_0 = \tau_1(w_0) \in U_1$. 
  Let $h_{\F, x_0} \colon \pi_1(L_0 , x_0) \to \GF^{w_0}$ denote the transverse holonomy representation, where \eqref{eq-holodef} defines the group $\GF^{w_0}$ of homotopy classes of closed paths based at $x_0$. 
  
  Since the map $h_{\F, x_0}$ is a homomorphism, its kernel $G_0 \subset \pi_1(L_{x_0} , x_0)$ is a normal subgroup, and we let   $\Pi \colon \wtL_{0} \to L_0$ denote the normal (holonomy) covering associated  to $G_0$. The leafwise Riemannian metric $\dF$ on   $L_0$ lifts to a Riemannian metric $\wtdF$ on $\wtL_0$ such that $\Pi$ is a local isometry.

Choose $\wtx_0 \in \wtL_0$ such that $\pi(\wtx_0) = x_0$. By definition, given any closed path $\wtgamma \colon [0,1] \to \wtL_0$ with basepoint $\wtx_0 = \wtgamma(0) = \wtgamma(1)$, the image of $\wtgamma$ in $L_0$ has trivial germinal holonomy as a leafwise path in $\fM$.
Then the transverse holonomy map defined by a path $\wtgamma$ in $\wtL_0$ starting at $\wtx_0$ is determined by the endpoint $\wtgamma(1)$.

We next select a collection of points in $\wtL_0$ which are sufficiently dense, so that homotopy classes of paths between the points capture all of the holonomy defined by the leaf $L_0$.
\begin{defn} \label{def-net}
Let $(X, d_X)$ be a complete separable metric space. Given
$0 < e_1 < e_2$, a subset $\cM \subset X$ is a \emph{$(e_1 , e_2)$-net}, or \emph{Delaunay set},  if:
\begin{enumerate}
\item $\cM$ is $e_1$-separated: for all $y \ne z \in \cM$, $e_1 \leq d_{X}(y,z)$;
\item $\cM$ is $e_2$-dense: for all $x \in X$, there exists some $z \in \cM$ such that $d_{X}(x,z) \leq e_2$.
\end{enumerate}
\end{defn}
 Next,  we  construct a $(e_1 , e_2)$-net in the given leaf $L_0$ in   accordance with the constants defined in section~\ref{sec-concepts}.
Recall that $\eFU$ defined by \eqref{eq-leafdiam} was chosen so that every leafwise disk of radius $\eFU$ is contained in a metric ball of $\fM$ of radius $\eU/4$. That is, for all $y \in \fM$, $D_{\F}(y, \eFU) \subset D_{\fM}(y, \eU/4)$.
Let $e_2 = \eFU/4$, then choose $\cM_0 \subset L_0$ an $(e_1 , e_2)$-net for $L_0$
for some $0 < e_1 < e_2$.
We can assume without loss of generality that $x_0 \in \cM_0$. Condition~(\ref{def-net}.2) implies that the collection of leafwise open disks
$\ds \{B_{\F}(z , \eFU/2) \mid z \in \cM_0 \}$ is an open covering of $L_0$.

\subsection{Delaunay sets and covers of $\fM$} 
 For each $z \in \cM_0$, choose an index $1 \leq i_z \leq \nu$ so that $B_{\fM}(z, \eU) \subset U_{i_z}$. Without loss of generality, we can assume that $B_{\fM}(x_0, \eU) \subset U_{1}$.
Then note that for all $z' \in D_{\F}(z, \eFU)$, we have $z' \in D_{\fM}(z, \eU/4)$ so
the triangle inequality implies that
\begin{equation}\label{eq-containment}
D_{\F}(z', \eFU) \subset D_{\fM}(z', \eU/4) \subset D_{\fM}(z, \eU/2) \subset B_{\fM}(z, \eU) \subset U_{i_z} .
\end{equation}

\begin{lemma}\label{lem-newcover}
If the leaf $L_0$ is dense, then the collection $\{U_{i_z} \mid z \in \cM_0\}$ is a covering of $\fM$  with Lebesgue number $\eU/2$.
\end{lemma}
\proof
Let $y \in \fM$, then $L_0$ is dense so there exists $y' \in L_0$ with $\dM(y,y') < \eU/4$.
Let $z \in \cM_0$ with $\dF(y', z) \leq e_2 = \eFU/4$.
Then $y' \in D_{\fM}(z, \eU/4)$ by \eqref{eq-leafdiam},
hence $y \in B_{\fM}(z, \eU/2) \subset U_{i_z}$ by \eqref{eq-containment}.

Next, we show that $\eU/2$ is a Lebesgue number for this covering. Let 
 $y'' \in B_{\fM}(y, \eU/2)$ then the above implies that  $y'' \in B_{\fM}(z, \eU) \subset U_{i_z}$ by the choice of $i_z$ above. Thus,   
$B_{\fM}(y, \eU/2) \subset U_{i_z}$.
\endproof

\subsection{Foliated microbundle}\label{subsec-folmicro}
We construct the foliated microbundle associated to the choices made above. 
Choose $\wtx_0 \in \wtcM_0$ with $\Pi(\wtx_0) = x_0 \in \cM_0$, and let 
 $\wtcM_0 = \Pi^{-1}(\cM_0)$ which is a $(e_1 , e_2)$-net for $\wtL_0$ for the Riemannian metric lifted from $L_0$.
The points of $\wtcM_0$ are denoted by $\wtz$,    where $\wtz$ is a lift of $z \in \cM_0$.
The idea is to associate to each $\wtz$,  a disjoint copy of the foliation chart $\oU_{i_z}$ then form the union for all $\wtz \in \wtcM_0$ with appropriate identifications. 

For each $\wtz \in \wtcM_0$,  set $\wtU_{\wtz} = \oU_{i_z} \times \{\wtz\}$.
For $(x, \wtz) \in \wtU_{\wtz}$ define $\Pi \colon \wtU_{\wtz} \to \oU_{i_z}$ by $\Pi(x, \wtz) = x$. 

The  leafwise plaques for $\wtU_{\wtz}$  are defined by $\wtcP_{\wtz}(\wty) = \cP_{i_z}(x) \times \{\wtz\}$   
for  $\wty = (x, \wtz) \in \wtU_{\wtz}$.

In the case where $x \in \cP_{i_z}(z)$,  then   $\wtcP_{\wtz}(\wty)$ is identified with the plaque of $\wtL_0$ containing $\wtz$, so the collection
$\ds \{\wtcP_{\wtz}(\wtz) \mid \wtz \in \cM_0 \}$ are open sets in  $\wtL_0$.
In fact,  as   $D_{\wtL_0}(\wtz , \eFU) \subset \wtcP_{\wtz}(\wtz)$ for each $\wtz \in \wtcM_0$, this is an open covering of $\wtL_0$. 
One thinks of the plaques $\wtcP_{\wtz}(\wtz)$ as ``convex tiles'', and the collection $\{ \wtcP_{\wtz}(\wtz) \mid \wtz \in \wtcM_0\}$ as a ``tiling'' of $\wtL_0$. The interiors of the plaques need not be disjoint, so this is not a proper tiling in the usual sense (for example see \cite{AP1998, BG2003}, or \cite[\S 11.3.C]{CandelConlon2000}).

\begin{defn}\label{fol-micro}
The \emph{foliated microbundle} over $\wtL_0$ is the space
\begin{equation}\label{eq-folmicro}
\wtfN_0 = \bigcup_{\wtz \in \wtcM_0} ~ \wtU_{\wtz} ~ {\Big \slash} \sim
\end{equation}
where $\wty \in \wtU_{\wtz}$ and $\wty' \in \wtU_{\wtz'}$ are identified
if $\Pi(\wty) = \Pi(\wty')$ and $\wtcP_{\wtz}(\wtz) \cap \wtcP_{\wtz'}(\wtz') \ne \emptyset$. The connected components of $\wtfN_0$ form the leaves of a foliation $\wtF$.
\end{defn}

Informally, the space $\ds \wtfN_0 $ is simply the union of copies of all flow boxes associated as above to  the points  $\wtz \in \wtcM_0$ and identified in the obvious way to obtain an ``open normal neighborhood'' of $\wtL_0$. 

For each $\wtz \in \wtcM_0$, 
the composition
$\ds \wtvp_{\wtz} \equiv \vp_{i_z} \circ \Pi \colon \wtU_{\wtz} \to [-1,1]^n \times \fT_{\wtz}$
defines a coordinate chart on $\wtfN_0$, making it into a foliated space with foliation denoted by $\wtF$.
Let $\wtpi_{\wtz} \colon \wtU_{\wtz} \to \fT_{\wtz}$ be the normal coordinate,
and $\wtlambda_{\wtz} \colon \wtU_{\wtz} \to [-1,1]^n$ be the leafwise coordinate.

 The   foliated microbundle   $\ds \wtfN_0 $ provides a uniform setting for  the holonomy maps of paths in the leaf $\wtL_0$.  
 Introduce the transversals to $\wtF$ which are the lifts of the transversals to $\F$, where 
for each $\wtz \in \wtcM_0$, let $\fT_{\wtz} = \fT_{i_z} \times \{\wtz\}$. 
Given $\wtz \in \wtcM_0$, subset $V \subset \fT_{\wtz}$ and $\xi \in [-1,1]^n$, we obtain a local section for $\wtF$ by
\begin{equation}
\wttau_{\wtz , \xi} \colon V \to \wtU_{\wtz} ~ , ~ \wttau_{\wtz, \xi}(w) = \wtvp_{\wtz}^{-1}(\xi, w) = (\vp_{i_z} ^{-1}(\xi, w) , \wtz).
\end{equation}

Note that while the core leaf $\wtL_0$ of the foliated microbundle    $\wtfN_0$ is a  regular covering of the leaf $L_0 \subset \fM$, the projection of other leaves $\wtL$ of $\wtF$ may not be coverings, as the leaves of $\wtF$ may ``escape'' from the flow boxes defining $\wtfN_0$.    
If the groupoid $\cGF^*$  has equicontinuous dynamics, then  with a suitable choice of transversals,  the leaves of $\wtF$ in $\wtfN_0$ are all coverings of leaves of $\F$ in $\fM$. We show this, but require a technical aside.

A path $\wtgamma \colon [0,1] \to \wtL_0$ is said to be \emph{nice}, if there exists a partition
$a = s_0 < s_1 < \cdots < s_{\alpha} = b$ such that for each $0 \leq \ell \leq \alpha$, the restriction
$\wtgamma \colon [s_{\ell}, s_{\ell + 1}] \to \wtL_0$ is a geodesic segment between points
$\wtz_{\ell} = \wtgamma(s_{\ell}), \wtz_{\ell +1} = \wtgamma(s_{\ell +1})\in \wtcM_0$ with $\dF(\wtz_{\ell}, \wtz_{\ell +1}) < \eFU$.
Then $\widetilde{\cI} = ({\wtz_{0}}, \ldots , {\wtz_{\alpha}})$ is an admissible sequence for $\wtF$,  and $\cI = (i_{\wtz_{0}}, \ldots , i_{\wtz_{\alpha}})$ is an admissible sequence for $\F$. The sequence $\widetilde{\cI}$ defines the holonomy maps $\wth_{\widetilde{\cI}}$ for $\wtF$, and $\cI$ defines the holonomy map $h_{\cI}$ for $\F$.
Clearly, $\wth_{\cI}$ is just the lift of $h_{\cI}$, and $h_{\cI}$ is the holonomy map for the leafwise path $\gamma = \Pi \circ \wtgamma$ constructed in section~\ref{sec-holonomy}.
As before, we note that $\wth_{\widetilde{\cI}}$ depends only on the endpoints of $\cI$. For $\wtz \in \wtcM_0$ let $\wth_{\wtz}$ denote the holonomy along some nice path $\wtgamma_{\wtz}$ from $\wtx_0$ to $\wtz$, considered as a transformation of the space $\wtfT$, which is the disjoint union of the local transversals $\fT_{\wtz}$.
Let $h_{\wtz}$ denote the holonomy along the   path  $\gamma_{\wtz} = \Pi \circ \wtgamma_{\wtz}$.

\subsection{Equicontinuous matchbox manifolds and Thomas tubes}\label{subsec-equicontinuous}
If $\fM$ is an equicontinuous matchbox manifold, then for any $\e > 0$, let $V_{\ell} \subset \fT_1$ be a clopen set with $w_0 \in V_{\ell}$ and 
${\rm diam}_{\fX}(V_{\ell}) < \e$ given by Theorem~\ref{thm-invariants}.

For $h_{\gamma} \in \cGF^*$ with $V_{\ell} \subset D(h_{\gamma})$, set $V_{\ell}^{\gamma} = h_{\gamma}(V_{\ell})$. Then either  $V_{\ell}^{\gamma} = V_{\ell}$ or $V_{\ell}^{\gamma} \cap V_{\ell} = \emptyset$, and the set 
$\ds \{V_{\ell}^{\gamma} \mid V_{\ell} \subset D(h_{\gamma})~ , ~ \gamma \in \cGF^*\}$ is a clopen partition of $\fT_*$. 

For $\wtz \in \wtcM_0$ there is a nice path $\gamma_{\wtz}$ from $\wtz_0$ to $\wtz$ which defines a holonomy map denoted by
$h_{\wtz} \equiv h_{\gamma_{\wtz}}$. Then for $\wtz \in \wtcM_0$ define
\begin{equation}\label{def-basicblocks}
V_{\ell, \wtz} = h_{\wtz}(V_{\ell}) \subset \fT_{i_{\wtz}} ~ , ~ \wtV_{\ell, \wtz} = \wth_{\wtz}(V_{\ell}) = V_{\ell, \wtz} \times \{\wtz\} \subset \fT_{{\wtz}} ~ .
\end{equation}
The union of the sets $V_{\ell, \wtz}$ is the saturation of $V_{\ell}$ under the action of the pseudogroup $\cGF^*$, and hence it forms a clopen partition of $\fT_*$.
Introduce the local coordinate chart saturations of these sets:
\begin{equation}\label{eq-coordinatecovers}
\fU^{V}_{\ell, \wtz} = \pi_{i_{\wtz}}^{-1}(V_{\ell, \wtz}) \subset \oU_{i_{\wtz}} ~, ~ \wtfU^{V}_{\ell, \wtz} = \fU^{V}_{\ell, \wtz} \times \{\wtz\} \subset \wtU_{\ell, \wtz} ~ .
\end{equation}
Then $\fU^{V}_{\ell, \wtz} $ is the union of the plaques in $\oU_{i_{\wtz}}$ through the points of $V_{\ell, \wtz}$.

\begin{defn} \label{def-thomas}
The \emph{Thomas tube} associated with $V_{\ell}$ is the subset of the microbundle $\wtfN_0$,
\begin{equation} \label{eq-thomas}
\wtfN(V_{\ell}) ~ = ~ \bigcup_{\wtz \in \wtcM_0} ~ \wtfU^{V}_{\ell, \wtz} ~ \subset ~ \wtfN_0 .
\end{equation}
\end{defn}
This construction generalizes that used by Thomas for equicontinuous flows in \cite{EThomas1973}, hence the name.

The image $\Pi(\wtfN(V_{\ell})) \subset \fM$ is the saturation by $\F$ of the clopen set $V_{\ell}$, hence $\Pi(\wtfN(V_{\ell})) = \fM$. Note that each leaf $\wtL$ of $\wtF$ in $\wtfN(V_{\ell})$ has no holonomy and is properly embedded by construction, though the projection $L$  in $\fM$  of $\wtL$ is recurrent, as  $\F$ is minimal by Theorem~\ref{thm-minimal}. 

In terms of shape theory, the above shows that an equicontinuous matchbox manifold admits a shape approximation of diameter less than $\e$ which is the image of the  foliated space $\ds \wtfN(V_{\ell})$.

\subsection{Reeb neighborhoods of compact sets} 
For the general case, where $\cGF^*$ is not assumed to be equicontinuous, we require a modification of the above construction.

Let $L_0 \subset \fM$ be a leaf,  and  $K \subset L_0$   a proper base with $z_0 \in K$, as in as in Definition~\ref{def-LBB}. Then by assumption,   $K$ is   a union of closed plaques in the  foliation charts $\{U_{i} \mid 1 \leq i \leq \nu\}$ and  there exists   $\wtK \subset \wtL_0$ which is a connected compact subset  of the holonomy covering $\Pi \colon \wtL_0 \to L_0$, such that  $K = \Pi(\wtK)$.  Choose a basepoint $\wtz_0 \in  \wtK \cap \wtcM_0$ with  $z_0 = \Pi(\wtz_0)$.
Note that $\wtK \cap \wtcM_0$ is  finite.

Suppose we are given a clopen neighborhood $z_0 \in V_{z_0} \subset \fT_*$ then for each $\wtz \in  \wtcM_0$ there is a well-defined holonomy map 
$h_{\wtz} \equiv h_{\gamma_{\wtz}}$ where $\gamma_{\wtz}$ is a path in $\wtK$ from $\wtz_0$ to $\wtz$. As above, if $V_{z_0} \subset D(h_{\wtz})$ then we define $\ds V_{\wtz} = h_{\wtz}(V_{z_0})$.

For the general case, there is no expectation that the sets $\ds \{V_{\wtz} \mid \wtz \in \wtK \cap \wtcM_0\}$ have any special relationship to each other. Instead, we formulate a condition ensuring they are defined and disjoint.

\begin{defn}\label{def-admissibledisjoint}
A clopen neighborhood $z_0 \in V_{z_0} \subset \fT_*$ is \emph{$\wtK$-admissible} if $\ds V_{z_0} \subset D(h_{\wtz})$ for each $\ds \wtz \in \wtK \cap \wtcM_0$,  and \emph{$\wtK$-disjoint} if, in addition,  $V_{\wtz} \cap V_{z_0} = \emptyset$ for all $\ds \wtz_0 \ne \wtz \in   \wtK \cap \wtcM_0$.
\end{defn}

 For a $\wtK$-admissible   clopen set with $z_0 \in V_{z_0} \subset \fT_*$ we can define the \emph{Reeb neighborhood} of $\wtK$      by
\begin{equation}\label{eq-tesselspace}
  \wtfN(\wtK, V_{z_0}) ~ \equiv ~ \bigcup_{\wtz \in \wtK \cap \wtcM_0} ~ \wtpi_{\wtz}^{-1}(\wtV_{\wtz}) ~ \subset ~ \wtfN_0 .
\end{equation}
If  $V_{z_0}$ is  also \emph{$\wtK$-disjoint} then we define a \emph{Reeb neighborhood}  of $K$ in $\fM$ by 
\begin{equation}\label{eq-tesselspace2}
  \fN(K, V_{z_0}) ~ \equiv ~ \Pi\left(  \wtfN(\wtK, V_{z_0}) \right) = \bigcup_{\wtz \in \wtK \cap \wtcM_0} ~ \Pi\left\{ \wtpi_{\wtz}^{-1}(\wtV_{\wtz}) \right\} ~ \subset ~ \fM .
\end{equation}

Note that each leaf of the restricted foliation $\wtF | \wtfN(\wtK, V_{z_0})$ is a properly embedded compact subset, and the holonomy $\wth_{\wtgamma}$ along any closed loop $\wtgamma$ contained in $\wtfN(\wtK, V_{z_0})$ is trivial. 
If $V_{z_0}$ is $\wtK$-disjoint,
then the same holds for each path component of $\fN(K, V_{z_0})$.

\section{Transverse Cantor foliations}\label{sec-cantor}

 For a smooth foliation $\F$ of a compact manifold $M$, a foliation $\cH$ on $M$ is transverse to $\F$ if the leaves of $\F$ and $\cH$ have complementary dimensions, and are everywhere transverse  as submanifolds of $M$. If $\F$ is a smooth foliation of codimension-one, then the normal distribution to the tangent bundle of $\F$ is always integrable, and the foliation $\cH$ this distribution defines is obviously transverse to $\F$. When $\F$ has codimension greater than one, then the integrability of the normal distribution is not guaranteed, and the existence of a transverse foliation $\cH$ is a strong assumption. For example, if $\F$ is a foliation constructed to be transverse to the fibers of a fibration $\pi \colon M \to B$, then the fibers of $\pi$ define a transverse foliation $\cH$. 
 
 A \emph{Cantor foliation} $\cH$ on a matchbox manifold $\fM$ is a ``foliation'' whose ``leaves'' are Cantor sets. This notion is problematic, as there is no inherent way to speak of the ``regularity'' of the leaves, analogous to the case of smooth foliations where the leaves are defined as integral manifolds for a distribution. It is, in fact, preferable to think of the leaves of such a foliation as defining an equivalence relation $\cong_{\cH}$ on $\fM$, where two points $x,y \in \fM$ satisfy $x \cong_{\cH} y$ if and only if they belong to the same leaf of $\cH$.
 
 For example, in the case of the attractors for Axiom A diffeomorphisms of a compact manifold $M$ studied by Williams   \cite{Williams1967,Williams1974}, the closed attractor $\Omega$ of such a map $f \colon M \to M$ has a local product structure satisfying the conditions to be a matchbox manifold. The differential $Df \colon TM \to TM$ is assumed to restrict to a hyperbolic map on an open neighborhood  $\Omega \subset U \subset M$, where the leaves of $\F$ in $\Omega$ are   unstable manifolds of the map, and the intersections of the stable leaves for $Df | TU$ with $\Omega$ are Cantor sets, and define the transverse Cantor foliation $\cH$. In this case, the regularity of $\cH$ is derived from the dynamical properties of the smooth map $f$.
 
The works by Putnam \cite{Putnam1996,PS1999} study the  Smale spaces  introduced by Ruelle \cite{Ruelle1988}, and the thesis of  Wieler \cite{Wieler2012a,Wieler2012b} studies a   generalization of these ideas to spaces  which are   matchbox manifolds.  They define     the ``leaves'' of a Cantor foliation  $\cH$   dynamically, as the stable or unstable manifolds for a ``locally hyperbolic'' action. Again, the the regularity of $\cH$ is derived from the dynamical properties of the   map $f$ and their assumptions on the local properties of this map.

For the general case of matchbox manifolds, there is no dynamical system associated to a locally hyperbolic map of the space available to define $\cH$. Rather, the approach we take begins with the   regularity inherent in a covering of $\fM$ by coordinate charts for $\F$.

 \subsection{Cantor foliations}
   Recall that we assume there is a fixed regular covering $\{U_{i} \mid 1 \leq i \leq \nu\}$ of $\fM$ by foliation charts, as in Proposition~\ref{prop-regular}, with charts $\vp_i \colon \oU_i \to [-1,1]^n \times \fT_i$ where $\fT_i \subset \fX$ is a clopen subset. Moreover,  there exists  a foliated extension $\whvarp_i \colon \whU_i \to (-2,2)^n \times \fT_i$ where $\oU_i \subset \whU_i \subset \fM$ is an open neighborhood of $\oU_i$ and $\whvarp_i | \oU_i = \vp_i$. 
   
On each chart $\oU_{i}$ there is a   Cantor foliation $\cH_i \equiv \cH | \oU_i$ whose leaves are the closed sets $\vp_i^{-1}(\{\xi\} \times \fT_i)$ for $\xi \in [-1,1]^n$. The problem is that for $\oU_i \cap \oU_j \ne \emptyset$ for $1 \leq i \ne j \leq \nu$, the two foliations $\cH_i$ and $\cH_j$ need not agree. The natural question is then, do there exists small perturbations $\wtcH_i$  of the   foliations $\cH_i$ on coordinate charts so that $\wtcH_i = \wtcH_j$ on all non-empty overlaps? If the leaves of $\F$ are defined by the action of $\mR^n$ on $\fM$, then this existence question can be formulated as a problem of solving a cocycle equation over the covering of $\fM$ with values in $\mR^n$.

 In the more general case of matchbox manifolds, the problem of existence of $\cH$ can be considered as asking for a solution of a  ``non-linear cocycle equation''.  
This is what we give in the later sections of this paper, essentially by showing that for sufficiently small domains, there is a sufficiently good linear approximation to the problem, and this can be solved by an explicit recursive procedure. Motivated by these considerations, we next define a transverse Cantor foliation $\cH$ to $\F$ as a perturbation of the   solutions $\cH_i$ in coordinate charts, so that they define a globally defined ``foliation'' or equivalence relation. The definition is actually given for a closed subset $\fB \subset \fM$, as our recursive procedure  will construct the solution over an increasing sequence of such subspaces.

 For a clopen set $V \subset \fT_i $ set
 $\ds  \fU^{V}_{i} = \pi_{i}^{-1}(V) \subset \oU_{i}$ and $\ds  \whfU^{V}_{i} = \whvarp_i^{-1}(V) \subset \whU_{i}$.

\begin{defn}\label{def-cantorfol}
Let $\fM$ be a matchbox manifold, and  $\fB \subset \fM$ a closed subset. An equivalence relation $\approx_{\cH}$ on $\fB$ 
is said to define a \emph{transverse Cantor foliation} $\cH$ of $\fB$ if for each 
 $x \in \fB$, the class $\cH_x = \{y \in \fB \mid y \approx_{\cH} x\}$ is a Cantor set. Moreover, we require that 
  there exists a covering of $\fM$ by foliation charts as above, such that for each $x \in \fB$,    there exists:
\begin{enumerate}
\item   $1 \leq i_x \leq \nu$ with $x \in U_{i_x}$, 

\item  a clopen subset $V_x \subset \fT_{i_x}$ with $w_x = \pi_{i_x}(x) \in V_x$;

\item    a homeomorphism into $\Phi_{x} \colon [-1,1]^n \times V_{x} \to \whU_{i_x}$ such that 
$$\Phi_{x}(\xi , w_x) = \whvarp_{i_x}^{-1}(\xi , w_x)  ~ {\rm for} ~   \xi \in [-1,1]^n ,$$
\item   for   $\xi \in [-1,1]^n$ and  $z = \whvarp^{-1}(\xi , w_x) \in \fB$, the image $\Phi_x(\{ \xi \}  \times V_{x}) = \cH_z \cap \whfU^{V_x}_{i_x}$.
\end{enumerate}
The leaves of the ``foliation'' $\cH$ are defined to be the equivalence classes $\cH_x$ of ~ $\approx_{\cH}$ in $\fB$.
\end{defn}

Condition~\ref{def-cantorfol}.1 specifies a coordinate chart covering a neighborhood of   $x \in \fB$, while 
Condition~\ref{def-cantorfol}.2 specifies the transverse projection of the leaf   $\cH_x$.
 Condition~\ref{def-cantorfol}.3 states that the reparametrization map $\Phi_{x}$ agrees with the coordinate chart $\whvarp_{i_x}$   on the horizontal slice though $x$. 
  Then  Condition~\ref{def-cantorfol}.4  is the main assertion, that the image of $V_x$ under $\Phi_x$ equals the leaf $\cH_x$ in the chart $\whfU^{V_x}_{i_x}$.

Note that the   images of the maps $\Phi_x$ are allowed to take values in the open neighborhood $\whU_{i_x}$ of $\ds \oU_{i_x}$ as the  ``leaf''  $\cH_z$ may not have constant horizontal coordinate $\lambda_{i_x}$ so that  for   a leafwise boundary point  $z \in \fB \cap \oU_{i_x}$, the equivalence class $\cH_z$ need not be   contained in $\oU_{i_x}$.

   The functions $\Phi_x$ are the adjustments to the local vertical foliations $\cH_{i_x}$   so that the perturbed leaves   are coordinate independent, hence are well defined on $\fB$. In this sense, they can be viewed as solutions of the non-linear cocycle problem mentioned above. 
The maps  $\Phi_x$ are not required to be leafwise smooth, and in fact, our solutions will be   piecewise-linear maps when restricted to leaves.
 
\bigskip

\subsection{Cantor foliations and microbundles}
The definition of transverse Cantor foliation applies equally to foliated microbundles introduced in section~\ref{sec-microbundles}. 
Given   a leaf $L_0 \subset \fM$ with holonomy covering $\Pi \colon \wtL_0 \to L_0$,   an  
   $(e_1 , e_2)$-net    $\cM_0$   for $L_0$, and   the lifted net   $\wtcM_0$   on   $\wtL_0$, 
     form the foliated microbundle  $\wtfN_0$.

  Let $K \subset L_0$ be a proper base. 
 Chose a basepoint $z_0 \in K \cap \cM_0$, with lift $\wtz_0 \in \wtK \cap \wtL_0$.
Let   $V_{z_0} \subset \fT_0$ be a $\wtK$-admissible clopen subset so that  $\wtfN(K, V_{z_0}) \subset \wtfN_0$  is defined by \eqref{eq-tesselspace}.
Assume that there is a transverse Cantor foliation for a closed subset $\fB$ with $\ds \wtfN(K, V_{z_0})  \subset \fB \subset \wtfN_0$.
 Condition~\ref{def-cantorfol}.2 is satisfied by  the collection of translates 
$\{ V_{\wtz} \mid \wtz \in  \wtK \cap \wtcM_0\}$, and we obtain:

 \begin{lemma} \label{lem-bifoliated}
 The set $\ds \fB$ is  a bi-foliated neighborhood of $\wtK$ in the foliated microbundle $\wtfN_0$ for which there exists a bi-foliated homeomorphism 
$\ds \wtPhi \colon \fB   \to  \wtK \times V_{z_0}$. \hfill $\Box$
 \end{lemma}
 \proof
Given $\wtx \in \wtK$ and $\wty \in \wtcH_{\wtx}$ the leaf $\wtL_{\wty}$ intersects the transversal $\wtcH_{\wtz_0} \cong V_{z_0}$ in a unique point $w_{\wty}$. We then set 
   $\ds \wtPhi(\wty) =   (\wtx  , w_{\wty}) \in \wtK \times V_{z_0}$. 
 \endproof

 \begin{defn} \label{def-holoequiv}
 A transverse Cantor foliation $\wtcH$ for a compact subset $\wtfB \subset \wtfN_0$ is said to be \emph{holonomy equivariant} if it defines a  transverse Cantor foliation $\cH$ on the image $\fB = \Pi(\wtfB) \subset \fM$.
 \end{defn}
 If the restriction $\Pi \colon \wtfB \to \fM$ is injective, then this is always the case. For example, if $V_{z_0}$ is $\wtK$-disjoint then any transverse Cantor foliation $\wtcH$ on $\ds \fU(\wtK , V_{z_0}) $ will be holonomy equivariant. On the other hand, if $V_{z_0}$ is $\wtK$-admissible but not $\wtK$-disjoint, then  $\wtcH$ is holonomy equivariant if the images of the equivalence classes $\Pi(\wtcH_{\wty})$ in $\fM$ agree on the overlap of any two coordinate charts on $\ds \fU(\wtK , V_{z_0})$. 
   This condition will be satisfied, for example, if the equivalence classes $\wtcH_{\wtx}$ for $\wtx \in \wtK$ are defined in terms of a transverse Cantor foliation $\cH$ on the image $\ds \Pi ( \fU(\wtK , V_{z_0})  ) \subset \fM$.

 \begin{defn} \label{def-complete}
 A pair $\{ \wtK, V_{z_0} \}$ is a \emph{complete model} for $\fM$,  if 
  $V_{z_0}$ is $\wtK$-admissible, and the map $\ds \Pi \colon   \wtfN(\wtK , V_{z_0})  \to \fM$ is surjective.
 \end{defn}
 
For example, we have:

\begin{lemma}\label{lem-equicdense}
Let $\fM$ be an equicontinuous matchbox manifold. Then there exists $\{ \wtK, V_{z_0} \}$ which  is a   complete model  for $\fM$.
\end{lemma}
\proof
Given $\e > 0$, choose  a $\cGF$-invariant clopen subset $\ds w_0 \in V_{\ell} \subset \fT_*$ as in Definition~\ref{def-thomas} with associated Thomas tube $\ds \wtfN(V_{\ell})$. As previously remarked, $\ds \Pi(\wtfN(V_{\ell})) = \fM$ so is a complete model.
\endproof

 The basic observation is that if $\{ \wtK, V_{z_0} \}$ is a \emph{complete model} for $\fM$, and $\wtcH$ is a holonomy equivariant transverse Cantor foliation for $\wtfB = \wtfN(\wtK , V_{z_0})$, then $\wtcH$ defines a transverse Cantor foliation on $\fM$.

\section{Voronoi tessellations} \label{sec-VT}

The concept of a Voronoi cell decomposition (or tessellation) of  Euclidean space   is extraordinarily useful for applications of geometry to a variety of problems, and is   very well-studied (for example,  see \cite[Introduction]{OBKC2000}). Associated to every tiling of $\mR^n$ is a Voronoi tesselation, and conversely a tesselation yields a tiling. 
In the next sections,   we   develop the basic concepts of Voronoi tessellations and Delaunay triangulations in a form applicable to  metric spaces derived from the leaves of  matchbox manifolds.   

Assume there is a fixed regular covering $\{\vp_i \colon \oU_i \to [-1,1]^n \times \fT_i \mid 1 \leq i \leq \nu\}$ of $\fM$ by foliation charts, where each $\fT_i \subset \fX$ is a clopen subset, as in   Proposition~\ref{prop-regular}.

Let  $L \subset \fM$ a leaf with induced leafwise Riemannian metric $d_L$. 
 Let $X \subset L$ be a closed connected set which is a union of plaques. The typical  examples we consider are for  $X = L$, or for $X$ a compact subset of $L$ which contains a proper base $K$ in its interior.

Recall that  $\lF > 0$ is the leafwise constant defined in Lemma~\ref{lem-stronglyconvex},   such that for all $x \in L$, the closed disk $D_{L}(x, \lF) \subset L$ is   strongly convex. 

Let $\cN_X$ be a given  $(d_1,d_2)$-net for $X$.   The value of the density constant $d_2$ will be fixed later, with $d_2 \leq \eFU/5 <   \eFU/4 = e_2$ and  depending on estimates derived from the geometry of the leaves and the metric distortion of the transverse holonomy maps.
Associated to the  net $\cN_X$ is the  \emph{Voronoi tessellation}, which is a partition of the space into compact star-like regions, called cells.

\subsection{Voronoi cells}

Introduce the ``leafwise nearest--neighbor distance'' function, where for $y \in  L$,
\begin{equation}\label{eq-kF2}
\kX(y) = \inf \left\{ d_{L}(x, y) \mid x \in \cN_X \right\}.
\end{equation}

Note that $\kX(y) = 0$ if and only if $y \in \cN_X$.

\begin{defn}
For $x \in \cN_X$, define its \emph{Dirichlet region}, or \emph{Voronoi cell},  in $L$ by
\begin{equation}\label{eq-voronoi2}
\cC_L(x) = \left\{ y \in L \mid d_{L}(x,y) = \kX(y)\right\}.
\end{equation}
\end{defn}
That is, for $x \in \cN_X$ the Voronoi cell $\cC_L(x)$ consists of the points $y \in L$ which are closer to $x$ in the leafwise metric than to any other point of $\cN_X$. Thus, for each $y \in L$ there exists some $x \in \cN_X$ with $y \in \cC_L(x)$. In particular, for $\cC_X(x) =  \cC_L(x) \cap X$, then the collection of closed subsets 
\begin{equation}
\left\{\cC_X(x)  \mid x \in \cN_X\right\}
\end{equation}
forms a closed covering of $X$, and we obtain the Voronoi decomposition of $X$
\begin{equation}\label{eq-Vdecomposition}
X ~ = ~ \bigcup_{x \in \cN_X}  ~ \cC_X(x)
\end{equation}
 Introduce the subset of $\cN_X$ consisting of net points whose Voronoi cells lie in $X$, 
 \begin{equation}\label{eq-interiornet}
\cN_X^* = \{ x \in \cN_X \mid \cC_L(x) \subset X \}
\end{equation}
We develop some of the   properties of the cells $\cC_L(x)$ for $x \in \cN_X$ and $\cN_X^*$.
In particular, Lemma~\ref{lem-netdensity} below and the assumptions that $X$ is a union of plaques  and that $d_2 \leq e_2$ implies $\cN_X^* $ is not empty.

\begin{lemma}\label{lem-celldiam2}
For each $x \in \cN_X^*$, 
\begin{equation}\label{eq-Vcellbounds}
D_{L}(x, d_1/2)   \subset \cC_L(x) \subset D_{L}(x, d_2)  
\end{equation}
 In particular, $\cC_L(x)$ has diameter at most $2d_2$. \hfill $\Box$
\end{lemma}

The upper bound estimate in \eqref{eq-Vcellbounds} need not hold if $\cC_L(x) \not\subset X$, for then the set $\cN_X$ is not   $d_2$-dense in all of $L$. However, we always have:

\begin{lemma}\label{lem-celldiam3}
For    $x  \in \cN_X$,   $\cC_X(x) \subset D_L(x,d_2)$.  \hfill $\Box$
\end{lemma}

A set $Y \subset L$ is \emph{star-like  with respect to $x \in Y$ }   if for all $y \in Y$, each geodesic ray  from $x$ to $y$ is contained in $Y$.

\begin{lemma} \label{lem-starlike2} 
For each $x \in \cN_X$, the set  $Y = \cC_L(x) \cap D_L(x,\lF)$ is star-like with respect to $x$. In particular, for all  $x \in \cN_X^*$ the set $\cC_L(x)$ is star-like with respect to $x$. \hfill $\Box$
\end{lemma}

The strong convexity of disks $D_L(x,\lF)$ also yields the following.
 \begin{lemma}\label{lem-netdensity}
 If $x \in \cN_X$ and there exists   $r > d_2$ for which $B_L(x,r) \subset X$, then $x \in \cN_X^*$. 
\end{lemma}
\proof
Suppose that  $y \in \cC_L(x)$ but $y \not\in X$. 
Let $\sigma_{x,y} \colon [0,1] \to   L$  be a geodesic segment with $\sigma_{x,y}(0) = x$ and $\sigma_{x,y}(1) = y$, and the length equal to $d_L(x,y)$.
Let $0 < s  < 1$ be the greatest value for which  $y' = \sigma_{x,y}(s) \in \cC_X(x)$, then $d_2(x,y') \geq r > d_2$ by assumption.
As  $y' \in X$, there exists $z \in \cN_X$ with    $d_L(y',z) < d_2$. 
Then 
$$d_L(y,z) \leq d_L(y,y') + d_L(y',z) < d_L(y,y') + d_2 <  d_L(y,y') + r \leq d_L(y,y') + d_L(y',x) = d_L(y,x)$$
which contradicts that $y \in \cC_L(x)$. Thus, $\cC_L(x) \subset X$ hence $x \in \cN_X^*$.
\endproof

Next, we introduce the \emph{star-neighborhoods} of Voronoi cells.  Given $x \in \cN_X$, introduce the \emph{vertex-sets}
\begin{equation}\label{eq-vertexset}
\cV_X(x) = \{ y \in \cN_X \mid \cC_X(y) \cap \cC_X(x) \ne \emptyset \} ~ ; \quad \cV_X^*(x) = \{y \in \cV_X(x) \mid y \ne x\}.
\end{equation}
Note that $\cV_X(x)$ is a finite set by the net condition on $\cN_X$, and $y \in \cV_X^*(x)$ if and only if $x \in \cV_X^*(y)$.

\begin{defn}\label{lem-star-nbhd}
For $x \in \cN_X$ the ``star-neighborhood'' of the Voronoi cell $\cC_X(x)$ is the set
\begin{equation}\label{eq-starset}
\cS_X(x) = \bigcup_{y \in \cV_X(x)} ~ \cC_X(y).
\end{equation}
\end{defn}

\begin{lemma}\label{lem-star2}
Assume that $d_2 \leq \lF/5$. 
For each $x \in \cN_X$, $\cS_X(x) \subset B_{L}(x, 3d_2) \subset B_{L}(x, \lF)$, hence $\cS_X(x)$ is contained in a strongly convex subset of $L$.
\end{lemma}
\proof
Suppose that $\cC_X(x) \cap \cC_X(y) \ne \emptyset$. As $\cC_X(z)$ has diameter at most $2d_2$ for all $z \in \cN_X$ by Lemma~\ref{lem-celldiam3}, 
we obtain $\cS_X(x) \subset B_{L}(x, 3d_2)$. As $d_2 \leq \lF/5$,   the claim follows.
\endproof

For $y \in \cV_X^*(x)$ set
\begin{equation}\label{eq-hyperplanes2}
H(x,y) = \{ z \in D_{L}(x, \lF) \mid d_{L}(x,z) \leq d_{L}(y,z)\}.
\end{equation}
Thus $H(x,y)$ contains the set of points in the closed disk $D_{L}(x, \lF)$ which are closer to $x$ than to $y$.

Clearly, each $H(x,y)$ is closed, and the strong convexity of $D_{L}(x, \lF)$ implies that for $x \in \cN_X^*$, 
\begin{equation}\label{eq-cellshyperplanes}
\cC_L(x) = \bigcap_{y \in \cV_X^*(x)} ~ H(x,y) .
\end{equation}
Conversely, $\cC_L(x) \cap \cC_L(y) \ne \emptyset$  implies that the intersection
\begin{equation}\label{eq-hyperplanes3}
L(x,y) = H(x,y) ~ \cap ~ H(y,x) ~ \ne ~ \emptyset .
\end{equation}
For example, if $L$ is isometric to Euclidean space $\mR^2$, then   $H(x,y)$ is the intersection of  the disk $D_{L}(x, \lF)$ with the half-plane   in $\mR^2$  consisting of the points which are closer to $x$ than to $y$. Thus,  $L(x,y)$ is a line segment. In the more general case, where $L$ is   a complete Riemannian manifold, then the local picture of $L(x,y)$ is similar to   the Euclidean case, as seen below. However, unless $L$ has a global convexity property,   $L$ may have focal points  and  the global structure of $L(x,y)$ is not so easily described. Thus, we   restrict consideration  to convex neighborhoods. 

\begin{lemma}\label{lem-hyperplanes}
For $x \in \cN_X$ and $y \in \cV_X^*(x)$,   $L(x,y) \cap D_{L}(x, \lF)$ is a codimension-one closed submanifold.
\end{lemma}
\proof
We have $x,y \in D_{L}(x, \lF)$ by Lemma~\ref{lem-star2}. As the metric $d_{L}$ is strongly convex when restricted to $D_{L}(x, \lF)$,  the functions $f_x(z) = d_{L}(x,z)^2$ and $f_y(z) = d_{L}(y,z)^2$ are both regular on $D_{L}(x, \lF)$, which implies that $L(x,y)  \cap D_{L}(x, \lF)$ is a codimension-one closed submanifold.
\endproof

Now restrict attention to $x \in \cN_X^*$ so that $\cC_L(x) \subset X$. Let $\cV_X^1(x) \subset \cV_X^*(x)$ be the subset corresponding to the codimension-one  faces of the boundary of $\cC_L(x)$.
That is, $y \in \cV_X^1 (x)$ if and only if  
$ \partial_y \cC_X(x) = \cC_L(x) \cap L(x,y)$ has non-trivial interior  as a subset of the submanifold $L(x,y)$.
Then the topological boundary $\partial \cC_X(x)$ is the finite union
\begin{equation}\label{eq-boundary2}
\partial \cC_X(x) = \bigcup_{y \in \cV_X^1(x)} ~ \partial_y \cC_X(x).
\end{equation}
We summarize the results of this section.
\begin{prop}\label{prop-voronoi1}
Let $\cN_X$ be an $(d_1,d_2)$-net in  $X$, such that $d_2 \leq \lF/5$. Then there exists a subset $\cN_X^* \subset \cN_X$ 
and a collection of closed sets $\ds   \{\cC_X(y) \mid y \in \cN_X \}$
satisfying:
\begin{enumerate}
\item $\cC_X(x) \subset X$ for each $x \in \cN_X$;
\item  $\cC_X(x) \subset D_{L}(x, d_2)$ for each $x \in \cN_X$;
\item  $int(\cC_X(x)) \cap int (\cC_X(y)) = \emptyset$ for each pair   $x \ne y \in \cN_X$; 
\item The collection $\{\cC_X(y) \mid y \in \cN_X \}$ is a closed covering of $X$.
\end{enumerate}
In addition, for $x \in \cN_X^*$ we have:
\begin{enumerate}\setcounter{enumi}{4}
\item $\cC_X(x) = \cC_L(x)$;
\item  $D_{L}(x, d_1/2) \subset \cC_X(x)$;
\item  $\cC_X(x)$ is star-like with respect to $x$;
\item $\partial \cC_X(x)$ is a union of codimension-one submanifolds with boundary.
\end{enumerate}

\end{prop}

The collection $\{\cC_X(y) \mid y \in \cN_X \}$ is called the \emph{Voronoi tessellation} of $X$ associated to $\cN_X$.

\section{Delaunay simplicial complex} \label{sec-DSC}

We next introduce the \emph{Delaunay simplicial complex}   obtained from a  $(d_1,d_2)$-net $\cN_X$ for $X \subset L$, using the ``circumscribed sphere'' characterization of the simplices. Recall that  $d_2 \leq \eFU/5 <   \eFU/4 = e_2$. 

Let $0 < r < \lF$. Then the leafwise sphere of radius $r$ centered at $z$ is
$$S_{L}(z,r) \equiv \{y \in L \mid d_{L}(z,y) = r\} = D_{L}(z,r) - B_{L}(z,r).$$
Note that if $B_{L}(x,r) \cap \cN_X = \emptyset$ for $x \in X$, then $r < d_2$ by the definition of $d_2$.

 \subsection{Definition of a simplicial complex}

The \emph{Delaunay complex} $\Delta(\cN_X)$ of $L$ derived from the net $\cN_X$ is defined by specifying the subsets of $\cN_X$ which form the vertices of the simplices in $\Delta(\cN_X)$. For $k \geq 0$, denote by $\Delta^{(k)}(\cN_X)$ the collection of $k$-simplices, defined as follows:

\begin{defn}\label{def-simplex} For each $z_0 \in \cN_X$, the set $\Delta(z_0) = \{z_0\}$ is a
$0$-simplex in $\Delta^{(0)}(\cN_X)$.

For $k > 0$, a $(k+1)$-tuple $\{z_0, \ldots , z_k\} \subset \cN_X$ forms a $k$-simplex $\Delta(z_0, \ldots , z_k) \in \Delta^{(k)}(\cN_X)$ if there exists $x \in L$ and $0 < r \leq d_2$ such that
$B_{L}(x,r) \cap \cN_X = \emptyset$, and $\{z_0, \ldots , z_k\} \subset S_{L}(z,r) \cap \cN_X$. Then $ S_{L}(x,r)$ is called the \emph{circumscribed sphere} of the simplex $\{z_0, \ldots , z_k\}$.
\end{defn}

If $\Delta(z_0, \ldots , z_k) \in \Delta^{(k)}(\cN_X)$, then every subset of $(\ell +1)$-points,
$\{z_{i_0}, \ldots , z_{i_{\ell}}\} \subset \{z_0, \ldots , z_k\}$ 
yields an $\ell$-simplex $\Delta(z_{i_0}, \ldots , z_{i_{\ell}}) \in \Delta^{(\ell)}(\cN_X)$, as the circumscribed sphere condition holds for all subsets. In particular, we have well-defined face and boundary operators defined on $\Delta(\cN_X)$.

\subsection{Realization of a Delaunay simplex}

If the manifold $L$ is Euclidean, then given a $k$-simplex $\Delta(z_0, \ldots , z_k) \in \Delta^{(k)}(\cN_X)$,
the convex hull of the  vertices defines a geometric $k$-simplex in $L$, which is its \emph{geometric realization}. 
For a non-Euclidean manifold, this elementary and intuitive approach need not work, as the  convex span of a 
$k$-simplex need not be a $k$-dimensional subset  if the leaves have curvature.
 Rather, one must choose a procedure for ``filling in'' the geometric simplex spanned by a set of vertices, in order to obtain a geometric realization. 

For a $1$-simplex $\Delta(z_0, z_1)$, there is a canonical ``filling in''   using the geodesic between $z_1$ and $z_0$, which is unique due to the strong convexity of $B_{L}(z_0,\lF)$. For    higher-dimensional simplices, we  use an inductive procedure  to   fill in the faces using the geodesic cone from each successive vertex.

Define the standard $k$-simplex $\Delta^k$ in $\mR^{k+1}$ by the barycentric coordinate approach,
$$
\Delta^k ~ = ~ \left\{ (t_0, \ldots , t_k) \mid t_{\ell} \geq 0 ~ , ~ t_0 + \cdots + t_k = 1 \right\}.
$$
The vertices of $\Delta^k$ are the coordinate vectors $\vec{e}_{\ell} = (0, \ldots , 1, \ldots , 0)$ where the unique non-zero entry   is   in the $(\ell +1)$-coordinate position.

\begin{lemma}\label{lem-geofill}
Let $\Delta(z_0, \ldots , z_k) \in \Delta^{(k)}(\cN_X)$ be given, so that $\{z_0, \ldots , z_k\} \subset B_L(z_0, \lF)$. Then there exists a diffeomorphism
$\sigma_k \colon \Delta^k \to L$ such that $\sigma_k(\vec{e}_{\ell}) = z_{\ell}$, and the maps $\{ \sigma_i \mid 0  \leq i \leq k\}$ are natural with respect to the face operators.
\end{lemma}
\proof
The map $\sigma_k$ is defined by induction on the dimensions of the faces of $\Delta^k$. Set $\sigma_k(\vec{e}_{\ell}) = z_{\ell}$.

Given a string $I = i_0 < i_1 < \cdots < i_{\nu}$ with $0 \leq i_0$ and $i_{\nu} \leq k$, define the $I$-face $\partial_I \Delta^k$ to be the subset consisting of points where the only non-zero entries are in the coordinates appearing in the string. For $\nu > 0$, let $I' = i_0 < i_1 < \cdots < i_{\nu - 1}$. By induction, we may assume that the map $\sigma_k \colon \partial_{I'} \Delta^k \to L$ has been defined.

Note that each point $\vec{v} \in \partial_I \Delta^k$ can be written  
$\vec{v} = (1-s) \cdot \vec{v}' + s \cdot \vec{e}_{i_{\nu}}$ where $ \vec{v}' \in \partial_{I'} \Delta^k $ and $0 \leq s \leq 1$.
The point $z' = \sigma_k(\vec{v}') \in L$ is defined by the inductive hypothesis, and so there exists a unique geodesic segment
$\tau \colon [0,1] \to B_L(z_0, \lF)$ such that $\tau(0) = z'$ and $\tau(1) =  z_{i_{\nu}}$. Then set $\sigma_k(\vec{v}) = \tau(s)$. The resulting map defined on  $\Delta^k$ satisfies the conclusions of   Lemma~\ref{lem-geofill}.
\endproof

\begin{defn}\label{def-realization}
Let $\Delta(z_0, \ldots , z_k) \in \Delta^{(k)}(\cN_X)$, then the \emph{geometric realization} is the set
\begin{equation}
\left| \Delta(z_0, \ldots , z_k) \right| = \sigma_k(\Delta^k) 
\end{equation}
\end{defn}

\medskip

\begin{lemma} For all $0 \leq \ell \leq k$, we have $\ds \left| \Delta(z_0, \ldots , z_k) \right| \subset B_L(z_{\ell}, \lF)$. 

\end{lemma}
\proof
Let $x \in L$ and $0 < r \leq d_2$ such that $\{z_0, \ldots , z_k\} \subset S_{L}(z,r) \cap \cN_X$. Thus, $d_L(z_{\ell}, z_{\ell'}) \leq 2d_2$ for all $0 \leq \ell' \leq k$, and so the set of vertices $\{z_0, \ldots , z_k\} \subset D(z_{\ell}, 2d_2) \subset B_{L}(z_{\ell}, \lF)$.
  As $\ds B_{L}(z_{\ell}, \lF)$ is strongly convex, 
the geodesic segment between any two vertices of $\Delta(z_0, \ldots , z_k)$ is also contained in $\ds B_{L}(z_{\ell}, \lF)$. Then proceed inductively, following the construction of $\sigma_k$  in the proof of Lemma~\ref{lem-geofill}, and it follows that the image of the map $\sigma_k$ is also contained in $\ds B_{L}(z_{\ell}, \lF)$.
\endproof

\begin{remark}\label{rem-nonuniquedelone}
{\rm
As was already mentioned, if the manifold $L$ is not flat, the map $\sigma_k$ may depend on the ordering of the set of vertices  $\{z_0, \ldots , z_k\}$   for $k>1$, except on the edges of a simplex $\Delta^k$. Indeed, the ordering of vertices in the string $I = i_0 < i_1 < \cdots < i_{\nu}$ defines a choice of geodesic spray from the vertex with the largest index $i_{\nu}$ to the vertices $i_{\ell}$ with $\ell < \nu$.
 As 
 the points of $\sigma_k(\partial_I \Delta^k)$ are obtained by flowing along geodesic curves, a different ordering of vertices defines a different choice of spanning geodesic rays. In case when $L$ is a surface, this simply results in different parametrizations of the set $\ds \left| \Delta(z_0, z_1 , z_2) \right|$, as the boundary is $1$-dimensional and so well-defined.    If $L$ is not flat and has dimension $n > 2$, the image of a point $\vec{v} \in \Delta^k$ need not be the same under the maps defined by these choices. In our applications, there will be given a ``local ordering'' of the points in $\cN_X$, which defines an ordering of the set of vertices of a given simplex, so that the geometric realization 
 $\left| \Delta(z_0, \ldots , z_k) \right| $  is thus well-defined. }
\end{remark}

The Voronoi cell decomposition and Delaunay triangulation of $L$ are closely related. For Euclidean space, one says that $\Delta(\cN_X)$ is dual to the Voronoi tessellation. For the general case of a Riemannian manifold with bounded geometry,    we have the following   results.

\begin{prop}\label{prop-duality} For $z_0 \in \cN_X^*$,
let $\{z_1, \ldots , z_k\} \subset \cV_X^1(z_0)$. Then
\begin{equation}\label{eq-duality}
L(z_0, z_1) \cap \cdots \cap L(z_0, z_k) \cap \cC_L(z_0) \ne \emptyset \quad \Longleftrightarrow \quad \Delta(z_0, \ldots , z_k) \in \Delta^{(k)}(\cN_X) .
\end{equation}
\end{prop}
\proof
Recall   that for $z  \in \cN_X^*$,
 $\cC_L(z) \subset D_L(z, \lF)$. 
Recall also that for $z  \in \cN_X^*$ and $y \in \cV_X^*(z)$,   $L(z,y) \cap D_{L}(z, \lF)$ is a smooth submanifold formed by the intersecting boundaries of the Voronoi cells $\cC_L(z)$ and $\cC_L(y)$, and  
 $\cV_X^1(z) \subset \cV_X^*(z)$ is the subset corresponding to the codimension-one  faces of the boundary of $\cC_X(z)$. 

Let $x \in L(z_0,z_1) \cap \cdots \cap L(z_0,z_k) \cap \cC_L(z_0)$ and set $r = d_{L}(x,z_0)$.
Then $d_{L}(x,z_i) = d_{L}(x,z_0) = r$ for each $1 \leq i \leq k$, and thus  $\{z_0, \ldots , z_k\} \subset S_L(x,r)$. 
As each $z_i \in \cV_X^1(z_0)$, we have $d_L(z_0,z_i) \leq 2 d_2$ and hence   $r \leq d_2$. 
By the definition of the Voronoi cells, $B_{L}(x,r) \cap \cN_X = \emptyset$. Suppose not, then   there exists $y \in B_{L}(x,r) \cap \cN_X$ with  $d_L(y, x) < r = d_L(z_i,x)$ for  $0 \leq i \leq k$,   and so  $x \not\in \cC_L(z_0) \subset X$. 
This implies $\Delta(z_0, \ldots , z_k) \in \Delta^{(k)}(\cN_X)$.

Conversely, for $\{z_1, \ldots , z_k\} \subset \cV_X^1(z_0)$ with  $\{z_0, \ldots , z_k\} \subset S_L(x,r)$, then $x$ is equidistant from each point $z_i$ and so $x \in L(z_0,z_j)$ for all $1 \leq i \leq k$.  Moreover, for all $0 \leq i \leq k$, $x \in \cC_L(z_i)$ as     $B_{L}(x,r) \cap \cN_X = \emptyset$ implies there is no $z \in \cN_X$ with $d_L(x,z) < d_L(x,z_i)$. 

Thus, $\ds L(z_0, z_1) \cap \cdots \cap L(z_0, z_k) \cap \cC_L(z_i) \ne \emptyset$.
\endproof

A point $x \in \partial \cC_L(z_0)$ is called \emph{extremal} if the distance function $y \mapsto d_L(z_0,y)$ has a local maximum on $\cC_L(z_0)$ at $y = x$. 
Let $z_0 \in \cN_X^*$, so that $\cC_X(z_0) = \cC_L(z_0)$. 
For $z_i \in \cV_X^1(z_0)$, the boundary component
$ \partial_{z_i} \cC_L(z_0) = \cC_L(z_0) \cap L(z_0, z_i)$ has codimension one.  Thus, for $z_0 \in \cN_X^*$, a point
$x \in \partial \cC_L(z_0)$ is extremal exactly when there is $\{z_1, \ldots , z_n\} \subset \cV_X^1(z_0)$ with
\begin{align*}
x = \omega(z_0, \ldots , z_n) ~ = ~ L(z_0, z_1) \cap \cdots \cap L(z_0, z_n) \cap \cC_L(z_0),
\end{align*}
and $\omega(z_0, \ldots , z_n) $ is the center of a circumscribed sphere containing $\{z_0, \ldots , z_n\}$ with radius
\begin{align*}
r(z_0, \ldots , z_n) ~ = ~ \dF(z_{\ell}, \omega(z_0, \ldots , z_n)) ~ , ~ 0 \leq \ell \leq n.
\end{align*}

\medskip

Now introduce the \emph{simplicial cone} of $z_0 \in \cN_X^*$
\begin{equation}\label{eq-simplicialcone}
\cC_{\Delta}(z) = \bigcup \left \{ \left| \Delta(z_0, \ldots , z_n) \right| ~ \mid \Delta(z_0, z_1, \ldots , z_n) \in \Delta^{(n)}(\cN_X) \right\} \subset B(z, \lF)
\end{equation}

\begin{prop}\label{prop-filling}
For all $z \in \cN_X^*$,  $ \cC_L(x) \subset \cC_{\Delta}(z)$.
\end{prop}
\proof
Let $\{x_1, \ldots , x_k\} \subset \partial \cC_L(z)$ denote the set of extremal points for the   distance function $d_L(z,y)$. 

For each $1 \leq i \leq k$, let 
$\ds \Delta(z, z_1^i, \ldots , z_n^i) \in \Delta^{(n)}(\cN_X)$ denote the $n$-simplex defined by the center $x_i$, so that 
$\ds \{z_1^i, \ldots , z_n^i\} \subset \cV_X^1(z)$. The claim is that the intersection
$|\Delta(z, z_1^i, \ldots , z_n^i)| \cap \partial \cC_L(z)$ is a topological ball   with ``center''  $x_i$, and the union of all these boundary regions for $1 \leq i \leq k$ is a closed covering of $\partial \cC_L(z)$. Thus, 
$$\cC_L(z) \subset \bigcup_{i=1}^k ~ |\Delta(z, z_1^i, \ldots , z_n^i)| $$
Each vertex $z_{\ell}^i$ corresponds to a face   of $\partial \cC_X(z)$  as in \eqref{eq-boundary2}  and an associated hyperplane,   denoted by 
 $$ \partial_{\ell}^i \cC_X(z) = \cC_L(z) ~ \cap ~ L(z,z_{\ell}^i)$$
 The geodesic segment from $z$ to $z_{\ell}^i$ intersects the face $\ds \partial_{\ell}^i \cC_X(z)$ in an interior point, denoted by $\whz_{\ell}^i$. This geodesic segment is a boundary $1$-simplex of each $n$-simplex that intersects this face. These $n$-simplices correspond to the extreme points for the   distance function $d_L(z,y)$ restricted to  $\ds \partial_{\ell}^i \cC_X(z)$, which we denote by 
 $\{x_{i_1}, \ldots , x_{i_m} \}$. Each such point $x_{i_j}$ then corresponds to an $n$-simplex, which contains both points $\{z, z_{\ell}^i\}$ as vertices by 
 Proposition~\ref{prop-duality}. Thus, the face $\partial_{\ell}^i \cC_X(z) $ is partitioned into closed regions corresponding to its intersection with the $n$-simplices determined by the extreme points $x_{i_j}$ for $1 \leq j \leq m$. Thus, each face $\ds \partial_{\ell}^i \cC_X(z)$ is contained in the union of the realizations of the simplices satisfying $\ds  \Delta(z, z_1, \ldots , z_n) \in \Delta^{(n)}(\cN_X)$. The inclusion $ \cC_L(x) \subset \cC_{\Delta}(z)$ follows.
 \endproof

\subsection{Regular Delaunay simplicial complex} The simplicial complex $\Delta(\cN_X)$ may have non-trivial $(n+1)$-simplices, where some collection of $(n+1)$-hyperplanes satisfy
$$ L(z_0, z_1) \cap \cdots \cap L(z_0, z_{n+1}) \cap D_{\F}(z_0, \lF) \ne \emptyset ~ .$$
This is a degenerate condition, as typically every collection of $(n+1)$-hyperplanes in $D_{L}(z_0, \lF)$
should have empty intersection. 
This motivates the following definition.
\begin{defn}\label{def-regularcomplex}
The simplicial complex $\Delta(\cN_X)$ is \emph{regular} if
$\Delta^{(n+1)}(\cN_X) = \emptyset$. We say that the net $\cN_X$ is \emph{regular} if $\Delta(\cN_X)$ is regular.
\end{defn}
Note that regularity of the net $\cN_X$, and hence the complex   $\Delta(\cN_X)$,  is an open condition.
Much of the technical work in later sections   is to give conditions on a regular net $\cN_X$  such that for a net $\cN_{X'}'$ sufficiently close to $\cN_X$, the Delaunay simplicial complex $\Delta(\cN_{X'}')$ is also   regular.

\section{Foliated Voronoi structure} \label{sec-FVS}

The goal of this section is to develop  the notion of a   \emph{nice stable transversal}  $\cX$ as in Definition~\ref{def-nst}.

The constructions above for the net $\cN_X$ for  $X \subset L$  is first extended to a net in   $\wtK \subset \wtL \subset \fN_0$, then the key idea is to introduce   a parametrized set of vertices in an open neighborhood of $\wtK$ in $\fN_0$. Such an extension gives rise to a collection of transversals to $\F$, for  which we impose regularity conditions with respect to the construction of the Voronoi cells and the Delaunay triangulation.

We assume a leaf $L_0 \subset \fM$ is given, with basepoint $x_0 \in L_0$ and holonomy covering $\wtL_0$. 
As before,  assume  that  $d_2 \leq \eFU/5 <   \eFU/4 = e_2$,   that $\cM_0$ is an  $(e_1 , e_2)$-net  for $L_0$, and 
 for each $z \in \cM_0$ there is an index  
$1 \leq i_z \leq \nu$ so that $D_{\fM}(z, \eU) \subset U_{i_z}$.

 Let     $\wtcM_0 = \Pi^{-1}(\cM_0)$   be the lifted   $(e_1 , e_2)$-net for $\wtL_0$. The points of $\wtcM_0$ are denoted by $\wtz$,    where $\wtz$ is a lift of $z \in \cM_0$, and  $\wtx_0 \in \wtcM_0$ is the lift of the   basepoint $x_0 \in L_0$. 

   Let $\wtfN_0$ be the   foliated microbundle associated to the net $\wtcM_0$. 
For each $\wtz \in \wtcM_0$,    $\wtU_{\wtz} = \oU_{i_z} \times \{\wtz\}$ is the corresponding foliation chart for $\wtfN_0$. The   Riemannian metric on leaves in $\wtfN_0$ is induced by the local covering maps to leaves in $\fM$.

For $\wtz \in \wtcM_0$ and $\wty = (x, \wtz) \in \wtU_{\wtz}$, let
$\wtcP_{\wtz}(\wty) = \cP_{i_z}(x) \times \{\wtz\}$ denote the plaque of $\wtU_{\wtz}$ containing $\wty$.  
Note that by choice, $D_{\wtL_0}(\wtz , \eFU) \subset \wtcP_{\wtz}(\wtz)$ for each $\wtz \in \wtcM_0$,
and as $\cM_0$ is $e_2$-dense,  the collection
$\ds \{\wtcP_{\wtz}(\wtz) \mid \wtz \in \cN_0 \}$ is an open covering of $\wtL_0$.

Let   $\wtK \subset \wtL_0$  be a connected compact subset   which is a union of plaques,   such that the composition
$\ds \iota_0 \colon \wtK \subset \wtL_0 \to L_0 \subset \fM$
is injective with image $K$.  Assume there is given a    $(d_1,d_2)$-net $\cN_K$ for $K$, which lifts to a $(d_1,d_2)$-net $\wtcN_K$ for $\wtK$.

  We next introduce a sequence of basic concepts used in  our constructions. First is the notion of   transversals which are in ``standard form'' with respect to the chosen foliation  covering of $\fM$.
  
  \begin{defn}\label{def-standtrans}
A closed subset $\cX \subset \fM$  is a \emph{standard transversal}  if $\cX = \cX_1 \cup \cdots \cup \cX_{p}$ is a disjoint union, where for each $1 \leq \ell \leq p$, there exists a foliation chart
$\vp_{i_{\ell}}$, clopen subset $X_{\ell} \subset \fT_{i_{\ell}}$ and basepoint $v_{\ell} \in (-1,1)^n$ such that
$\ds \cX_{\ell} = \vp_{i_{\ell}}^{-1}\left(v_{\ell}, X_{\ell} \right)$. 
\end{defn}

  \begin{defn}\label{def-standtranscover}
A closed subset $\wtcX \subset \wtfN_0$   is a \emph{standard transversal}  if $\wtcX = \wtcX_1 \cup \cdots \cup \wtcX_{p}$  is a disjoint union, where for each $1 \leq \ell \leq p$, there exists $\wtz_{\ell} \in \wtcM_0$ so that for the  foliation chart
$\ds \wtvp_{\wtz_{\ell}}  \colon \wtU_{\wtz_{\ell}} \to [-1,1]^n \times \fT_{\wtz_{\ell}}$ there is a 
 clopen subset $X_{\ell} \subset \fT_{\wtz_{\ell}}$ and basepoint $v_{\ell} \in (-1,1)^n$ such that
 $\ds \wtcX_{\ell} = \wtvp_{\wtz_{\ell}}^{-1}\left(v_{\ell}, X_{\ell} \right)$. 
\end{defn}
Note that this definition just ensures that the sets $\wtcX_{\ell}$ have a standard form in local coordinates, but does not assert that the sets form a complete transversal for $\wtfN_0$.
We consider next the standard transversals which are ``holonomy invariant'' in $\ds \wtfN_0$.
  \begin{defn}\label{def-standtranscover2}
A  standard transversal   with $\wtcX = \wtcX_1 \cup \cdots \cup \wtcX_{p} \subset \wtfN_0$  is \emph{$\cGF$-invariant} if for each $1 < \ell \leq p$,  there exists a leafwise path $\wtgamma_{\ell} \colon [0,1] \to   \wtL_0$ with $\wtgamma_{\ell}(0) = \wtz_1$ and $\wtgamma_{\ell}(1) = \wtz_{\ell}$ such that
$\ds X_{1} \subset D(h_{\wtgamma_{\ell}})$ and $\ds X_{\ell} = h_{\wtgamma_{\ell}}(X_1)$.
\end{defn}
Recall that the induced foliation $\wtF$ on $\wtfN_0$ is without germinal holonomy, so this notion is independent of the chosen paths $\ds \wtgamma_{\ell}$, as long as the domain conditions are satisfied. 
The next conditions are concerned with the extensions of the notions of sections~\ref{sec-VT} and \ref{sec-DSC}. Recall that if  $x,y \in \fM$ and are not on the same leaf, then   $\dF(x,y) = \infty$, and similarly for  $\wtx, \wty \in \wtfN_0$.

\begin{defn}\label{def-reguniformtransversal}
Let $\wtfR \subset \wtfN_0$ be a given closed subset, and suppose that $\wtcX \subset \wtfR$ is a standard transversal, defined as above. 
Then $\wtcX$ is $(d_1, d_2)$--\emph{uniform} on  $\wtfR$ if there exists $0 < d_1 < d_2 \leq \lF/5$ such that 
for each $\wtx \ne \wty \in \wtcX$ we have $\dF(\wtx, \wty) \geq d_1$, and for each $\wty \in \wtfR$ there exists $\wtx \in \wtcX$ with 
$\dF(\wtx, \wty) \leq d_2$.
\end{defn}
The $(d_1, d_2)$--uniform assumption above implies that $\wtcX$ is a complete transversal for $\wtfR$, as every point   lies within leafwise distance $d_2$ of a point of $\wtcX$.

Now assume there is given a $(d_1, d_2)$-uniform   transversal $\wtcX$ for $\wtfR$.
The nearest--neighbor distance function $\kL$ extends to a leafwise function, 
\begin{equation}\label{eq-kF}
\wtkF(\wty) = \inf \left\{ \wtdF(\wtx, \wty) \mid \wtx \in \wtcX \right\},
\end{equation}
 Then we extend the definition of the Voronoi cells  to holonomy coverings by setting, for $\wty \in \wtcX$, 
\begin{equation}\label{eq-voronoileafwise}
\cC_{\wtfR}(\wty) = \left\{ \wtz \in \wtfR  \mid \wtdF(\wtz, \wty) = \wtkF(\wtz)\right\} .
\end{equation}
In other words,   $\cC_{\wtfR}(\wtx)$ is the Voronoi cell in $\wtL_{\wtx} \cap \wtfR$ defined by the net $\wtcX(x) = \wtcX \cap \wtfR$, which consists of  $\wty \in \wtL_{\wtx} \cap \wtfR$ which are closer to $\wtx$ in the leafwise metric $\wtdF$ than to any other point of $\wtcX$. By the definition of the $(d_1, d_2)$-uniform   transversal $\wtcX$ for $\wtfR$, each $\wty \in \wtfR$ belongs to at least one such cell.

 The leafwise Voronoi cells $\cC_{\wtfR}(\wtx)$  can be organized into  \emph{Voronoi cylinders} using the decomposition 
 $\ds \wtcX = \wtcX_1 \cup \cdots \cup \wtcX_{p}$. 
For $1 \leq \ell \leq p$ define the \emph{Voronoi cylinder} by
\begin{equation}\label{eq-voronoicylinder}
\fC^{\ell}_{\wtfR} = \bigcup_{\wtx \in \wtcX_{\ell}} ~ \cC_{\wtfR}(\wtx)
\end{equation}
 We thus obtain the \emph{Voronoi decomposition} $\ds \wtfR = \fC_{\wtfR}^1 \cup \cdots \cup \fC_{\wtfR}^p$ associated to $\wtcX$.

The notion of the \emph{star-neighborhood} of a Voronoi cell, given in Definition~\ref{lem-star-nbhd}, extends immediately to the Voronoi cells in each   leaf. 
First, for  $\wtx \in \wtcX$ introduce the \emph{vertex-set}
\begin{equation}\label{eq-vertexsetF}
\cV_{\wtfR}(\wtx) = \{ \wty \in \wtcN \mid \cC_{\wtfR}(\wty) \cap \cC_{\wtfR}(\wtx) \ne \emptyset \} .
\end{equation}
The  star-neighborhood  of the Voronoi cell $\cC_{\wtfR}(\wtx)$ is the set
$\ds \cS_{\wtfR}(\wtx) =\bigcup_{\wtx \in \cV_{\wtfR}(\wtx)} ~ \cC_{\wtfR}(\wty)$. 
Then for each $1 \leq \ell \leq p$,  define the  \emph{star-neighborhood} of the cylinder $\ds \fC^{\ell}_{\wtfR}$ by
\begin{equation} \label{eq-starneighborhoodF}
\fS_{\wtfR}^{\ell} =\bigcup_{\wtx \in \wtcX_{\ell}} ~ \cS_{\wtfR}(\wtx) .
\end{equation}
Lemma~\ref{lem-star2} shows that for  a point $x \in \cN_K$, the   star-neighborhood $\cS_K(x)  \subset B_{L}(x, 3d_2)$ so that $d_2 \leq \lF/5$  implies $\cS_K(x)$ is contained in some  coordinate chart $U_{z}$ for $z \in \cM_0$. 
 For the   star-neighborhood  $\ds \fS_{\wtfR}^{\ell}$ of the cylinder $\ds \fC^{\ell}_{\wtfR}$, the conclusion that it is contained in some foliation chart for $\fM$ is not a priori satisfied,  and this condition is imposed as one of our assumptions. It will    later be checked that it is satisfied for the   transversals constructed.

 \begin{defn}\label{def-starcontainment}
Let $\wtfR \subset \wtfN_0$ be a given closed subset, and suppose that $\ds \wtcX = \wtcX_1 \cup \cdots \cup \wtcX_{p}$ is a uniform standard transversal. 
Then we say that $\wtcX$ is \emph{centered} if for each $1 \leq \ell \leq p$ there is a coordinate chart $U_{\nu_{\ell}}$ for   $\nu_{\ell} = i_{\wtx}$ for some $\wtx  \in \cM_0$ such that 
$\ds \fS_{\wtfR}^{\ell} \subset U_{\nu_{\ell}}$. 

A  standard  transversal $\wtcX$  for $\wtfR$  is said to be \emph{nice} if it is   $(d_1, d_2)$-uniform, invariant  and centered.
\end{defn}

Next, define the leafwise simplicial complex $\Delta_{\F}(\wtcX)$ associated to a nice    transversal $\wtcX$ for $\wtfR$.  The method of circumscribed spheres adapts immediately, as follows.
\begin{defn}
Let $\wtcX$ be a  nice  transversal  for $\wtfR$. 
The collection of points $\ds \{\wtz_0 , \ldots, \wtz_k\} \subset \wtcX$ defines a 
 $k$-simplex $\ds \Delta(\wtz_0 , \ldots, \wtz_k) \in \Delta_{\F}(\wtcX)$ if
 there exists $\ds \wtz \in \wtL_{\wtz_0}$ and $r \leq d_2$ such that
 \begin{equation}
\{\wtz_0, \ldots , \wtz_k\} \subset S_{\wtF}(\wtz,r) \cap \wtcX \quad , \quad 
B_{\wtF}(\wtz,r) \cap \wtcX = \emptyset
\end{equation}
\end{defn}
Again, note that $\wtdF(\wtx, \wty) = \infty$ if $\wtx$ and $\wty$ lie on distinct leaves, so $\ds \Delta(\wtz_0 , \ldots, \wtz_k) \in \Delta_{\F}(\wtcX)$ implies that 
$\ds \{\wtz_0 , \ldots, \wtz_k\} \subset \wtL_{z_0}$. Consequently,  $\ds \Delta(\wtz_0 , \ldots, \wtz_k)$ can also be considered as a $k$-simplex for  $\wtcX \cap \wtL_{z_0}$, so that $\Delta_{\F}(\wtcX)$ consists of a union of   simplices contained in the leaves of $\wtF$.  
The key question is then,   given 
$\ds \Delta(\wtz_0 , \ldots, \wtz_k) \in \Delta_{\F}(\wtcX)$,  is it contained in a transverse family of simplices? We make this property precise, as it is   fundamental.

 Let $1 \leq j_0, \ldots , j_k \leq p$ be indices such that $\wtz_{\ell} \in \wtcX_{j_{\ell}}$ for $0 \leq \ell \leq k$.
 In particular,  $\wtz_0 \in \wtcX_{j_{0}} \subset U_{i_{j_{0}}}$ and thus $\wtz_0 \in \wtcP_{i_{j_{0}}}(\wtz_0)$. 
As   $\wtcX$ is centered,   
   we can assume    $z_{\ell} \in \cP_{i_{j_{0}}}(z_0)$ for $1 \leq \ell \leq k$  
and   $\cX_{j_{\ell}} \subset U_{i_{j_{0}}}$. 

Note  that
for $\ell \ne \ell'$ the sets $\cX_{j_{\ell}}$ are $\cX_{j_{\ell'}}$ are disjoint by the $(d_1, d_2)$-net hypothesis.

Let $z_0' \in \cX_{j_0}$.
Let $\cP_{i_{j_{0}}}(z_0')$ denote the plaque of $U_{i_{j_{0}}}$ containing $z_0'$.
For each $1 \leq \ell \leq n$, let $z_{\ell}' = \cX_{i_{\ell}} \cap \cP_{i_{j_{0}}}(\xi_0')$ be the unique point of $\cX_{i_{\ell}}$ contained in the plaque defined by $z_0'$. Observe that the points $z_{\ell}'$ depend continuously on $z_0' \in \cX_{j_0}$.
\begin{defn}\label{def-nst}
Let $\wtcX$ be a nice  transversal  for $\wtfR$. Then $\wtcX$ is   \emph{stable} if for each $k$-simplex 
$\Delta(\wtz_0 , \ldots, \wtz_k) \in \Delta_{\F}(\wtcX)$ and
$\wtz_0' \in \cX_{j_0}$, we have $\Delta(\wtz_0' , \ldots, \wtz_k') \in \Delta_{\F}(\wtcX)$.
\end{defn}

At first inspection, stability of simplices for a Delaunay triangulation associated to a $\wtcX$ seems to be intuitively clear, and in fact this is basically correct for dimension $n \leq 2$ as the $(d_1, d_2)$-net hypothesis implies stability for \emph{planar} tessellations.
The difficulty is that for $n > 2$, as the transverse coordinate $\wtz_0' \in \cX_{j_0}$  varies,   ``small variations'' of the spacings of the net points of $\wtcX \cap \wtL_{\wtz_0' }$   may result in an abrupt change in the Delaunay simplicial structure, if some face of a Voronoi cell has too small of a diameter relative to the size of the variation. Consequently,  the  existence of a nice stable transversal   for $n > 2$ requires the delicate estimates in its   construction in later sections.

\section{Constructions of transverse Cantor foliations} \label{sec-proofstableappr}

In this section, we show how   the existence of  a nice stable transversal    is used to construct  a transverse Cantor foliation $\cH$ on a given set $\fB$ and consequently a product structure on an open neighborhood. We use notation   as in  the previous sections.

For $x \in \fM$,  assume there is given  a connected compact subset  $K_x \subset L_x$ such that there is    $\wtK_x \subset \wtL_x$   such that 
$\ds \iota_x \colon \wtK_x \subset \wtL_x \to L_x \subset \fM$
is injective with image $K_x$. That is, $K_x$ is a proper base as in Definition~\ref{def-LBB}. We introduce below an extension $\wtK_x \subset \whK_x$ and choose a $\whK_x$-admissible transversal $V_x$ containing $x$ so that $\ds  \fN(\whK_x, V_{x})$ is well-defined. 
Here is the main result:

 \begin{thm}\label{thm-stableapprox}
For     $0 < d_1 < d_2 \leq \lF/5$,  assume there is given a nice stable $(d_1, d_2)$--uniform transversal $\wtcX$ for $\ds  \fN(\whK_x, V_{x})$. Then there exists a foliated homeomorphism into, 
\begin{equation}\label{eq-localproduct}
\Phi \colon  \wtK_x \times V_x \to  \wtfN(\whK_x, V_{x})
\end{equation}
 such that the images $\ds \Phi \left(  \{\wty\} \times V_x \right)$, for  $ \wty \in \wtK_x$  define a continuous family of Cantor transversals for $\wtF | \wtfN(\whK_x, V_{x})$ which extend the transversals in $\wtcX$. 
 \end{thm}
Thus, the assumption there is a nice stable transversal for $\whK_0 \subset \fN(\whK_x, V_{x})$ implies there is a  transverse Cantor foliation $\wtcH$ defined on some open neighborhood of   $\wtK_x$ in $\ds \fN(\whK_x, V_{x})$.
The proof of Theorem~\ref{thm-stableapprox} occupies the rest of this section.

We first introduce a sequence of modifications of the set $\wtK_x$, first to expand the set, then translate it to a leaf $L_0$ without holonomy, resulting in the set $\whK_0 \subset \wtL_0$.  
Without loss of generality, we may assume that  $x \in K_x$, and then let $\wtx \in \wtK_x$ be the lift of $x$, which is unique as $\wtK_x$ injects into $\fM$. 

Let  $U_{i_x}$ be the foliation chart with  $B_{\fM}(x, \eU) \subset U_{i_x}$, and  set    $w_x = \pi_{i_x}(x) \in \fT_{i_x}$. 
Since $\fM$ is minimal, $\fT_{i_x}$ is a Cantor set, 
and thus  $w_x$ is not an isolated point.

 Given complete separable metric space $(X, d_X)$,  a proper subset $Y \subset X$  and $\e > 0$, introduce  the notion of the \emph{$\e$-penumbra} of $Y$ in $X$,  
\begin{equation}\label{eq-penumbra}
\PX(Y, \e) = \{x \in X \mid d_X(x, Y) \leq \e\}.
\end{equation}
 That is, $\PX(Y, \e)$ is the closed subset of $X$ consisting of all points within distance $\e$ of $Y$. 
 We apply this construction for $\e = \lF$ to $\wtK_x \subset X = \wtL_x$. Then by definition,   for every point $\wty \in \wtK_x$ we have  $D_{\F}(\wty, \lF) \subset \PF(\wtK_x, \lF)$. Let  $\whK_x$ be the plaque saturation of $\PF(\wtK_x, \lF)$ in $\wtL_x$, so 
\begin{equation}\label{eq-translationK0}
\whK_x~ = ~ \bigcup  ~ \left\{ \wtcP_{\wtz}(\wtz) \mid \wtz \in \wtL_x ~ , ~ \wtcP_{\wtz}(\wtz) \cap  \PF(\wtK_x, \lF) \ne \emptyset \right\}
\end{equation}
Let $\whR_K$ denote the diameter of the set $\whK_x$ in $\wtL_x$. It follows that  $\whK_x \subset D_{\wtF}(\wtx, \whR_K)$.

Recall from Proposition~\ref{prop-domest} that given   $\epsilon > 0$, there exists
$0 < \delta(\epsilon, \whR_K) \leq \epsilon$ so that for any clopen neighborhood $V_x$ with diameter  at most $\delta(\epsilon, \whR_K)$ in $\fT_{i_x}$,   the set $V_x$ is   $\whK_x$-admissible, as defined by Definition~\ref{eq-tesselspace2}. The  choice of $\e >0$ and the clopen neighborhood $V_x$ will be specified in later sections, based on the radius $\whR_K$ and estimates derived from  the leafwise Riemannian geometry.  
For now, we assume they are   given.

The germinal holonomy of the leaf $L_x$ is given by the isotropy subgroup    $\G_{\F}^{w_x}$ defined in \eqref{eq-holodef}, which  is represented by the   elements of the pseudogroup $\cGF^*$ which fix $w_x$.  
Then for the clopen neighborhood $w_x \in V_x \subset \fT_{i_x}$,   Theorem~\ref{thm-emt} implies there exists $w_0 \in V_x$ such that the leaf $L_0$ corresponding to  $w_0$ is without holonomy. If  $L_x$ is without holonomy, then we may take $w_0 = w_x$.  Let $\Pi \colon \wtL_0 \to L_0$ be the holonomy cover, which is a diffeomorphism.

 Form the Reeb neighborhood $\ds \wtfN(\whK_x, V_{x}) \subset  \wtfN_0 $ as in Definition~\ref{def-admissibledisjoint}. 
Let $\whK_0$ be the connected compact subset of the holonomy cover $\wtL_0$ obtained by taking the union of the plaques in $\wtL_0$ which are contained in $\ds \wtfN(\whK_x, V_{x})$, so  $\whK_0$ is a ``translation'' of  $\whK_x$ to  the  leaf $\wtL_0$. Note that $V_x$ is also $\whK_0$-admissible by Lemma~\ref{lem-domainconst} and that $\ds  \wtfN(\whK_0, V_{x}) = \wtfN(\whK_x, V_{x})$.

The assumption of Theorem~\ref{thm-stableapprox} is that 
  there is   a nice stable $(d_1, d_2)$--uniform transversal  $\ds \wtcX = \wtcX_1 \cup \cdots \cup \wtcX_{p}$ for $\ds \whK_0 \subset \fN(\whK_x, V_{x})$. 
 That is, we have a product structure for $\ds \fN(\whK_0, V_{x})$ defined on the net $\wtcX \cap \whK_0$ and this must be extended to all of $\ds \fN(\whK_0, V_{x})$.

Let $\Delta(\wtz_0 , \ldots, \wtz_k) \in \Delta_{\F}(\wtcX)$ be given, with $\wtz_{\ell} \in \wtcX_{i_{\ell}}$. Then  
$\{\wtz_0, \ldots , \wtz_k\} \subset \wtU_{i_0}$ so we have $\ell \ne \ell'$ implies that $i_{\ell} \ne i_{\ell'}$.
Without loss of generality, we may re-index the vertices so that $\ell < \ell'$ implies $i_{\ell} < i_{\ell'}$. The indexing  of the sets $\wtcX_{\ell}$ yields an ordering of the vertices of $\Delta(\wtz_0 , \ldots, \wtz_k)$. This is the   ``local ordering'' referred to in Remark~\ref{rem-nonuniquedelone}.

The transversal $\ds \wtcX$ defines the leaves of the Cantor foliation $\wtcH$ through each vertex $\{\wtz_0 , \ldots, \wtz_k\}$. We next   show how to extend this finite collection of leaves  to a foliation through the faces and interior of the simplex $\Delta(\wtz_0 , \ldots, \wtz_k)$.

For each $\wtz_{0}' \in \wtcX_{i_0}$ and $1 \leq \ell \leq k$, let $\wtz_{\ell}' = \wtcP_{i_0}(\wtz_0') \cap \wtcX_{i_{\ell}}$. The stable hypothesis then implies that
$\Delta(\wtz_0', \ldots , \wtz_k') \in \Delta_{\F}(\wtcX)$.
By Lemma~\ref{lem-geofill}, for each $\wtz_0' \in \wtcX_{i_0}$ there exists a geodesic filling map $\sigma_{k, \wtz_0'} \colon \Delta^k \to \wtcP_{i_0}(\wtz_0') \subset \wtL_{\wtz_0'}$ associated to $\Delta(\wtz_0', \ldots , \wtz_k')$ which is natural with respect to the face maps. We thus obtain a continuous  map
\begin{equation}\label{eq-paramfilling}
\Sigma_{i_0} : \Delta^k \times \wtcX_{i_0} \to \fR: (\vec{v}, \wtz_0') \mapsto \sigma_{k,\wtz_0'} (\vec{v}) \quad , \quad \wtz_0'  \in \wtcX_{i_0} ~ , ~   \vec{v} \in  \Delta^k 
\end{equation}
For each   $\vec{v} \in \Delta^k$ 
and $\wtz_0' , \wtz_0'' \in \wtcX_{i_0}$   define $\ds \sigma_{k, \wtz_0'}(\vec{v}) \approx \sigma_{k, \wtz_0''}(\vec{v})$. The equivalence class of $\ds \sigma_{k, \wtz_0}(\vec{v})$ defines a Cantor transversal through the point, which is   a leaf of $\wtcH$.

The foliation $\wtcH$ is defined by the equivalence classes of points in    the interiors of the geometric realizations of the simplices in  $\ds \Delta_{\F}(\wtcX)$.  On the faces of adjacent simplices, the local orderings are compatible, so the geodesic filling maps agree, and thus so does the equivalence relation $\approx$.

We underline some   points of this construction.
First, for each $1 \leq \ell \leq p$ and $\wtz, \wtz' \in \wtcX_{\ell}$ then $\wtz \approx \wtz'$. That is, each transversal $\wtcX_{\ell}$ is a  leaf of the foliation $\wtcH$.

Second, for each $1$-simplex $\Delta(\wtz_0, \wtz_1) \in \Delta_{\F}(\wtcX)$ the equivalence relation $\approx$ identifies points with the same barycentric coordinate on the unique geodesic ray joining $\wtz_1'$ to $\wtz_0'$ where $\Delta(\wtz_0', \wtz_1')$ is the transverse transport of the given $1$-simplex. Thus, $\approx$ is independent of the ordering when restricted  to the $1$-skeleton of the leafwise triangulation $ \Delta_{\F}(\wtcX)$.
If $\F$ is an orientable foliation by 1-dimensional leaves, that is, it is defined by a flow, then we are done, and the equivalence relation $\approx$ depends canonically on the choice of the uniform transversal $\wtcX$, but is independent of its ordering.

If the leaves of $\F$ have dimension $n > 1$, then the ``spanning geodesic procedure'' in the proof of Lemma~\ref{lem-geofill} may well depend upon the ordering of the vertices in each simplex. However, the ``local ordering'' of the vertices in simplices is determined by the choice of the transversal $\wtcX = \wtcX_1 \cup \cdots \cup \wtcX_{p}$. Thus $\approx$ is well-defined, assuming the choice of the transversal $\wtcX$ with its ordering.

Finally, we must show the map \eqref{eq-localproduct} is well defined. For this, we show  that   $\wtK_x$ is contained in the domain  of the equivalence relation $\approx$.

 For each $\wty \in \wtcX$, with leaf $\wtL_{\wty}$ containing it,  
 the intersection $\wtcN_{\wty} = \wtcX \cap \wtL_{\wty}$ is a $(d_1 , d_2)$-net for $\fN(\whK_x, V_{x}) \cap \wtL_{\wty}$
by Definition~\ref{def-reguniformtransversal}.
  Define the function $\wtkF$ as in \eqref{eq-kF}  and the Voronoi cell, as in \eqref{eq-voronoileafwise}, 
  $$\cC(\wty) =  \{ \wtz \in \wtL_{\wty} \mid \dF(\wtz, \wty) = \wtkF(\wtz) \}$$ 
  Let $\wtcN_{\wty}^* \subset \wtcN_{\wty}$ be the subset of points for which $\cC(\wty) \subset \fN(\whK_x, V_{x})$, and    
    $\wtcN_0 = \wtcX \cap \wtL_0$ with   $\wtcN_{0}^* \subset \wtcN_{0}$.

 Let    $\cC_{\Delta}(\wty)$ be the simplicial cone of $\wty \in \wtcX$ in the complex $\Delta_{\F}(\wtcX)$ as defined leafwise by  \eqref{eq-simplicialcone}. 
Then the stability assumption on $\wtcX$ implies the simplicial complex $\Delta(\wtcN_{\wty})$ is stable, and we define
 \begin{equation}\label{eq-B}
\fB ~ = ~ \bigcup_{\wtz_i \in \wtcN_0^*} ~ \bigcup_{\wty \in \wtcX_i} ~\cC_{\Delta}(\wty)
\end{equation}
Then the equivalence relation $\approx$ is defined on $\fB$ by definition, so $\fB$ admits a Cantor foliation $\cH$. 

The proof of Theorem~\ref{thm-stableapprox}   then follows from   
 \begin{lemma}\label{lem-KinB}
  $\wtK_x \subset \fB$.
 \end{lemma}
\proof
By Definition~\ref{def-reguniformtransversal} and the construction of $\fN(\whK_x, V_{x})$,   
for each $\wty \in \wtK_x$ there exists $\wtz \in \wtcX_{\ell}$ such that $d_{\wtF}(\wty, \wtz) \leq d_2$. 
We may assume in addition that $\wtz$ is a closest point in $\wtcX$, so that $\wty \in \cC(\wtz)$. 
Then note  that 
$$D_{\wtF}(\wty, \lF) \subset \PF(\wtK_x, \lF)  \subset \whK_x$$ 
by construction, and $d_2 \leq \lF/5$ implies that 
$D_{\wtF}(\wtz, 4d_2) \subset D_{\wtF}(\wtz, 4\lF/5) \subset \whK_x$ as well. 

In particular, $D_{\wtF}(\wtz, 4d_2) \subset \whK_x$, so for each $\wty' \in \cV_{\wtfR}(\wtz) $ as defined in \eqref{eq-vertexsetF}, we have 
$D_{\wtF}(\wty', 2d_2) \subset \whK_x$. This implies $\wty' \in \wtcN_{\wty}^*$ by the extension of Lemma~\ref{lem-netdensity}.
Consequently,  the star-neighborhood,  as defined in \eqref{eq-starneighborhoodF}  for $\ds \wtfR = \fN(\whK_x, V_{x})$, satisfies 
$\ds \fS_{\wtfR}^{\ell}(\wtz)  \subset B_{\wtL_x}(\wtz, 3d_2) \subset \whK_x$.

 By Proposition~\ref{prop-filling}, for $\wtz \in \wtcN_0^*$ we have    $ \cC(\wtz) \subset \cC_{\Delta}(\wtz)$,  
 where $\cC_{\Delta}(\wtz)$ is the simplicial cone of $\wtz$ in the complex $\Delta_{\F}(\wtcX)$ as defined by  \eqref{eq-simplicialcone}. 
Thus we have 
 \begin{equation}
\wty \in  \cC(\wtz) \subset \cC_{\Delta}(\wtz) \subset  \fB
\end{equation}
which completes the proof of Lemma~\ref{lem-KinB} and so also Theorem~\ref{thm-stableapprox}.
\endproof

\section{Delaunay simplices in Euclidean geometry} \label{sec-euclidean}

The construction of a Delaunay triangulation from a point-set in $\mR^n$ using the Voronoi tessellation it defines is well-known, and a fundamental tool in computational geometry \cite{OBKC2000}. The application of this method in the case of complete Riemannian manifolds is less well developed, except for hyperbolic space and some other variants of the standard Euclidean metric. See  \cite{Clarkson2006} and  \cite{LeibonLetscher2000} for discussions of some of the aspects of adapting the Euclidean methods to a non-Euclidean framework. 

Our construction of a nice stable transversals uses  the construction of Voronoi tessellations of the leaves to obtain   stable Delaunay triangulations, and for this  we require detailed estimates on the properties  of the construction, especially with respect to a transverse parameter. In the next few sections, we develop the estimates required. The techniques are almost all based on methods of ``elementary''  linear algebra \cite{Hogben2007}, but the applications to our situation are more specialized.  For example, given a collection of vectors $\{\vec{y}_{0}, \ldots , \vec{y}_n\}$ in $\mR^n$ which admit a circumscribed sphere with
center $ \omega(\vec{y}_0, \ldots , \vec{y}_n) $ and radius $r(\vec{y}_0, \ldots , \vec{y}_n)$, we  derive a stability criterion, conditions for which a small displacement
$\{\vec{z}_{0}, \ldots , \vec{z}_n\}$ of $\{\vec{y}_{0}, \ldots , \vec{y}_n\}$  still  uniquely defines a   circumscribed sphere.

\subsection{Preliminaries}
Let $\mR^n$ have the standard Euclidean metric $d_{\mR^n}$ and associated norm $\| \cdot \|$. To fix notation, we consider   $\vec{x} \in \mR^n$ as a column vector, and let
$\vec{x} \bullet \vec{y} = \vec{x}^t \cdot \vec{y}$ denote the ``dot-product" of two vectors, where $\vec{x}^t$ denotes the matrix transpose of $\vec{y}$, and $\vec{x}^t \cdot \vec{y}$ denotes the   matrix product.

Given a collection of $n$ vectors, $\{\vec{a}_{1}, \ldots , \vec{a}_n\} \subset \mR^n$, let $\bA$ denote the $n \times n$ matrix with these vectors as \emph{rows}. Let 
$\ds  \| \bA \| ~ = ~ \max \left\{\| \bA \cdot \vec{x} \| ~ | ~ \vec{x} \in \mR^n ~ , ~ \|\vec{x} \| = 1\right\}$ be the operator norm for $A$.
If $\bA$ is a diagonal matrix with entries $\{\lambda_1 , \ldots , \lambda_n\}$, then the norm is calculated  by
\begin{align}\label{eq-diagdef}
\|\bA \| ~ = ~ \max~ \{ |\lambda_1| , \ldots , |\lambda_n|\}
\end{align}
and in general, the Cauchy-Schwartz inequality   yields the estimate
\begin{align}\label{eq-normdef}
\|\bA \|^2 ~ \leq ~ \|\vec{a}_1\|^2 + \cdots + \|\vec{a}_n\|^2.
\end{align}

We recall an elementary result of linear algebra:
\begin{lemma}\label{lem-volume}
Let $\{\vec{y}_{0}, \ldots , \vec{y}_n\} \subset \mR^n$ and form their convex hull
 $$\Delta (\vec{y}_{0}, \ldots , \vec{y}_n ) = \left\{ t_0 \  \vec{y}_{0} + \cdots +  t_n \ \vec{y}_{n}  \mid t_0 + \cdots + t_n = 1 ~ , ~ t_i \geq 0 \right\}$$
Fix $0 \leq \ell \leq n$, and let  $D_{\ell}(\vec{y}_{0}, \ldots , \vec{y}_n)$ be the $n \times n$-matrix whose rows are the transposes of the vectors $\vec{y}_i - \vec{y}_{\ell}$ for $i \ne \ell$. Then
\begin{equation}
| \det  D_{\ell}(\vec{y}_{0}, \ldots , \vec{y}_n) | =   n! \cdot {\rm Vol}(\Delta (\vec{y}_{0}, \ldots , \vec{y}_n ))
\end{equation}
\end{lemma}

Now assume we are given a collection   $\{\vec{y}_{0}, \ldots , \vec{y}_n\} \subset \mR^n$ which admit a circumscribed sphere with center $ \omega(\vec{y}_0, \ldots , \vec{y}_n)$ and radius $r(\vec{y}_0, \ldots , \vec{y}_n)$. The point $\omega(\vec{y}_0, \ldots , \vec{y}_n)$ is characterized as being equidistant from all of the points $\vec{y}_i$, and thus it lies on each perpendicular bisector hyperplane  $L(\vec{y}_i, \vec{y}_j)$ equidistant between the two points.

To derive the equations for the center and radius, we   the vector  $\vec{v}_{n}$ as a ``base point'', though clearly the solutions do not depend on which vertex is chosen.  
For  each $1 \leq k \leq n$, set $\vec{u}_{k} = (\vec{y}_{k-1} - \vec{y}_n)$. Then $\|\vec{u}_k \| \leq 2 \, r(\vec{y}_0, \ldots , \vec{y}_n)$ as all vectors
$ \vec{y}_i$ are contained in a set with diameter $2 \, r(\vec{y}_0, \ldots , \vec{y}_n)$.
Let $\bU$ denote the $n \times n$ matrix whose rows are the transposes of the vectors $\vec{u}_{k}$.
Let $| \bU | = |\det \bU |$ denote the absolute value of the determinant of $\bU$. 
Then  $| \bU | = n! \cdot {\rm Vol}(\Delta (\vec{y}_{0}, \ldots , \vec{y}_n ))$ by Lemma~\ref{lem-volume}.

The hyperplanes $L(\vec{y}_{k-1}, \vec{y}_{n})$, for $1 \leq k \leq n$,  are defined by the equations
\begin{eqnarray}
L(\vec{y}_{k-1}, \vec{y}_{n}) ~ & = & ~ \left \{ \vec{x} \in \mR^n \mid ( \vec{y}_{k-1} - \vec{y}_n) \bullet (\vec{x} - ( \vec{y}_{k-1} + \vec{y}_n)/2 ) = 0 \right \} \nonumber \\
~ & = & ~ \left \{\vec{x} \in \mR^n \mid \vec{u}_{k} \bullet \vec{x} = 1/2 \cdot \vec{u}_{k} \bullet \vec{u}_{k} + \vec{u}_{k} \bullet \vec{y}_n \right \} \nonumber \\
~ & = & ~ \left \{ \xi + \vec{y}_n \mid \vec{u}_{k} \bullet \vec{\xi} = 1/2 \cdot \|\vec{u}_{k}\|^2 \right \} \label{eq-solutiony},
\end{eqnarray}
where $\xi = \vec{x} - \vec{y}_n$ represents the coordinates for $L(\vec{y}_k, \vec{y}_{n})$ with $\vec{y}_n$ translated to the origin.
The center $\vec{\omega}(\vec{y}_0, \ldots, \vec{y}_{n}) \in \mR^n$ is given by the intersection of these hyperplanes, and thus is the solution of the system of equations
\begin{align}\label{eq-solU}
\bU \cdot \xi ~ = ~ \frac{1}{2} \cdot \vl(\bU) \quad \text{so} \quad
\vec{\omega}(\vec{y}_0, \ldots, \vec{y}_{n}) ~ = ~ \frac{1}{2} \cdot \left\{ \bU^{-1} \cdot \vl(\bU) \right\} + \vec{y}_n,
\end{align}
where $\vl(\bU) = (\|\vec{u}_{1}\|^2, \ldots , \|\vec{u}_{n}\|^2)^t$ is the column vector with entries $\|\vec{u}_{k}\|^2$.

 \subsection{Effective estimates}\label{subsec-effective}
Now  suppose there exists constants $0 < e_1< e_2$ and $\ve, \delta > 0$ such that
\begin{enumerate}
\itemsep 4pt
\item $e_1 \leq \| \vec{y}_i - \vec{y}_j \|$ ~ for all $0 \leq i \ne j \leq n$,
\item $e_1/2 \leq r(\vec{y}_0, \ldots , \vec{y}_n) \leq e_2$ and hence $ \| \vec{y}_i - \vec{y}_j \| \leq 2e_2$,
\item $| \bU | \geq \delta$.
\end{enumerate}
Assume there is also    given vectors $\{\vec{z}_{0}, \ldots , \vec{z}_n\} \subset \mR^n$ such that
\begin{enumerate}\setcounter{enumi}{3}
\item $\| \vec{y}_i - \vec{z}_i \| < \ve$ ~ for all $0 \leq i \leq n$.
\end{enumerate}
We determine values of the constants $\ve, \delta > 0$ such that the points $\{\vec{z}_{0}, \ldots , \vec{z}_n\}$ admit a unique circumscribed sphere with center $ \omega(\vec{z}_0, \ldots , \vec{z}_n) $ and radius $r(\vec{z}_0, \ldots , \vec{z}_n)$, and obtain estimates for
$$\| \omega(\vec{z}_0, \ldots , \vec{z}_n) - \omega(\vec{y}_0, \ldots , \vec{y}_n) \| ~ \text{ and } ~ | r(\vec{z}_0, \ldots , \vec{z}_n) - r(\vec{y}_0, \ldots , \vec{y}_n) | .$$
Let $\bV$ denote the $n \times n$ matrix whose rows are the transposes of the vectors $\vec{v}_{k} = \vec{z}_{k-1} - \vec{z}_n$ for $1 \leq \ell \leq n$, and set
$\vl(\bV) = (\|\vec{v}_{1}\|^2, \ldots , \|\vec{v}_{n}\|^2)^t$.
Assuming that $\bV^{-1}$ exists, then for $\zeta = \vec{x} - \vec{z}_n$, the solution of the matrix equations
\begin{align}\label{eq-solV}
\bV \cdot \zeta~ = ~ \frac{1}{2} \cdot \vl( \bV) \quad , \quad \vec{\omega}(\vec{z}_0, \ldots, \vec{z}_{n}) ~ = ~ \frac{1}{2} \cdot \left\{ \bV^{-1} \cdot \vl(\bV) \right\} + \vec{z}_n
\end{align}
is the center for a unique circumscribed sphere containing the points $\{\vec{z}_{0}, \ldots , \vec{z}_n\}$.  

 Our next goal is to obtain an effective estimate on $ \| \vec{\omega}(\vec{z}_0, \ldots, \vec{z}_{n}) - \vec{\omega}(\vec{y}_0, \ldots, \vec{y}_{n}) \|$ as given by \eqref{eq-godawfulest} below. Using \eqref{eq-solU} and \eqref{eq-solV}, this will be based upon obtaining effective estimates for the matrix norms $\| \bU^{-1}\|$ and $\| \bV^{-1} \|$.
Let $\bW = \bV - \bU$ so $\bV = \bU + \bW$, and set $\bQ = \bW \bU^{-1}$.
\begin{lemma}\label{lem-Vinvest}
Assume that $\|\bQ\| \leq 1/2$, then $\bV^{-1}$ exists, and $\ds \| \bV^{-1}\| ~ \leq ~ 2 \|\bU^{-1}\|$. \hfill $\Box$
\end{lemma}

Next, the triangle inequality and our given data yield the following estimates, where $e_3 = e_2 + \ve$,
\begin{equation}
\|\vec{v}_{k} - \vec{u}_{k}\| ~ \leq ~ \|\vec{z}_{k-1} - \vec{y}_{k-1}\| + \|\vec{z}_{n} - \vec{y}_{n}\| ~ \leq ~ 2\ve,
\end{equation}
\begin{equation} \label{eq-messyest2}
e_1 - 2\ve ~ \leq ~ \|\vec{u}_{k} \| - \|\vec{v}_{k} - \vec{u}_{k}\| ~ \leq ~ \|\vec{v}_{k} \| ~ \leq ~ \|\vec{u}_{k} \| + \|\vec{v}_{k} - \vec{u}_{k}\| ~ \leq ~ 2e_2 + 2\ve = 2 e_3,
\end{equation}
then by assumption \ref{subsec-effective}.2  and \eqref{eq-messyest2}
\begin{equation}\label{eq-messyest3}
\left | \|\vec{v}_{k}\|^2 - \|\vec{u}_{k}\|^2 \right| ~ = ~ \left | (\vec{v}_{k} - \vec{u}_{k} ) \bullet (\vec{v}_{k} + \vec{u}_{k} ) \right | ~ \leq ~
\| \vec{v}_{k} - \vec{u}_{k} \| \cdot \left ( \|\vec{v}_{k} \| + \|\vec{u}_{k}\| \right ) ~ \leq ~ 4 \ve (e_2 + e_3).
\end{equation}
so that \eqref{eq-normdef} and \eqref{eq-messyest2} imply
\begin{equation}\label{eq-messyest5}
\|\bW\| ~ = ~ \| \bV - \bU \| ~ \leq ~ \sqrt{ \|\vec{v}_{1} - \vec{u}_{1}\|^2 + \cdots + \|\vec{v}_{n} - \vec{u}_{n}\|^2} ~ \leq ~ 2 \ve \sqrt{n}.
\end{equation}

We next estimate the norm $\|\bU^{-1} \|$.
Our colleague Shmuel Friedland suggested the use of the \emph{Hadamard determinantal inequality} to obtain an estimate for $\|\bA^{-1}\|$. As this is a fundamental estimate for deriving our estimates, we include a proof.
\begin{lemma} \label{lem-Aest}
Let $\bA$ be an $n \times n$-matrix whose determinant has absolute value $| \bA | > 0$, and such that each column of $\bA$ has norm at most $C$. Then
\begin{align}\label{eq-Aest}
\|\bA^{-1}\| ~ \leq ~ n \cdot C^{n-1}/| \bA |.
\end{align}
\end{lemma}
\proof For an invertible $n \times n$-matrix $\bC$, let
$0 < |\sigma_n(\bC)| \le \cdots \le |\sigma_1(\bC)|$ denote the singular values of $\bC$, ordered by their norms.
Recall that $\| \bC \|^2 = \|\bC^t \cdot \bC \| = |\sigma_1(\bC)|^2$.

Let $\adj (\bA)$ denote the adjoint of $\bA$.
Since $\bA^{-1}=\frac{1}{|\bA|} \cdot \adj (\bA)$
it follows that
the singular values of $\adj (\bA)$ are all the $(n-1)$ products of the $n$
singular values of $\bA$.
Hence the largest singular value for $\adj (\bA)$ is
$\ds \sigma_1(\adj (\bA))=\sigma_1(\bA) \cdots \sigma_{n-1}(\bA)$.

Each entry of $\adj (\bA)$ is an $(n-1)$
minor of $\bA$, and thus its
absolute value is less or equal to $C^{n-1}$ by the Hadamard
determinantal inequality.
Now if $\bB=[b_{ij}] \in \mR^{n \times n}$ is such that
the absolute value of each entry is bounded above by $\alpha > 0$,
then $\| \bB \| \le n \alpha$, since each $L^2$-norm of the column of $\bB$ is bounded
by $\alpha \sqrt{n}$ and we apply \eqref{eq-normdef}.

Thus $|\sigma_1(\adj (\bA))| = | \sigma_1(\bA)...\sigma_{n-1}(\bA) |\le n \cdot C^{n-1}$, and the claim \eqref{eq-Aest}  follows.
\endproof

\medskip

\begin{cor} \label{cor-Uest}
Let $\{\vec{u}_{1}, \ldots , \vec{u}_n\} \subset \mR^n$ satisfy $\| \vec{u}_k \| \leq 2e_2$ for $1 \leq k \leq n$, and $| \bU | \geq \delta$. Then
\begin{align}\label{eq-Uest}
\|\bU^{-1} \| \leq n (2e_2)^{n-1}/| \bU | \leq n \cdot (2e_2)^{n-1}/\delta.
\end{align}
\end{cor}
The estimates \eqref{eq-messyest5} and \eqref{eq-Uest}  yield
\begin{eqnarray}
\|\bQ\| ~ = ~ \| \bW \cdot \bU^{-1} \| ~& \leq & ~ \|\bW \| \cdot \|\bU^{-1}\| \leq  \ve \cdot 2^n n^{3/2} (e_2)^{n-1}/\delta \label{eq-messyest7}.
\end{eqnarray}
and so  Lemma~\ref{lem-Vinvest} and Corollary~\ref{cor-Uest} imply:
\begin{cor}\label{cor-invest}
Assume that $ \ve < \delta/ 2^{n+1} n^{3/2} (e_2)^{n-1} $, then $\|\bQ\| < 1/2$ and so $\bV^{-1}$ exists.
Moreover, we have the estimate
$\| \bV^{-1} \| \leq n \cdot 2^n (e_2)^{n-1}/\delta$.
\end{cor}

We now return to the task of estimating $ \| \vec{\omega}(\vec{z}_0, \ldots, \vec{z}_{n}) - \vec{\omega}(\vec{y}_0, \ldots, \vec{y}_{n}) \|$ which will follow from an estimate the remaining terms in the equations \eqref{eq-solU} and \eqref{eq-solV}. Note that by \eqref{eq-messyest2}  and \eqref{eq-messyest3},
\begin{equation}\label{eq-messyest4}
\| \vl(\bV) \| ~ \leq ~ (2 e_3)^2 \sqrt{n} \quad , \quad \| \vl(\bV) - \vl(\bU) \| ~ \leq ~ 4 \ve (e_2 + e_3)\sqrt{n}.
\end{equation}

Using \eqref{eq-solU} and \eqref{eq-solV}, and the Taylor expansion of $(I+\bQ)^{-1}$, we then calculate
\begin{eqnarray}
2 ~\lefteqn{ \left \| \vec{\omega}(\vec{z}_0, \ldots, \vec{z}_{n}) - \vec{\omega}(\vec{y}_0, \ldots, \vec{y}_{n}) \right \| \quad } \nonumber\\
& \leq & ~ 2 \| \vec{z}_n - \vec{y}_n \| + \|\bU^{-1} \| \cdot \left\{ \| \vl(\bV) - \vl(\bU) \| ~ + ~ \| \vl(\bV) \| \cdot \|\bQ\|/(1 - \|\bQ\| ) \right\}. \label{eq-messyest1}
\end{eqnarray}

Assume that $ \ve < \delta/ 2^{n+1} n^{3/2} (e_2)^{n-1} $,
hence $\|\bQ\| < 1/2$ by Corollary~\ref{cor-invest} and so $\|\bQ\|/(1 - \|\bQ\| ) < 1$ and thus $\bV^{-1}$ exists. We use the more accurate estimate
$\|\bQ\| \leq \ve 2^n n^{3/2} (e_2)^{n-1}/\delta$ from \eqref{eq-messyest7} which, combined with
the previous estimates $\| \vec{y}_n - \vec{z}_n \| < \ve$, \eqref{eq-Uest} and \eqref{eq-messyest4}, then \eqref{eq-messyest1} becomes
\begin{eqnarray*}
\lefteqn{ 2 ~ \left \| \vec{\omega}(\vec{z}_0, \ldots, \vec{z}_{n}) - \vec{\omega}(\vec{y}_0, \ldots, \vec{y}_{n}) \right \| \quad } \\
& \leq & ~ 2 \ve + \left\{ n \cdot (2e_2)^{n-1}/\delta \right\} \cdot \left\{ 4 \ve (e_2 + e_3) \sqrt{n} ~ + ~ (2e_3)^2 \cdot 2 \ve \cdot 2^n n^{2} (e_2)^{n-1}/\delta \right\}.
\end{eqnarray*}
Then using that $e_3 = e_2 + \ve > e_2$ we have
\begin{equation}\label{eq-godawfulest}
\left \| \vec{\omega}(\vec{z}_0, \ldots, \vec{z}_{n}) - \vec{\omega}(\vec{y}_0, \ldots, \vec{y}_{n}) \right \| ~ < ~
\ve \cdot \left\{ 1 + n^{3/2} \, 2^{n+1} \, (e_3)^n/\delta ~ + ~ 2 \, n^{3} \, 2^{2n} \, (e_3)^{2n}/\delta^2 \right\}.
\end{equation}
It is important to note that the ratios $(e_3)^n/\delta$ and $ (e_3)^{2n}/\delta^2$ are ``dimensionless'', so the estimate \eqref{eq-godawfulest} is scale invariant, in that the expression in brackets on the right hand side is unchanged by scalar multiplication on $\mR^n$.

\subsection{Robustness of simplices}

 We next give an estimate for $\delta$, the constant in \eqref{eq-godawfulest}  which is a lower bound on   $|\bU|$, or equivalently on the volume of the simplex 
 $\Delta (\vec{y}_{0}, \ldots , \vec{y}_n )$. Since the edges of the simplex have lengths bounded by $2e_2$, this  condition guarantees that the vertex $\vec{y}_k$ in a simplex is not too close to a $k-1$-dimensional subspace defined by $\{\vec{y}_0,\ldots,\vec{y}_{k-1}\}$, and so ensures that a small perturbation of vertices does not change drastically the geometry of the simplicial complex.

\begin{defn}\label{def-robust1}
Let $\rho > 0$ and $1 \leq m \leq n$. A collection of vectors $\{\vec{y}_{0}, \ldots , \vec{y}_m\} \subset \mR^n$ is said to be \emph{$\rho$-robust} if for each $0 \leq k < m$, the distance from the point $\vec{y}_{k+1}$ to the affine subspace spanned by the vertices
$\ds \{\vec{y}_{0}, \ldots , \vec{y}_{k}\}$
is at least $\rho$.
\end{defn}
The significance of this definition is seen from an elementary estimation, whose proof follows by induction and standard Euclidean geometry. Let $P (\vec{y}_{0}, \ldots , \vec{y}_n )$ denote the parallelepiped with edges $\vec{v}_i - \vec{v}_0$ for $1 \leq i \leq n$, and note that its volume is equal to 
$n! \cdot {\rm Vol} (\Delta (\vec{y}_{0}, \ldots , \vec{y}_n ))$. 

\begin{lemma} \label{lem-volest1}
Let  $\{\vec{y}_{0}, \ldots , \vec{y}_n\} \subset \mR^n$ be a  $\rho$-robust collection, then  $P (\vec{y}_{0}, \ldots , \vec{y}_n )$
has volume at least $ \rho^{n-1} \cdot \|\vec{y}_1 - \vec{y}_0\|$. \hfill $\Box$
\end{lemma}
This   estimate can be improved when the vertices are lattice points on a circumscribed sphere:
\begin{lemma} \label{lem-volest2}
For $0 < e_1 < e_2$, there exists $V_2(e_1 , e_2) > 0$ such that given $\{\vec{y}_{0}, \ldots , \vec{y}_n\} \subset \mR^n$, and $0 < r \leq e_2$ satisfying:
\begin{enumerate}
\item $e_1 \leq \| \vec{y}_{k} - \vec{y}_{j} \|$ for $0 \leq j \ne k \leq n$,
\item $\| \vec{y}_{k} \| = r$ for all $0 \leq k \leq n$,
\item $\{\vec{y}_{0}, \ldots , \vec{y}_n\}$ is $\rho$-robust.
\end{enumerate}
Then  
$P(\vec{y}_{0}, \vec{y}_1 , \ldots , \vec{y}_n)$
has volume at least $V_2(e_1 , e_2) \cdot \rho^{n-2}$.
\end{lemma}
\proof
First, note that the vectors $\{\vec{y}_{0}, \vec{y}_{1} , \vec{y}_2\} \subset \mR^n$ cannot be collinear, as they lie on a sphere of radius $r \leq e_2$. Also, the vectors $\vec{\sigma}_1 = \vec{y}_1 - \vec{y}_0$ and $\vec{\sigma}_2 = \vec{y}_2 - \vec{y}_0$ have lengths greater than $e_1$ by (\ref{lem-volest2}.1), and thus define a non-degenerate parallelogram $P(\vec{y}_{0}, \vec{y}_1 , \vec{y}_2)$. The minimum for the area over all such parallelograms must be positive, as these conditions define a compact set of such, all of which have positive area. Let $V_2(e_1 , e_2) > 0$ denote this minimum.

Next, the vector $ \vec{y}_3$ lies at distance at least $\rho$ from the plane spanned by $\{\vec{y}_{0}, \vec{y}_1 , \vec{y}_2\}$ by the $\rho$-robust assumption. As $\vec{y}_0$ lies on this plane, $\vec{\sigma}_3 = \vec{y}_3 - \vec{y}_0$ must also lie distance at least $\rho$ from it.
Thus, $P(\vec{y}_{0}, \vec{y}_1 , \vec{y}_2, \vec{y}_3)$ with edges by $\{\vec{\sigma}_{1}, \vec{\sigma}_2 , \vec{\sigma}_3\}$ has $3$-volume bounded below by $V_2(e_1 , e_2) \cdot \rho$.

Continuing by induction, one has that the parallelepiped
$P(\vec{y}_{0}, \vec{y}_1 , \ldots , \vec{y}_{k})$ with edges $\{\vec{\sigma}_{1}, \ldots , \vec{\sigma}_{k}\}$
has $k$-volume bounded below by $V_2(e_1 , e_2) \cdot \rho^{k-2}$ for all $2< k \leq n$.  
\endproof

Lemma~\ref{lem-volest2} hints at a fundamental difference between the study of Delaunay triangulations in dimension $2$, and the theory for dimensions greater than two. The volume estimate for   simplices in dimension two admits a uniform lower positive bound depending only on the constants $0 < e_1 < e_2$. For higher dimensions, there is an additional restriction required to obtain   an estimate, the \emph{robustness} of the vertices, or some equivalent version of this condition. 
For example, if bounds are given on the interior angles of the simplex, then   this observation is surely well known.

We combine the above results to obtain the final form \eqref{eq-godawfulest2} of the desired estimate:
\begin{prop}\label{prop-varestimate}
Let $\{\vec{y}_{0}, \ldots , \vec{y}_n\} \subset \mR^n$ be $\rho > 0$ robust, and admit a circumscribed sphere with
center $ \omega(\vec{y}_0, \ldots , \vec{y}_n) $ and radius $r(\vec{y}_0, \ldots , \vec{y}_n)$.
Given $0 < e_1 < e_2$, set $\delta = V_2(e_1 , e_2) \cdot \rho^{n-2}$, and let $\ve > 0$. Suppose that, in addition, we have:
\begin{enumerate}\itemsep4pt
\item $e_1 \leq \| \vec{y}_i - \vec{y}_j \|$ \quad for all $0 \leq i \ne j \leq n$,
\item $e_1/2 \leq r(\vec{y}_0, \ldots , \vec{y}_n) \leq e_2$,
\item $\ds \ve ~ \leq ~ \delta/ 2^{n+1} n^{3/2} (e_2)^{n-1} ~ \leq ~ \frac{V_2(e_1 , e_2) \cdot \rho^{n-2}}{ 2^{n+1} n^{3/2} (e_2)^{n-1}} $.
\end{enumerate}
Let $\{\vec{z}_{0}, \ldots , \vec{z}_n\} \subset \mR^n$ satisfy
\begin{enumerate}\setcounter{enumi}{3}
\item $\| \vec{y}_i - \vec{z}_i \| \leq \ve$ \quad for all $0 \leq i \leq n$,
\end{enumerate}
then $\{\vec{z}_{0}, \ldots , \vec{z}_n\}$ has a circumscribed sphere with center $ \omega(\vec{z}_0, \ldots , \vec{z}_n) $ so that for $e_3 = e_2 + \ve$ ,
\begin{equation}\label{eq-godawfulest2}
\left \| \vec{\omega}(\vec{z}_0, \ldots, \vec{z}_{n}) - \vec{\omega}(\vec{y}_0, \ldots, \vec{y}_{n}) \right \| ~ < ~
\ve \cdot \left\{ 1 + n^{3/2} \, 2^{n+1} \, (e_3)^n/\delta ~ + ~ 2 \, n^{3} \, 2^{2n} \, (e_3)^{2n}/\delta^2 \right\}.
\end{equation}
\end{prop}

 \bigskip

\section{Circumscribed spheres via inequalities} \label{sec-inequalities}

We develop   an alternative approach to deriving the equations  of a circumscribed sphere for a given collection of  points
$\{\vec{z}_{0}, \ldots , \vec{z}_n\}$. The method assumes that  a system of inequalities is given, which defines an ``approximate solution'', and that there is a perturbation to an actual solution as described in the last section. This approach is advantageous when considering perturbations of a given triangulation, and we develop some key estimates which are used in later sections.

\subsection{Approximating centers of circumscribed spheres}

Given vectors $\{\vec{z}_{0}, \ldots , \vec{z}_n\} \subset \mR^n$, let $\bV$ denote the $n \times n$ matrix whose rows are the transposes of the vectors $\vec{v}_{k} = \vec{z}_{k-1} - \vec{z}_n$ for $1 \leq \ell \leq n$, and set $\vl(\bV) = (\|\vec{v}_{1}\|^2, \ldots , \|\vec{v}_{n}\|^2)^t$. Assuming that $\bV$ is invertible, the first result gives an estimate on the distance between an approximate center for the points and the actual center.
\begin{prop}\label{prop-varestimate2}
Suppose that we are given vectors
$\{\vec{z}_{0}, \ldots , \vec{z}_n\} \subset \mR^n$, $\omega \in \mR^n$ and constants $0 < C_1 < r$ and $C_2 > 0$ such that
  $\ds r - C_1 < \| \vec{z}_k - \omega \| < r + C_1 $   for all $0 \leq k \leq n$, and 
  $\| \bV^{-1} \| \leq C_2$.
Then $\{\vec{z}_{0}, \ldots , \vec{z}_n\}$ has a circumscribed sphere with center $ \omega(\vec{z}_0, \ldots , \vec{z}_n)$
such that
\begin{equation}\label{eq-approxsol}
\| \omega - \omega(\vec{z}_0, \ldots , \vec{z}_n)\| < 2\sqrt{n} \cdot r \, C_1 \, C_2
\end{equation}
\end{prop}
\proof
 The center $ \omega(\vec{z}_0, \ldots , \vec{z}_n)$ lies in the common intersection of the hyperplanes
\begin{eqnarray*}
L(\vec{z}_{k-1}, \vec{z}_{n}) ~ & = & ~ \left \{ \vec{x} \in \mR^n \mid ( \vec{z}_{k-1} - \vec{z}_n) \bullet (\vec{x} - ( \vec{z}_{k-1} + \vec{z}_n)/2 ) = 0 \right \} \\
~ & = & ~ \{ \zeta + \vec{z}_n \mid \vec{v}_{k} \bullet \vec{\zeta} = 1/2 \cdot \|\vec{v}_{k}\|^2 \}
\end{eqnarray*}
where $\zeta = \vec{x} - \vec{z}_n$. Thus, the solution $\ds \vec{\omega}(\vec{z}_0, \ldots, \vec{z}_{n})$ of the matrix equation \eqref{eq-solV} 
is the center for a circumscribed sphere containing the points $\{\vec{z}_{0}, \ldots , \vec{z}_n\}$.
We   estimate $\|\omega - \omega(\vec{z}_0, \ldots , \vec{z}_n)\|$.

As $r - C_1 > 0$, the vector $\omega$ satisfies the inequalities
\begin{equation}\label{eq-approxineq}
(r - C_1)^2 < \| \vec{z}_k - \omega \|^2 < (r + C_1)^2
\end{equation}
Make the change of variables $ \vec{v}_k = \vec{z}_{k-1} - \vec{z}_n$ and $\zeta = \omega - \vec{z}_n$ and 
subtract the inequalities \eqref{eq-approxineq} for $k=n+1$ from those for $1 \leq k \leq n$. Using that $\vec{v}_{n+1} = \vec{z}_{n} - \vec{z}_n = \vec{0}$, and expanding and canceling terms then yields
$$
\begin{array}{ccccc}
- 4 r C_1 & ~ < ~ & (\vec{v}_k - \zeta) \bullet (\vec{v}_k - \zeta) - (\vec{v}_{n+1} - \zeta) \bullet (\vec{v}_{n+1} - \zeta) & ~ < ~ & 4 r C_1 \\
- 4 r C_1 & ~ < ~ & \vec{v}_k \bullet \vec{v}_k - 2 \vec{v}_k \bullet \zeta & ~ < ~ & 4 r C_1
\end{array}
$$
Conditions of Proposition (\ref{prop-varestimate2}) and the above implies that $\zeta = \omega - \vec{z}_n$ is a solution of the matrix inequality
\begin{equation}\label{eq-matrixineq1}
\bV \cdot \zeta - \frac{1}{2} \vl(\bV) \in B(0, 2\sqrt{n} \cdot r C_1),
\end{equation}
and using the equation \eqref{eq-solV} we obtain that 
$\omega' = \omega - \omega(\vec{z}_0, \ldots , \vec{z}_n)$ is a solution of the matrix inequality
\begin{equation}\label{eq-matrixineq3}
\bV \cdot \omega' \in B(0, 2\sqrt{n} \cdot r \, C_1)
\end{equation}
We are given that $\|\bV^{-1}\| \leq C_2$ hence we obtain the estimate \eqref{eq-approxsol}.
\endproof

\medskip

\subsection{Stability of Delaunay triangulations}

Stability of the Delaunay triangulation associated to a net $\cN \subset \mR^n$ under perturbation of $\cN$ is equivalent to the stability of the circumscribed spheres for the vertices of a simplex. The following result shows the existence of circumscribed spheres based on estimates which are almost ``stable under sufficiently small'' perturbation.

\begin{prop}\label{prop-varestimate3}
Let $\{\vec{z}_{0}, \ldots , \vec{z}_n\} \subset \mR^n$ be $\rho$-robust, for $\rho > 0$. Assume there are constants
$0 < e_1 < e_2$ and $0 < C_1 < r < e_1$, and  that there exists $\omega \in \mR^n$ such that
\begin{enumerate}\itemsep4pt
\item $e_1 < \| \vec{z}_i - \vec{z}_j \| < 2 e_2$ \quad for all $0 \leq i \ne j \leq n$
\item $r - C_1 < \| \vec{z}_k - \omega \| < r + C_1 $ \quad for all $0 \leq k \leq n$,
\end{enumerate}
Then $\{\vec{z}_{0}, \ldots , \vec{z}_n\}$ has a circumscribed sphere with center $ \omega(\vec{z}_0, \ldots , \vec{z}_n) $ so that for
\begin{equation} \label{eq-varestimate3}
\left \| \omega - \vec{\omega}(\vec{z}_0, \ldots, \vec{z}_{n}) \right \| ~ \leq ~ C_1 \cdot n^{3/2} (2e_2)^{n-1}/\rho^{n-1}
\end{equation}
\end{prop}
\proof
Lemma~\ref{lem-volest1} implies that the volume of the parallelepiped
$P(\vec{z}_{0}, \ldots , \vec{z}_{n})$ with edges $\{\vec{v}_{1}, \ldots , \vec{v}_{n}\}$ is bounded below by $e_1 \rho^{n-1}$, and hence $|\bV | \geq e_1 \rho^{n-1}$.
Thus by Corollary~\ref{cor-Uest}, we have
\begin{align}\label{eq-Vest3}
\|\bV^{-1} \| \leq n (2e_2)^{n-1}/| \bV | \leq n \cdot (2e_2)^{n-1}/e_1 \rho^{n-1}
\end{align}
Then \eqref{eq-varestimate3}  follows from estimate \eqref{eq-approxsol} of Proposition~\ref{prop-varestimate2} and the hypotheses $r \leq e_1$.
\endproof

Propositions~\ref{prop-varestimate} and \ref{prop-varestimate3} show the importance of the robustness condition in Definition~\ref{def-robust1} for estimating the stability of solutions for the equations \eqref{eq-solU}. Our next result shows that a small perturbation of a robust simplex is also robust.

\begin{prop}\label{prop-robustvarest}
Let $1 \leq m \leq n$, and
assume that $\{\vec{y}_{0}, \ldots , \vec{y}_m\} \subset \mR^n$ is $\rho$-robust. Let $\{\vec{z}_{0}, \ldots , \vec{z}_m\} \subset \mR^n$ be also given, along with the constants
$0 < e_1 < e_2$ and $0 < \ve < e_1/4$ such that
\begin{enumerate}\itemsep4pt
\item $e_1 \leq \| \vec{y}_i - \vec{y}_j \| \leq 2 e_2$ ~ for all $0 \leq i \ne j \leq m$
\item $\| \vec{y}_i - \vec{z}_i \| \leq \ve$ ~ for all $0 \leq i \leq m$.
\end{enumerate}
Then $\{\vec{z}_{0}, \ldots , \vec{z}_m\}$ is $\rho_m$-robust, for
$\rho_m = \rho_m(\rho, \ve, e_1, e_2)$ as defined below. Moreover,
$\rho_m(\rho, \ve, e_1, e_2)$ is monotone increasing in $e_2$ and $\rho$, and monotone decreasing in $e_1$ and $\ve$, and is scale-invariant. That is, for $s > 0$,
$\rho_m(s \cdot \rho, s \cdot \ve, s \cdot e_1, s \cdot e_2)= s \cdot \rho_m(\rho, \ve, e_1, e_2)$.
\end{prop}
\proof
Set $e_1' = e_1 - 2\ve$, $e_2' = e_2 + \ve$ and $e_4 = 4(e_2 + e_1)$. Then for all $0 \leq i \ne j \leq m$,
$$e_1/2 < e_1' < \| \vec{z}_i - \vec{z}_j \| < 2 e_2' < e_4$$
For each $0 \leq k \leq m$, let ${\rm Span}(\vec{y}_0, \ldots, \vec{y}_k) \subset \mR^n$ denote the affine subspace spanned by the vectors, and let $\xi_{k} \in {\rm Span}(\vec{y}_0, \ldots, \vec{y}_{k-1})$ be the point closest to $\vec{y}_k$.
Then $\rho \leq \|\vec{y}_k - \xi_{k}\| \leq \|\vec{y}_k - \vec{y}_{0}\| \leq 2e_2$.

Similarly, let ${\rm Span}(\vec{z}_0, \ldots, \vec{z}_{k-1}) \subset \mR^n$ denote the affine subspace spanned by the vectors, and
$\zeta_{k} \in {\rm Span}(\vec{z}_0, \ldots, \vec{z}_{k-1})$ be the point closest to $\vec{z}_k$.
Then $\|\vec{z}_k - \zeta_{k}\| \leq \|\vec{z}_k - \vec{z}_{j}\| < 2e_2'$ for $ j \leq k-1$.

The triangle inequality yields a lower bound
\begin{eqnarray}
d_{\mR^m}\left( \vec{z}_k , {\rm Span}(\vec{z}_0, \ldots, \vec{z}_{k-1}) \right) = \| \vec{z}_k - \zeta_{k} \|
& ~ \geq ~ & \| \vec{y}_k - \xi_k \| ~ - ~ \| \vec{z}_k - \vec{y}_k \| ~ - ~ \|\xi_k - \zeta_k \| \nonumber \\
& ~ \geq ~ & \rho ~ - ~ \ve ~ - ~ \|\xi_k - \zeta_k \| \label{eq-rhoest}
\end{eqnarray}
We  develop an upper bound estimate for $\| \xi_k - \zeta_k \|$.

For the case $k =1$,
note that ${\rm Span}(\vec{z}_0) = \{\vec{z}_0\}$ is just the single point, so $\xi_1 = \vec{y}_0$ and $\zeta_1 = \vec{z}_0$, and
$\|\zeta_1 - \xi_1 \| = \|\vec{z}_0 - \vec{y}_0\| \leq \ve$, so in terms of the estimate \eqref{eq-rhoest} we have
$\ds d_{\mR^m}\left( \vec{z}_1 , {\rm Span}(\vec{z}_0) \right) ~ \geq ~ \rho - 2 \ve$.
Set $\delta_1 = 2$, then $\rho_1 = \rho - \ve \cdot \delta_1$.
This completes the proof of  Proposition~\ref{prop-robustvarest} for the case $m = 1$.

When $m > 1$ and  $2 \leq k \leq m$, an  upper bound estimate on $\| \xi_k - \zeta_k \|$ requires more delicate arguments. 

We are given that $\vec{y}_j , \vec{z}_j \in D_{\mR^n}(\vec{y}_k , 2e_2 + \ve)$ for each $0 \leq j \leq m$.
Since the distance from $\vec{y}_k$ to $\xi_k$ is at most that from $\vec{y}_k$ to $\vec{y}_0$ we also have
$\xi_k \in D_{\mR^n}(\vec{y}_k , 2e_2)$. The analogous estimate is true for $d_{\mR^n}(\vec{z}_k  , \zeta_k)$,
 and since $\|\vec{y}_k - \vec{z}_k \| \leq \ve$ we have that $\zeta_k \in D_{\mR^n}(\vec{y}_k , 2e_2')$. 
 It follows that all of the points in consideration, $\vec{y}_j,\vec{z}_j,\xi_j,\zeta_j$, $1 \leq j \leq k$, lie in the closed disk $D_{\mR^n}(\vec{y}_k , e_4)$ with radius $e_4 = 4(e_2 + e_1)$. This compactness estimate is fundamental.

Let $\ds {\rm Span}_k(\vec{y}_0, \ldots, \vec{y}_{k-1}) ~ = ~ 
{\rm Span}(\vec{y}_0, \ldots, \vec{y}_{k-1}) ~ \cap ~ D_{\mR^n}(\vec{y}_k , 2e_2')$, the restricted subdisk of radius $2e_2'$. 
Note that we showed above that $\{\vec{y}_0, \ldots, \vec{y}_{k-1}, \xi_1, \ldots, \xi_k\} \subset {\rm Span}_k(\vec{y}_0, \ldots, \vec{y}_{k-1})$.

For the case $k =2$, note that $\| \vec{y_1} - \vec{y}_0 \| \geq e_1$ and $\| \vec{z_1} - \vec{z}_0 \| \geq e_1' > e_1/2$, and using that the disk $D_{\mR^n}(\vec{y}_2 , 2e_2')$ has diameter at most $e_4$, we have
\begin{eqnarray}
{\rm Span}_2(\vec{y}_0, \vec{y}_1) & \subset & \{ \vec{y}_0 + t_1 (\vec{y_1} - \vec{y}_0) \mid - e_4/e_1 \leq t_1 \leq e_4/e_1 \} \label{eq-spany2}\\
{\rm Span}_2(\vec{z}_0, \vec{z}_1) & \subset & \{ \vec{z}_0 + s_1 (\vec{z_1} - \vec{z}_0) \mid - e_4/e_1' \leq s_1 \leq e_4/e_1' \} \label{eq-spanz2}
\end{eqnarray}

\begin{lemma}\label{lem-distanceest2}
Given $\vec{z} \in {\rm Span}_2(\vec{z}_0, \vec{z}_1) $, there exists $\vec{y} \in {\rm Span}(\vec{y}_0, \vec{y}_1)$ so that
\begin{equation}\label{eq-distanceest2}
\| \vec{z} - \vec{y} \| \leq \ve \cdot (1 + 4e_4/e_1)
\end{equation}
\end{lemma}
\proof
Write down $\vec{z} \in {\rm Span}_2(\vec{z}_0, \vec{z}_1) $ as
$\vec{z} = \vec{z}_0 + s_1 \cdot (\vec{z_1} - \vec{z}_0) \in {\rm Span}_2(\vec{z}_0, \vec{z}_1)$, then for
$\vec{y} = \vec{y}_0 + s_1 \cdot (\vec{y_1} - \vec{y}_0) \in {\rm Span}_2(\vec{y}_0, \vec{y}_1)$ we have
$\ds  \| \vec{z} - \vec{y} \|  \leq \ve \cdot (1 + 4e_4/e_1)$. 
Thus, every point of ${\rm Span}_2(\vec{z}_0, \vec{z}_1) $ has distance at most $\ve \cdot (1 + 4e_4/e_1)$ from a point of
$ {\rm Span}(\vec{y}_0, \vec{y}_1)$.
\endproof

Lemma~\ref{lem-distanceest2} implies that
$\|\xi_2 - \zeta_2\| \leq \ve \cdot (1 + 4e_4/e_1)$, hence $\|\vec{z}_2 - \zeta_2 \| \geq \rho_2$ by \eqref{eq-rhoest}, where
\begin{equation}\label{eq-delta2}
\rho_2 = \rho - \ve \cdot (2 + 4 e_4/e_1) = \rho - \ve \cdot \delta_2(\rho, e_1, e_2)
\end{equation}
Note that $\delta_2(\rho, e_1, e_2) = (2 + 4 e_4/e_1)$ depends only on the constants $e_1$ and $e_2$, and as the ratio $e_4/e_1$ is scale invariant, thus $\rho_2$ is also scale invariant. If $m =2$ then we are done.

Next, consider the case $k=3$. The estimate $\rho_3$ in this case is obtained from \eqref{eq-rhoest}  by subtracting from $\rho$ a term which involves linear combinations of $\vec{y}_2$ with points of the line ${\rm Span}(\vec{y}_0, \vec{y}_1)$, and the closer that $\vec{y}_2$ lies to this line, the larger the possible error, and likewise for ${\rm Span}_3(\vec{z}_0, \vec{z}_1, \vec{z}_2)$.

As seen before for $k=2$, the strategy is to estimate the parameters used to describe the planar region ${\rm Span}_3(\vec{y}_0, \vec{y}_1, \vec{y}_2)$ as in \eqref{eq-spany2},
and similarly for ${\rm Span}_3(\vec{z}_0, \vec{z}_1, \vec{z}_2)$ as in \eqref{eq-spanz2}.

Recall that $\xi_2 \in {\rm Span}(\vec{y}_0, \vec{y}_1)$ is the point on the line closest to
$\vec{y_2}$, and $\rho \leq \| \vec{y}_2 - \xi_2\| \leq 2e_2 < e_4$.

Likewise, the point $\zeta_2 \in {\rm Span}(\vec{z}_0, \vec{z}_1)$ closest to $\vec{z}_2$
satisfies $\rho_2 \leq \|\vec{z}_2 - \zeta_2 \| \leq 2e_2' <e_4$.

Now let $\xi_2' \in {\rm Span}(\vec{y}_0, \vec{y}_1)$ be the point closest to $\zeta_2$.
Then $\|\vec{y_2} - \xi_2'\| \geq \|\vec{y_2} - \xi_2\| \geq \rho > \rho_2$.

Furthermore, from the case $k=2$, we have that $\|\xi_2' - \zeta_2\| \leq \ve \cdot \delta_2(\rho, e_1, e_2)$.

The key idea is to bound the space ${\rm Span}_3(\vec{y}_0, \vec{y}_1, \vec{y}_2)$ using linear combinations with $(\vec{y_2} - \xi_2')$ and parameter bounds invoking $\rho$ and $\rho_2$:
\begin{eqnarray*}
{\rm Span}_3(\vec{y}_0, \vec{y}_1, \vec{y}_2) & \subset & \{ \vec{y}_0 + t_1 (\vec{y_1} - \vec{y}_0) + t_2 (\vec{y_2} - \xi_2') \mid - e_4/e_1 \leq t_1 \leq e_4/e_1 , - e_4/\rho \leq t_2 \leq e_4/\rho \} \label{eq-spany3}\\
{\rm Span}_3(\vec{z}_0, \vec{z}_1, \vec{z}_2) & \subset & \{ \vec{z}_0 + s_1 (\vec{z_1} - \vec{z}_0) + s_2 (\vec{z_2} - \zeta_2) \mid - e_4/e_1' \leq s_1 \leq e_4/e_1', - e_4/\rho_2 \leq s_2 \leq e_4/\rho_2 \} \label{eq-spanz3}
\end{eqnarray*}
As in the proof of Lemma~\ref{lem-distanceest2}, every point of ${\rm Span}_3(\vec{z}_0, \vec{z}_1, \vec{z}_2) $ thus lies a distance at most
$$ \ve \cdot \left\{ 1 + 2e_4/e_1 \cdot (1+1) + 2 e_4/\rho_2 \cdot (1 + \delta_2 ) \right\} $$
from a point of $ {\rm Span}_3(\vec{y}_0, \vec{y}_1, \vec{y}_2)$, and in particular this estimate holds for $\|\xi_3 - \zeta_3\|$.
Set
\begin{equation}\label{eq-delta3}
\delta_3 = \delta_3(\rho, e_1, e_2) = 2 + 4 e_4/e_1 + (1 + \delta_2 ) \cdot 2e_4/\rho_2
\end{equation}
Note that the ratio $e_4/\rho_2$ is scale-invariant, as is $\delta_2 $, and thus $\delta_3$ is scale-invariant.

Then for $\rho_3 = \rho - \ve \cdot \delta_3 $ by \eqref{eq-rhoest} we have $ \|\vec{z}_3 - \zeta_3 \| \geq \rho_3$.

Continuing in the way, given $\rho_k$ and $\delta_k$ for $2 \leq k < m$, define inductively
\begin{eqnarray}
\delta_{k+1} & ~ = ~ & 1 + \left\{ 1 + 2 \cdot 2e_4/e_1 + (1 + \delta_2 ) \cdot 2e_4/\rho_2 + \cdots + (1 + \delta_k ) \cdot 2 e_4/\rho_k \right\} \label{eq-deltak} \\
\rho_{k+1} & ~ = ~ & \rho - \ve \cdot \delta_{k+1} \label{eq-robust-k}
\end{eqnarray}
Then we have
$ \|\vec{z}_{k+1} - \zeta_{k+1} \| \geq \rho_{k+1}$. Continuing until $k+1 = m$, we obtain
\begin{eqnarray}
\delta_m(\rho, e_1, e_2) & ~ = ~ &2 + 4 e_4/e_1 + 2 \cdot \sum_{k=2}^{m-1} ~ \frac{(1 + \delta_k) e_4}{\rho_k} \label{eq-delta-k=m}\\
\rho_m(\rho, \ve, e_1, e_2) & ~ = ~ & \rho - \ve \cdot \delta_m(\rho, e_1, e_2) \label{eq-robust-k=m}
\end{eqnarray}
for which $\ds d_{\mR^m}\left( \vec{z}_m , {\rm Span}(\vec{z}_0, \ldots , \vec{z}_{m-1}) \right) = \|\vec{z}_m - \zeta_m \| \geq \rho_m(\rho, \ve, e_1, e_2)$.

Observe that by the inductive definition \eqref{eq-robust-k}, the values $\rho> \rho_1 > \cdots > \rho_m$ are monotone decreasing.
Furthermore, by an inductive argument, for each $1 \leq k < m$ the value of $\rho_k$ is a monotone increasing function of $e_2$ and $\rho$, and monotone decreasing for $e_1$, and thus each term
$(1 + \delta_k)e_4/\rho_k$ in the sum \eqref{eq-delta-k=m} is also monotone increasing, hence the same holds for
$\rho_m(\rho, \ve, e_1, e_2)$.
Also note that each additional term $(1 + \delta_k ) \cdot 2e_4/\rho_k$ in \eqref{eq-deltak}  is scale-invariant, so the sum \eqref{eq-robust-k=m}  is scale-invariant.
\endproof

\section{Micro-local foliation geometry} \label{sec-approx}

The construction of the Vononoi cells in the sections above uses the distance function on $\mR^n$ to derive the linear equations which define the circumscribed spheres that define the Delaunay triangulation. The extensions of these ideas to leaves of foliations requires working with the given Riemannian metric on the leaves, as these define the leafwise distance functions. If the foliation $\F$ is defined by a free action of $\mR^n$, such as for the case of a tiling space associated to a tiling, then there is a natural Euclidean metric on each leaf. 
However, when the leaves are just assumed to be smooth Riemannian manifolds, then the construction of the Delaunay triangulation associated to an arbitrary  leafwise net is problematic, as discussed in \cite{LeibonLetscher2000}, for example. 

Our approach is to introduce, given a matchbox manifold $\fM$ with smooth leafwise Riemannian metric,  a ``distance scale'' which is sufficiently small so that the local coordinate charts for this scale are ``almost Euclidean'', and then choose the leafwise net to be adapted to this scale. The definition of the Delaunay simplices will then be obtained by adapting the methods of the previous sections to this almost non-Euclidean context. 

This study of local properties of a Riemannian manifold is best done in adapted geodesic coordinates, and this requires the introduction of some standard ideas of Riemannian geometry such as local orthonormal frames and their Christoffel symbols. A good reference for this material is \cite[\S 9]{Helgason1978}. We recall the necessary material below, which leads to the   choices of constants in the next section, and the corresponding estimates in later sections which guarantee the stability of the simplices in the leafwise Delaunay triangulation.

For each $1 \leq i \leq \nu$,  there is given the coordinate chart $\vp_i \colon \oU_i \to [-1, 1]^n \times \fT_i$ and 
for $x \in \oU_i$,  the plaque for the chart $\vp_i$ containing $x$ is denoted by  $\cP_i(x)$. 
The horizontal coordinate function $\lambda_i \colon \oU_i \to [-1,1]^n$ is defined by setting $\vp_i(x) = (\lambda_i(x) , w_x) \in [-1,1]^n \times \fT_{i}$.
Also, recall that $\lF > 0$ was chosen in Lemma~\ref{lem-stronglyconvex} so that for all $x \in \fM$, the closed leafwise disk $D_{\F}(x, \lF)$ is strongly convex, and     $2\dFU < \lF/2$ bounds the diameter of the plaques in the foliation covering.

For $x \in \oU_i$ define the transversal section, for $ \fU \subset \fT_i$
\begin{equation}\label{eq-transsec1}
\fZ(x, i, \fU) \equiv \vp_i^{-1}\left(\lambda_i(x) , \fU \right) ~ ; ~ \fZ(x, i) \equiv \vp_i^{-1}\left(\lambda_i(x) , \fT_i \right) = \lambda_i^{-1} \circ \lambda_i(x)
\end{equation}
As a special case, for $r \geq 0$, define the compact ``disk section''
\begin{equation}\label{eq-transsec2}
\fZ(x, i, r) \equiv \vp_{i}^{-1}\left(\lambda_i(x) , D_{\fX}(w_x, r) \cap \fT_i\right) \subset \oU_i
\end{equation}
The local coordinate charts $\vp_{i} \colon \oU_i \to [-1,1]^n \times \fT_{i}$ are used to define a local ``vertical translation'' between plaques, which will be fundamental in the following.
For $x' \in \fZ(x, i)$, define
\begin{equation}\label{eq-transvar1}
\phi_i(x,x') \colon \cP_i(x) \to \cP_i(x') ~, ~ \xi' = \phi(x,x')(\xi) = \fZ(\xi, i ) \cap \cP_i(x')
\end{equation}
When expressed in coordinates, 
\begin{equation}\label{eq-coordsphi}
\vp_i \circ \phi_i (x,x') \circ  \vp_i^{-1}(\lambda_i(x), w_x) = (\lambda_i(x), w_{x'})
\end{equation}
which is just the constant map in the first coordinate.
Thus $\phi_i(x',x'') \circ \phi_i(x,x') = \phi_i(x,x'')$, and
the maps $\phi_i (x,x') $ are homeomorphisms which depend continuously on $x' \in \fT_i$ in the $C^0$-topology.

\subsection{Leafwise metric distortions}
The leafwise Riemannian metric on$\F$ depends continuously on the transverse coordinate in local charts. We use this to define    distortion estimates for the metric.

First, we introduce estimates on the \emph{leafwise metric} distortions of the maps $\phi_i(x,x')$, which 
 compare the leafwise Riemannian distance functions induced on differing plaques in the same chart $\oU_i$. Set:
\begin{eqnarray}
var(i, r) & = & \max \left\{ |\dF(x,y) - \dF(x', y')| \mid x \in \oU_i , ~ x' \in \fZ(x, i, r), ~ y \in \cP_i(x), ~ y' = \phi(x,x')(y) \right\} \nonumber \\
& = & \max \left\{ \{ |\dF(y,z) - \dF(\phi_i(x,x')(y), \phi_i(x,x')(z))| \} \mid y, z \in \cP_i(x) , ~ x' \in \fZ(x, i, r) \right\} \label{eq-var}
\end{eqnarray}
Note that $var(i,r)$ depends continuously on $r$, that $var(i,0) = 0$, and $var(i, r) \leq 2\dFU$ as $\cP_i(x)$ is contained in a disk in $L_x$ of radius $\dFU$.

There is another measure of the metric distortion between plaques, this time in terms of the variation due to differing coordinate systems.
For $z \in \oU_i \cap \oU_j$  we obtain two standard transversals $\fZ(z, i, r)$ and $\fZ(z, j, r)$ in $\fM$ through $z$. 
Define the \emph{divergence} between these two transversals by
\begin{equation}\label{eq-divpairwise}
div(z,i,j,r) = \max \left\{ \dF(x',y') \mid x' \in \fZ(z, i, r) ~, ~ y' \in \fZ(z, j, r) ~, ~ \cP_i(x') \cap \cP_j(y') \ne \emptyset \right\}
\end{equation}
The assumption $\dFU < \lF/4$ implies that $div(z,i,j,r) < \lF$.
Note that $div(z,i,i,r) = 0$ and that $div(z,i,j,0) = 0$. Define: 
\begin{equation}\label{eq-divpairwise2}
div(z, r) = \max \left\{ div(z,i,j,r) \mid z \in \oU_i \cap \oU_j \right\}
\end{equation}
The condition $ \cP_i(x') \cap \cP_j(y') \ne \emptyset$ is closed in $x',y'$, and hence $div(z,r)$ is an upper semi-continuous function of both $z$ and $r$.
In terms of the transverse translation maps $\phi_i$, for $\ve = div(z,r)$, the condition \eqref{eq-divpairwise2}  implies that the compositions $\phi_i(x',z) \circ \phi_j(z,y')$ are $\ve$-close to the identity.

\subsection{Adapted geodesic coordinate systems}
If $\F$ is defined by an isometric free action of $\mR^n$, as in the case of tiling spaces for example, then the leaves of $\F$ have natural isometric coordinate systems. However, for  a manifold of non-zero curvature, it is necessary to introduce geodesic coordinates based at the points of $\fM$.  The books  \cite{BC1964} and \cite{doCarmo1992} are suitable references.

Let $\whe \equiv \{\vec{e}_1 , \ldots, \vec{e}_n \}$ denote the standard orthonormal basis of $\mR^n$.
A point $\vec{x} \in \mR^n$ is then written in coordinates as $\vec{a} = (a_1 , \ldots , a_n)$, where $\vec{x} = \whe \cdot \vec{a} = a_1 \vec{e}_1 + \cdots + a_n \vec{e}_n$.
Recall that the closed ball of radius $\lambda$ about the origin in the standard metric is denoted by  $D(\lambda)$, or $D_{\mR^n}(\lambda)$ when it is better to emphasize that the   disk is defined using the standard norm $\| \cdot \|_{\mR^n}$.

For $x \in \fM$, and coordinate system $\vp_i$ with $x \in U_i$ the basis $\whe$ of $\mR^n$ defines a framing
$\whe_w$ of $T_{\vec{0}} (-1,1)^n \times \{w\}$ for each $w \in \fT_i$.  For each $x \in U_i$, the differential  of the coordinate map $\vp_i$ at $x$ defines  a    linear isomorphism $d_x\vp_i \colon T_x\F \cong \mR^n$, by which  $\whe_w$ induces a framing $\whe_x$ for $T_x\F$.
If the curvature of leaves is non-zero, then  the  tangent map  $d_x\vp_i$  is typically not   an isometry, and thus 
 the framing $\whe_x$ is typically not orthonormal for the leafwise Riemannian metric.

The leafwise Riemannian metric on $T\F$ induces on each plaque $\cP_i(w) = \varphi_i^{-1}((-1,1)^n \times \{w\})$ of $U_i$ a family of inner products on its tangent space, which    in terms of the framing $\whe_x$ at $x \in \cP_i(w)$ is denoted by the matrix $g_{jk}(x)$. By Theorem~\ref{thm-riemannian}, the tensor $g_{jk}(x)$ varies continuously in $w \in \fT_i$ for the $C^{\infty}$-topology on functions on $\cP_i(w)$.

Given an arbitrary orthonormal frame $\whu = \{\vec{u}_1, \ldots, \vec{u}_n\} \subset T_x\F$ for the leafwise Riemannian metric, define a linear isomorphism
\begin{equation}
F_{\whu} \colon \mR^n \to T_x\F \cong \mR^n ~ , ~ F_{\whu}(a_1, \ldots , a_n) = \whu \cdot \vec{a}
\end{equation}
where we  adopt the ``matrix notation'' $\ds \whu \cdot \vec{a} \equiv a_1 \vec{u}_1 + \cdots + a_n \vec{u}_n \in T_x\F$.
Via the coordinate isomorphism $d_x\vp_i$, the tangent vectors $\vec{u}_k$ form an orthonormal set 
$\wthu \subset \mR^n$ for the inner product $g_{jk}(x)$, and in this sense, $\whu \cdot \vec{a}$ is precisely a matrix product.
To simplify notation, we let $\whu = \wthu$ also denote this framing, as it is clear from context whether we consider the framing as in $T_x\F$ or in $\mR^n$.
Then $F_{\whu}$ is a linear isometry between $\{ \mR^n , \| \cdot \| \}$ and $\{\mR^n , \| \cdot \|_{\whu} \}$, where  $\| \cdot \|_{\whu}$ denotes the norm on $T_x\F \cong \mR^n$ induced by the inner product $g_{ i j}(x)$.

Recall that $\exp^{\F}_x \colon T_x\F \to L_x$ is the leafwise geodesic map at $x$. Given an orthonormal framing $\whu$ of $T_x\F$ and $0 < \lambda \leq \lF$, the \emph{leafwise geodesic coordinates} at $x$ are defined by
\begin{equation}\label{eq-normalcoords}
\psi^g_{x, \whu} \colon D_{\mR^n}(\lambda) \to D_{\F}(x, \lambda) \subset L_x ~, ~ \psi^g_{x, \whu}(\vec{a}) = \exp^{\F}_x(\whu \cdot \vec{a})
\end{equation}
Assume that $D_{\F}(x, \lambda) \subset U_i$, and let $\wtx = \lambda_i(x) \in (-1,1)^n \subset \mR^n$.
Then we have a second coordinate system on the neighborhood $D_{\F}(x, \lambda)$ of $x$, which is also ``adapted'' to the leafwise Riemannian metric on the disk $D_{\F}(x, \lambda)$.
Define $T_{\wtx} \colon \mR^n \to \mR^n$ by $T_{\wtx}(\vec{y}) =  \wtx + \vec{y}$, and compose $T_{\wtx}$ with the framing map $F_{\whu}$ to obtain:
\begin{equation}
\psi^i_{x, \whu} \equiv \vp_i^{-1}(T_{\wtx} \circ F_{\whu} , w_x) \colon D_{\mR^n}(\lambda) \to \cP_i(x) ~ , ~ \psi^i_{x, \whu}(\vec{y}) = \vp_i^{-1}(\wtx + \whu \cdot \vec{y} , w_x)
\end{equation}
Then $\psi^i_{x, \whu}$  is the    geodesic coordinate system for the \emph{flat metric} on $\cP_i(x)$ associated to $\| \cdot \|_{\whu}$.

\subsection{Comparison of geodesic coordinate systems}\label{subsubsec-compargeod}

We compare the \emph{affine geometries} defined by these two sets of coordinates, $\psi^g_{x, \whu}$ and $\psi^i_{x, \whu}$, using the coordinate system $\vp_i$ to convert the comparison  to a local problem on $\mR^n$ involving    differential equations on $\mR^n$.

Recall that $D(\lambda) = D_{\mR^n}(\lambda)$ is the Euclidean disk centered at the origin,  $\wtx = \lambda_i(x)$, and  $D_{\wtg}(\wtx, s)$ denotes the closed disk of radius $s$ about $\wtx$ for the metric $\wtg$.

Let $\wtD_i(\wtx , \lambda) = \vp_i(D_{\F}(x,\lambda)) \subset (-1,1)^n \times \{w_x\}$ 
denote the image of the   disk in the leafwise metric.

Let $\wtd$ denote the distance function on $\wtD_i(\wtx , \lambda)$ defined by the leafwise metric $\dF$.
That is, for $\vec{y}, \vec{z} \in \wtD_i(\wtx , \lambda)$,
$\ds \wtd(\, \vec{y}, \vec{z} \, ) ~ = ~ \dF (\vp_i^{-1}(\vec{y} ,w_x) , \vp_i^{-1}(\vec{z} ,w_x))$.

Let $\wtg$ denote the metric tensor on $\wtD_i(x, \lambda)$ in the coordinates $\vp_i$.
Note that the image under $\vp_i$ of a geodesic segment for $g$ is a geodesic segment for $\wtg$, and as $D_{\F}(x,\lambda)$ is strongly convex for $\lambda \leq \lF$, the same holds for the region $\wtD_i(\wtx , \lambda)$ with the metric $\wtg$.

Let $ \wtexp_{\wtx}$ denote the geodesic map associated to $\wtg$, centered at $\wtx$. Then for the orthonormal basis $\whu \subset T_x\F$ considered as a frame for $T_{\wtx}\mR^n$,  we set
\begin{equation}
\wtexp_{\wtx,\whu} \colon D(\lambda) \to \wtD_i(\wtx , \lambda) ~, ~ \wtexp_{\wtx,\whu}(\vec{a}) = \wtexp_{\wtx}(\whu \cdot \vec{a})
\end{equation}
Recall that we also have a linear map $T_{\wtx} \circ F_{\whu}$, which is a linear isometry between $\{ \mR^n , \| \cdot \| \}$ and $\{\mR^n , \| \cdot \|_{\whu} \}$, and satisfies $T_{\wtx} \circ F_{\whu}(\vec{0}) = \wtx = \wtexp_{\wtx,\whu}(\vec{0})$.
Let $g^{\whu} = (T_{\wtx} \circ F_{\whu})^*(\wtg)$ denote the metric $\wtg$ near $\wtx$ pulled back to $D(\lambda)$ via the isometry $T_{\wtx} \circ F_{\whu}$. Then $g_{j k}^{\whu}(\vec{a}) = \delta_{j k}$ for $\vec{a} = \vec{0}$ by definition of $\whu$. 
The metric $g^{\whu}$ is in \emph{Gauss normal form}, and its metric tensor  $g_{j k}^{\whu}$ consequently has further special properties in a neighborhood of $\vec{0}$.

\begin{defn}\label{def-approxeuclid}
Let $x \in \fM$ and $0 < \lambda \leq \lF/2$. Assume that $D_{\F}(x, \lambda) \subset \cP_i(x)$, and let $\whu$ be
an orthonormal frame for $T_x\F$. For $\ve > 0$, we say that $\psi^g_{x, \whu} \colon D(\lambda) \to D_{\F}(x, \lambda)$ is \emph{$\ve$-approximately Euclidean} if the following hold (in the coordinate system $\vp_i$):
\begin{enumerate}
\item For all $\vec{a} \in D(\lambda)$,
\begin{equation}\label{eq-dgs0}
\| \, g_{j k}^{\whu}(\vec{a}) - \delta_{j k} \, \| ~ \leq ~ \ve/n^2
\end{equation}
\item
For all $\vec{a} \in D(\lambda)$,
\begin{equation}\label{eq-dgs1}
\wtd(\wtexp_{\wtx,\whu}(\vec{a}) , T_{\wtx} \circ F_{\whu}(\vec{a})) ~ \leq ~ \ve \cdot \| \vec{a}\|
\end{equation}
\item
For a geodesic $\wtsigma \colon [0,1] \to \wtD_i(\wtx , \lambda)$ in the $\wtd$ metric, with $\wtsigma(0) = \vec{y}_0$ and $\wtsigma(1) = \vec{y}_1$,
set $\wttau(t) = t \cdot ( \vec{y}_1 - \vec{y}_0) + \vec{y}_0$, then
\begin{equation}\label{eq-dgs2}
\wtd( \wtsigma(t) , \wttau(t) ) ~ \leq ~ \ve \cdot \wtd(\vec{y}_0 , \vec{y}_1) ~, ~ \text{for all } 0 \leq t \leq 1
\end{equation}
\item
For $s \leq \lambda$, the Riemannian volume of leafwise disks satisfies
\begin{equation}\label{eq-dgs3}
\left | {\rm Vol}(D(s)) - {\rm Vol}_{\wtg}\,(D_{\wtg}(\wtx, s)) \right | \leq \ve \cdot s^n
\end{equation}
where ${\rm Vol}$ denotes the Euclidean volume and $ {\rm Vol}_{\wtg}$ is the volume form for the metric $\wtg$. 
More generally, given an open set $U \subset D(s)$, for $s \leq \lambda$, we require that
\begin{equation}\label{eq-dgs3a}
\left | {\rm Vol}(U) - {\rm Vol}_{\wtg}\,( \wtexp_{\wtx,\whu}(U)) ~ \right | \leq \ve \cdot s^n
\end{equation}
\end{enumerate}
\end{defn}
Conditions (\ref{def-approxeuclid}.1) and (\ref{def-approxeuclid}.4) concern the continuity of the metric tensor $\wtg$, while
conditions (\ref{def-approxeuclid}.2-3) concern the behavior of geodesics for the metric $\wtg$, so also require control on the first and second order derivatives of $\wtg$.
The condition (\ref{def-approxeuclid}.3) is simply that the geodesics for the metric $\wtg$ and the flat metric defined by $\whu$ ``stay close''.
The conditions (\ref{def-approxeuclid}.1-5) are closely related, but are formulated separately in the form they will be used later.

\begin{lemma}\label{lem-approxeuclid2}
Assume that $\psi^g_{x, \whu} \colon D(\lambda) \to D_{\F}(x, \lambda)$ is $\ve$-approximately Euclidean. Then \begin{equation}\label{eq-dgs0a}
| \, \wtd(\vec{z}, \vec{y}) - \| \, \vec{z} - \vec{y} \, \|\,_{\whu} | ~ \leq ~ \ve \cdot \| \, \vec{z} - \vec{y} \, \|\,_{\whu} \quad \text{for ~ all} \quad \vec{y}, \vec{z} \in \wtD_i(\wtx , \lambda)
\end{equation}
Thus, conditions (\ref{def-approxeuclid}.1) and(\ref{def-approxeuclid}.2) yield, for all $\vec{a} \in D(\lambda)$,
\begin{equation}\label{eq-dgs1a}
\| \, \wtexp_{\wtx,\whu}(\vec{a}) - T_{\wtx} \circ F_{\whu}(\vec{a}) \, \|_{\whu} ~ \leq ~ 2\ve \lambda
\end{equation}
\end{lemma}
\proof
Condition (\ref{def-approxeuclid}.1) and the estimate \eqref{eq-normdef}  imply the bound $\|g^{\whu} - \delta\|_{\whu} \leq \ve$ on the matrix norm,
from which \eqref{eq-dgs0a} follows.
Condition \eqref{eq-dgs1a} then follows, as $\| \, \vec{a} \, \| \leq \lambda$.
\endproof

 Here is a key technical result about coordinate charts in Gauss normal form.
 \begin{prop}\label{prop-approxeuclid}
Given  $\ve > 0$,  there exists $\lambda_{\ve} > 0$ such that for all $x \in \fM$ with $D_{\F}(x,\lambda_{\ve}) \subset U_i$ and orthonormal frame $\whu$ of $T_x\F$, the chart $\psi^g_{x, \whu} \colon D(\lambda_{\epsilon}) \to L_x$ is $\ve$-approximately Euclidean.
\end{prop}
\proof
The claim is that $\wtexp_{\wtx,\whu}$ is well-approximated by the affine map
$T_{\wtx} \circ F_{\whu}$ for $\lambda_{\ve}$ sufficiently small. This follows from standard facts about the geodesic charts for smooth metrics, where we use the continuity of the Riemannian metric and its derivatives as functions of $x \in \fM$ to obtain uniform estimates, for all $x \in \fM$. We give a brief sketch the proof.

Let $\wth = \wth_{jk}(\xi)$ denote the Riemannian tensor on $D(\lambda)$ induced from $\wtg$ by the geodesic map $\wtexp_{\wtx,\whu}$.
Note that geodesic coordinates have the property that
$\wth_{jk}(\vec{0}) = \delta_{jk}$, the Dirac $\delta$-function.
Moreover, the Riemannian Christoffel symbols $\wtGamma_{jk}^{\ell}(\xi)$ of the metric $\wth$ also vanish at the origin.

The tensor $\wtGamma_{j k}^{\ell}(\xi)$ is $C^{\ell -1}$-continuous as a function of the metric tensor in the $C^{\ell}$ topology, for $\ell \geq 1$, so the first derivatives of $\wtGamma_{i j}^k(\xi)$ vary continuously with the metric in the $C^2$ topology, hence its curvature tensor $\wtR(\xi)$ varies continuously in the $C^2$-topology as well.
Thus, by choosing $\lambda > 0$ sufficiently small, we can assume the quantities $\| \wth_{jk}(\xi) - \delta_{jk} \|$ and $|\wtGamma_{jk}^{\ell}(\xi) |$ are arbitrarily small on the disk $ D(\lambda)$, and moreover the norm of the curvature tensor $| \wtR(\xi) |$ is uniformly bounded.

Standard results of Riemannian geometry show that the second derivatives of the geodesic map $\wtexp_{x,\whu}$ at the origin are bounded by the norms of the Christoffel symbols $\wtGamma_{jk}^{\ell}(\xi)$, of their derivatives,  and of the curvature terms $\wtR(\xi)$. (For example, see \cite[Chapter 5, Remark 2.11]{doCarmo1992}.)
Thus, given $\ve' > 0$, there exists $\lambda_{x, \ve'} > 0$ such that $\wtexp_{\wtx,\whu}$ is $\ve'$-close to its linear approximation $T_{\wtx} \circ F_{\whu}$ in the Euclidean norm on $\mR^n$. This yields the estimate \eqref{eq-dgs1}  of Definition~\ref{def-approxeuclid}.2.

The condition \eqref{eq-dgs2} of Definition~\ref{def-approxeuclid}.3   follows, as   the local expressions of the Christoffel symbols $\wtGamma_{jk}^{\ell}$ are sufficiently small on $D(\lambda)$ and the quantities $| \wtGamma_{j k}^{\ell} |$ are uniformly bounded. Conditions (\ref{def-approxeuclid}.1) and (\ref{def-approxeuclid}.4-5) follow from the continuity of the metric tensor $\wtg$, as noted above.

For each $\ve > 0$, choose $\lambda_{\ve} > 0$ so that the conditions of Definition~\ref{def-approxeuclid} holds for all $x \in \fM$,  
and any choice of orthonormal frame $\whu$ for $T_x\F$.

There is one further subtlety, which  is that the error estimates in formulas \eqref{eq-dgs1}  and \eqref{eq-dgs2}  are in terms of the leafwise distance function $\dF$, while the error $\ve'$ above is in terms of the Euclidean norm $\| \cdot \|$ on $D(\lambda)$. Introduce the constant
$$
\|\dF\| = \max \left\{ \frac{\dF(\psi^g_{x, \whu}(\vec{b}) , \psi^g_{x, \whu}(\vec{a}))}{\|\vec{b} - \vec{a}\|} ~ , \frac{\|\vec{b} - \vec{a}\|}{\dF(\psi^g_{x, \whu}(\vec{b}) , \psi^g_{x, \whu}(\vec{a}))} \mid x \in \fM ~ , ~ \whu ~ , ~ \vec{a} \ne \vec{b} \in D(\lF) \right\}
$$
Given $\ve$, let $\ve' = \ve/\|\dF\|$ and choose $\lambda_{\epsilon}$ for  as above for the error $\ve'$. Thus,  by the compactness of $\fM$ and the continuity of the metric in the $C^2$-norm,
given $\ve > 0$ there exists an $\lambda_{\ve} > 0$, so that for all $x \in \fM$, for all  $0 < \lambda \leq \lambda_{\epsilon}$,  and for all coordinate chart indices $1\leq i \leq \nu$ with $D_{\F}(x, \lambda) \subset \cP_i(x)$, the estimates \eqref{eq-dgs1} to \eqref{eq-dgs3}  of Definition~\ref{def-approxeuclid} are satisfied.
\endproof

\begin{remark}
{\rm If the leaves of $\F$ are isometric to Euclidean space $\mR^n$, such as when $\F$ is defined by a free action of $\mR^n$, then $\lambda_{\ve}$ may be chosen arbitrarily large. Otherwise, if the leaves of $\F$ have large sectional curvatures and $\ve$ is small, then $\lambda_{\ve}$ may be quite small. The points in the leafwise nets constructed in section~\ref{sec-existence} will be small compared to $\lambda_{\ve}$ so if the curvature of the leaves is ``very large'', then the net spacing will be ``very small''.}
\end{remark}

\section{Setting the constants}\label{sec-constants}

Our ultimate ``affine approximation'' results are given by Propositions~\ref{prop-framed} and \ref{prop-robustnets} in the next section, which extend Proposition~\ref{prop-approxeuclid} above. However, in order to state and prove these results, it is necessary to specify the ``universal scale constant'' $\ve_0 > 0$ for which these results are valid. It is absolutely fundamental that the estimates provided by these propositions are independent of the choices made in their applications. That is, we must \emph{a priori} define the constant $\ve_0$ as well as error bounds $\ve_1$, $\ve_2$, $\ve_3$ and $\ve_4$ which are required. For this reason, in this section we prescribe these geometric constraints, and then make our construction using these fixed choices. This section can be skipped at first reading if desired, and then consulted later, though the process of making these choices, and some of their implications as pointed out below, are a key part of the construction. 

\subsection{Number of circumscribed spheres}
The first constant to define is a very large number, based on the combinatorics of nets in regions of $\mR^n$, which is the reason for role of the dimension number $n$ in the following. It may be possible to give a much more refined value to the constant, but for our purposes, the following suffices. Set:
\begin{align}\label{eq-Cn}
C_n = \frac{10^n!}{1! \, (10^n -1)!} + \frac{10^n!}{2! \, (10^n -2)!} + \cdots + \frac{10^n!}{n! \, (10^n -n)!} + \frac{10^n!}{(n+1)! \, (10^n -n-1)!}
\end{align}

Given a finite subset $\Omega \subset D_{\F}(\xi, \rF)$ with \emph{cardinality bounded above by $10^n$},
then $C_{n}$ is an upper bound for the number of distinct subsets of $\Omega$ consisting of at most $(n+1)$-distinct points.
In particular, $C_n$ is an upper bound on the number of distinct $n$-simplices, defined by $(n+1)$-vertices in $\Omega$. Thus, $C_n$ is an upper bound on the number of circumscribed spheres for the set $\Omega$.

\subsection{Geometric constants}
Next,  introduce four additional ``geometric constants''. The purpose of these choices is briefly indicated, and their precise roles will be defined later.
The constants are ``scale-invariant'', and in their applications are multiplied by the scale $\rF$    defined in \eqref{eq-convex} below.

The \emph{width of the annular regions} appearing in Lemma~\ref{lem-volest3} will be chosen   bounded above by
\begin{equation} \label{eq-epsilon1}
\ve_1 ~ = ~ 1/(C_n \cdot 1000 n \cdot 100^n) . 
\end{equation}
The \emph{thickness of the rectangular regions} appearing in the robustness condition \eqref{eq-slabvol} will be chosen   bounded above by
\begin{equation}\label{eq-epsilon2}
\ve_2 = 1/(C_n \cdot 2000 \cdot 2^n) .
\end{equation}
The  bound on the \emph{translation distance of the centers of circumscribed spheres} for a perturbed net is 
\begin{equation}\label{eq-epsilon3}
\ve_3 ~ = \ve_1/10 .
\end{equation}
The constant  $\ve_3$  first appears in the statement and proof of Proposition~\ref{prop-inductsphere}. We   repeatedly use the implication $\ve_3 < \ve_1/4$. 

The  \emph{error of the affine approximation}   in Proposition~\ref {prop-framed} is bounded by a   constant $\ve_4$, which 
 determines the recursive decrease in the robustness estimates in
Propositions~\ref{prop-robustvarest}, \ref{prop-robustnets} and \ref{prop-inductsphere}. The value of $\ve_4$ is defined by a recursive process, depending on the dimension $n$, which we recall from the proof of Proposition~\ref{prop-robustvarest}.

Proposition~\ref{prop-robustvarest} gives a recursive definition for the functions $\rho_m(\rho, \ve, e_1, e_2)$ for $1 \leq m \leq n$. As noted there, the function $\rho_m(\rho, \ve, e_1, e_2)$ is monotone increasing in $e_2$ and $\rho$, and monotone decreasing in $e_1$ and $\ve$, and satisfies
$\rho_m(s \cdot \rho, s \cdot \ve, s \cdot e_1, s \cdot e_2) = s \cdot \rho_m(\rho, \ve, e_1, e_2)$ for $s > 0$.
Moreover, for all $1 \leq m \leq n$,
$\rho_m(\rho, 0, e_1, e_2) = \rho$.
For the normalized values $e_1 = 1$, $e_2 = 2$, $e_4 = 4(e_2 + e_1) = 12$, and $\rho = \rho_0$, define functions $\rho_m(\rho_0, \ve)$ recursively by
$$ \rho_0(\rho_0, \ve) = \rho_0 ~ , ~ \delta_1 = 2 ~ , ~ \rho_1(\rho_0, \ve) = \rho_0 - 2 \ve ~ , ~ \delta_2 = 50 ~ , ~ \rho_2(\rho_0 , \ve) = \rho - 50 \ve$$
and for $1 < m \leq n$, by
\begin{equation} \label{eq-robust-norm}
\rho_m(\rho_0, \ve) ~ = ~ \rho_0 - \ve \cdot \delta_m \quad , \quad \delta_m ~ = ~ 50 + 24 \cdot \sum_{k=2}^{m-1} ~ \frac{(1 + \delta_k)}{\rho_k(\ve)}
\end{equation}
Note that each $\rho_m(\ve)$ is a continuous function of $\ve$. Also, for fixed initial data $(\rho_0, \ve)$, the sequence of values is monotone decreasing in $m$:
$$\rho_0 = \rho_0(\rho_0, \ve) > \rho_1(\rho_0, \ve) > \rho_2(\rho_0, \ve) > \cdots > \rho_n(\rho_0, \ve) > 0$$

At a key stage of the induction process, we introduce the following constants, 
for each $0 \leq k \leq n$:
\begin{eqnarray*}
\whrho_k ~ & = & ~ (18 - 2k/3n ) \cdot \ve_2\\
\whrhop_k ~ & = & ~ (18 - (2k+1)/3n ) \cdot \ve_2
\end{eqnarray*}
 Then we have
\begin{equation}\label{eq-epsilon4h}
18 \ve_2 = \whrho_0 > \whrhop_0 > \whrho_1 > \whrhop_1 > \cdots > \whrho_n > \whrhop_n > \whrho_{n+1} > \whrhop_{n+1} > 15 \ve_2
\end{equation}

Finally,   choose $\ve_4 > 0$ sufficiently small so that the following $2n+2$ inequalities hold, for $1 \leq k \leq n+1$:
\begin{equation}\label{eq-epsilon4a}
\whrho_k > \rho_n(\whrho_k, 10\ve_4) > \whrhop_k + \ve_2/100  
\end{equation}
\begin{equation}\label{eq-epsilon4b}
  \whrhop_k > \rho_n(\whrhop_k, 10\ve_4) > \whrho_{k+1} + \ve_2/100 
\end{equation}
The full set of these inequalities are used in the proofs of Propositions~\ref{prop-inductrobust} and ~\ref{prop-inductsphere}, where they are multiplied by the scale $s = \rF/10$.

\subsection{Error of transverse computations}
Finally, $\ve_0$ is the ``basic error'' appearing in almost every transverse translation calculation and estimate, so is restricted by multiple conditions. The following restrictions  are informally summarized by saying ``it is intuitively clear that there exists $\ve_0$ \emph{sufficiently small} so that all of these conditions are satisfied''. We make this intuition precise:
\begin{equation}\label{eq-epsilon0}
\text{Choose}~ \ve_0 > 0 ~ \text{which satisfies the following   conditions:}   \hspace{160pt}
\end{equation}
\begin{enumerate}
\item $\ve_0 < 1/2000$ -- used in equations \eqref{eq-netuniform1} and \eqref{eq-netuniform2}
\item $ \ve_0 \leq 50 \, n \, (2/5)^n \ve_1 $ -- used in equation \eqref{eq-volestep0}
\item $\ve_0 < \ve_2/2000$ -- used in equations \eqref{eq-slabs} and \eqref{eq-slabsest} and in proof of Proposition~\ref{prop-inductrobust}
\item $\ve_0 < \ve_3/4$ -- used in equations \eqref{eq-leafvar8}, \eqref{eq-leafvar12}), \eqref{eq-radest6}) and in proof of Proposition~\ref{prop-inductsphere}
\item $ \ve_0 < \ve_1/2 < \ve_1 - 2\ve_3$ -- used in \eqref{eq-leafvar20})
\item $ \ve_0 < \ve_3 /2 \{1 + 35 n^{3/2} \cdot (4/15\ve_2 )^{n-1} \} $ -- used in equation \eqref{eq-leafvar6cat} 
\item $ \ve_0 < \ve_4/20$ -- used in Proposition~\ref{prop-framed}
\item $ \ve_0 < \delta_n(\ve_4 )/100$ for $ \delta_n$ defined in Lemma~\ref{lem-GS}
\end{enumerate}
Note that the function $\delta_n$ in estimate (8) above, as defined in Lemma~\ref{lem-GS}, is   independent of all other choices, so that $\ve_0$ is well-defined.

\subsection{Leafwise constants}

Recall that $\lambda_{\ve_0}$ was defined in the proof of Proposition~\ref{prop-approxeuclid}. Introduce the fundamental ``leafwise'' constant:
\begin{equation}\label{eq-convex}
\rF = \min \{\dFU , \lF/5 , \lambda_{\ve_0}, 1\}
\end{equation}
which is the basic distance scale for all of our subsequent constructions, chosen so that the leafwise balls $D_{\F}(\xi, \rF)$ are ``$\ve_0$-approximately Euclidean''.
For example, if the leaves of $\F$ are isometric to Euclidean space $\mR^n$, then $\rF = \min \{\dFU , \lF/5, 1\}$. Otherwise, if the leaves of $\F$ have large sectional curvatures, then $\rF$ may be quite small.

\subsection{Variations of the metric in charts}\label{subsec-r*}

Recall the definitions of the functions $var$ in \eqref{eq-var}) and $div$ in \eqref{eq-divpairwise2}, and choose the ``transverse'' scale constant
$r_* > 0$ so that $div(z, r_*) \leq \ve_0 \rF$ for all $z \in \fM$,
and also $var(i,r_*) \leq \ve_0 \rF$ for all $1 \leq i \leq \nu$.

\section{Affine distortion estimates}\label{sec-distortions}

The stability of the Delaunay triangulation associated to a net in $\mR^n$ follows from   delicate linear algebra estimates in sections~\ref{sec-euclidean} and \ref{sec-inequalities}. It is thus natural that the analogs of these estimates for non-Euclidean Riemannian geometry, as given in Propositions~\ref{prop-framed} and \ref{prop-robustnets}, are even more delicate. In this section, we show how these Propositions follow from the   choices in the previous section~\ref{sec-constants}.

\begin{defn}\label{def-epsilon-isometry}
Let $\{X, d_X\}$ and $\{Y, d_Y\}$ be metric spaces, and $\e > 0$.
A homeomorphism into $\phi \colon X \to Y$ is said to be an \emph{$\e$-isometry} if
\begin{equation}\label{eq-epsilon-isometry}
d_X(x,x') - \e ~ \leq ~ d_{Y}(\phi(x), \phi(x')) ~ \leq ~ d_X(x,x') + \e \quad \text{for all} ~ x,x' \in X
\end{equation}
\end{defn}

\medskip

\begin{lemma}\label{lem-geodesicisom}
For $x \in \fM$ and orthonormal frame $\whu$ for $T_x\F$, the geodesic normal coordinate map
$\ds \psi^g_{x, \whu} \colon D(\rF/2) \to D_{\F}(x, \rF/2)$
is an $(\ve_0  \rF)$-isometry from the metric $\| \cdot \|$ to the metric $\dF$.
\end{lemma}
\proof
As the Euclidean disk $D(\rF/2)$ has diameter $ \rF$,   the claim then follows from the estimate $\|\wtg^{\whu} - \delta\|_{\whu} \leq \ve_0$ as in the proof of Lemma~\ref{lem-approxeuclid2}.
\endproof

Recall that the disk section $\fZ(y, i, r_*)$ of radius $r_* > 0$ was defined by formula \eqref{eq-transsec2}).
\begin{lemma}\label{lem-secdiv}
Let $x \in U_i$ and $y \in \cP_i(x) \cap U_j$ for some $1 \leq i, j \leq \nu$.
Assume that $x' \in \fZ(x, i, r_*)$ and $y' = \fZ(y, j, r_*) \cap \cP_i(x')$.
Then
\begin{equation}\label{eq-error1}
\dF(x,y) - 2\ve_0 \, \rF~ \leq ~ \dF(x',y') ~ \leq ~ \dF(x,y) + 2\ve_0 \, \rF
\end{equation}
If either $i=j$ or $x=y$, then a more strict estimate holds:
\begin{equation}\label{eq-error2}
\dF(x,y) - \ve_0 \, \rF~ \leq ~ \dF(x',y') ~ \leq ~ \dF(x,y) + \ve_0 \, \rF
\end{equation}
\end{lemma}
\proof
Let $y'' = \fZ(y, i, r_*) \cap \cP_i(x')$. Then $\dF(y',y'') \leq \ve_0 \, \rF$
by the definition of the divergence \eqref{eq-divpairwise2} and of $r_*$.
Thus, $|\dF(x',y'') - \dF(x',y')| \leq \ve_0 \, \rF$ by the triangle inequality.

Then by the definition of the variation \eqref{eq-var} and $r_*$ we also have
$| \dF(x',y'') - \dF(x,y) | \leq \ve_0 \, \rF$. The estimates \eqref{eq-error1} and \eqref{eq-error2} then follow.
\endproof

\subsection{Estimates for variations of affine geometries}

We next derive estimates comparing  the local \emph{affine  geometry} of geodesic coordinates in nearby plaques at nearby points.

For $x \in \fM$ with  $D_{\F}(x, \rF) \subset U_i$ let $y \in D_{\F}(x, \rF/2)$ so that $D_{\F}(y, \rF/2) \subset D_{\F}(x, \rF)$.

Let $x' \in \fZ(x, i, r_*)$ and set $y' = \phi_i(x,x')(y)$.
Choose orthonormal frames $\whu$ for $T_{x}\F$ and $\whvp$ for $T_{y'}\F$,
with corresponding geodesic coordinates $\psi^g_{x,\whu}$ and $\psi^g_{y',\whvp}$. Consider the composition
\begin{equation}\label{eq-Psixy0}
\Psi_{x,y'}' \equiv (\psi^g_{y',\whvp})^{-1} \circ \phi_i(x,x') \circ \psi^g_{x,\whu} \circ T_{\xi} \colon D(\rF) \to \mR^n
\end{equation}
where $\xi = (\psi^g_{x, \whu})^{-1}(y)$, and $T_{\xi} \colon \mR^n \to \mR^n$ denotes the affine transformation   $T_{\xi}(\vec{a}) = \vec{a} + \xi$. Then
$$\Psi_{x,y'}'(\vec{0}) = (\psi^g_{y',\whvp})^{-1} \circ \phi_i(x,x') \circ \psi^g_{x,\whu}(\xi) = (\psi^g_{y',\whvp})^{-1} \circ \phi_i(x,x')(y) = (\psi^g_{y',\whvp})^{-1}(y') = \vec{0}$$
The map $\Psi_{x,y'}'$ compares two coordinate systems about the point $y'$: one is the translate of the geodesic coordinates $\psi^g_{x,\whu}$ centered at $x$ but restricted to a neighborhood of $y$ in its domain, and the other is centered at the translated point $y'$.
Each coordinate system defines an ``affine structure'' in a neighborhood of $y'$. The next result shows that $\Psi_{x,y'} $ can be made ``almost the identity'' by the proper choice of the framing $\whvp$, so that these affine structures are arbitrarily close. The proof of this result is surprisingly complex. 

\begin{prop}\label{prop-framed}
There exists a choice of orthonormal frame $\whv$ for $T_{y'}\F$ so that
\begin{equation}\label{eq-Psixy}
\Psi_{x,y'} \equiv (\psi^g_{y',\whv})^{-1} \circ \phi_i(x,x') \circ \psi^g_{x,\whu} \circ T_{\xi} \colon D(\rF/2) \to \mR^n
\end{equation}
is ~ $\ve_4   \rF$-close to the identity, for $\ve_4$ chosen to satisfy    \eqref{eq-epsilon4a} and \eqref{eq-epsilon4b}.
\end{prop}
\proof
We are given    $x \in D_{\F}(x, \rF) \subset U_i$,  $y \in D_{\F}(x, \rF/2)$ so that $D_{\F}(y, \rF/2) \subset D_{\F}(x, \rF)$, and 
 $x' \in \fZ(x, i, r_*)$ and set $y' = \phi_i(x,x')(y)$. Also given are   frames $\whu$ for $T_{x}\F$ and $\whvp$ for $T_{y'}\F$. 

The idea of the proof is simple, in that we express both geodesic coordinate maps $\psi^g_{x,\whu}$ 
and $\psi^g_{y',\whvp}$ in the local coordinates $\vp_i$, 
as in the proof of Proposition~\ref{prop-approxeuclid}. The key point of the proof follows from some delicate linear algebra, used to choose the new framing $\whv$,  and then    estimating the   distortion  of the geodesic coordinates for this frame.

Let $\vp_i(x) = (\wtx, w_x)$,  $\vp_i(y) = (\wty, w_x)$ and $\vp_i(y') = (\wtyp, w_{y'})$ for $w_x , w_{y'} \in \fT_i$, where as before in section~\ref{subsubsec-compargeod}, $\wtx = \lambda_i(x)$ and $\wty = \lambda_i(y)$. By definition,  $\phi_i(x,x')$ is the identity map when expressed in the coordinate system $\vp_i$, so the assumption $y' = \phi_i(x,x')(y)$ implies $\wtyp = \wty$. 

We mention a point of notation established in section~\ref{sec-approx} and used repeatedly below. The ``tilde'' notation, $\wtx \in (-1,1)^n$ for example,   denotes the horizontal coordinates of a point or set;  the ``prime'' notation denotes a point  or set in the translated plaque $\cP_i(y')$; while $\vec{v} \in \mR^n$   denotes a vector in the vector space $\mR^n$, typically obtained from the inverse of the geodesic coordinates $\psi^g$.

Set $d_2 = \rF/5$. 
The restriction $ \phi_i(x,x') \colon D_{\F}(x,\rF/2) \to \cP_i(w_{y'})$ is an $\ve_0\rF$-isometry by
Lemma~\ref{lem-secdiv}, and thus 
$\ds \phi_i(x,x') (D_{\F}(y,2d_2)) \subset D_{\F}(y',\rF/2))$.
Indeed,
   $$\phi_i(x,x')(D_{\F}(y,2d_2)) \subset D_{\F}(y',2\rF/5 + \ve_0\rF),$$
and since by assumption $\ve_0 < 1/2000$ we have
   $$2 \rF/5 + 2 \ve_0 \rF \leq 401/1000 \rF < \rF/2.$$
This implies that the composition \eqref{eq-Psixy} is well-defined. 
Recall from section~\ref{subsubsec-compargeod} that we denoted
\begin{eqnarray*}
\wtD_i(\wtx , \rF/2) & = & \vp_i(D_{\F}(x,\rF/2)) \subset (-1,1)^n \times \{w_x\},\\
\wtD_i'(\wtyp , \rF/2) & = & \vp_i(D_{\F}(y',\rF/2)) \subset (-1,1)^n \times \{w_{y'}\}.
\end{eqnarray*}

Recall   from section~\ref{subsubsec-compargeod} that $\wtd$ denotes the distance function on $\wtD_i(\wtx , \rF/2)$ defined by the leafwise metric $\dF$ via the coordinates $\vp_i$, and $\wtg$ denotes the induced Riemannian metric on $\wtD_i(\wtx , \rF/2)$.
The geodesic coordinates about a neighborhood of $\wtx \in \mR^n$ associated to $\wtg$ and $\whu$,  are denoted by $\wtexp_{\wtx,\whu} \colon D(\rF/2) \to \wtD_i(\wtx, \rF/2)$. So for  $\xi = (\psi^g_{x, \whu})^{-1}(y)$, then $\wty = \wtexp_{\wtx ,\whu}(\xi)$ by definition.  

Similarly, $\wtdp$ denotes the distance function induced on $\wtD_i'(\wtyp , \rF/2)$, and $\wtgp$ denotes the induced metric tensor on $\wtD_i'(\wtyp, \rF/2)$.
The geodesic coordinates associated to $\wtgp$ and $\whvp$, centered at $\wtyp$, are denoted by $\wtexp_{\wtyp,\whvp} \colon D(\rF/2) \to \wtD_i'(\wtyp , \rF/2)$.
Then the map $\Psi_{x,y'}$ from \eqref{eq-Psixy} can be expressed by 
$$\Psi_{x,y'} = \wtexp_{\wtyp,\whvp}^{-1} \circ \wtexp_{\wtx ,\whu} \circ T_{\xi}$$ 
and there is the diagram of maps:
\begin{equation}\label{eq-diagram}
\begin{array}{ccccc}
T_x\F \cong \mR^n \supset D(\rF/2) & \stackrel{\wtexp_{\wtx,\whu} \circ T_{\xi}}{\longrightarrow} & \wtD_i(\wtx , \rF/2) \subset (-1,1)^n \times \{w_x\} & \stackrel{\vp_i^{-1}}{\longrightarrow} & \cP_i(w_x) \\
& & & & \\
\Psi_{x,y'} \downarrow & & = ~ \downarrow & & \downarrow \phi_i(x,x') \\
& & & & \\
T_{y'}\F \cong \mR^n \supset D(\rF/2) & \stackrel{\wtexp_{\wtyp,\whvp}}{\longrightarrow} & \wtD_i'(\wtyp , \rF/2) \subset (-1,1)^n \times \{w_{y'}\} & \stackrel{\vp_i^{-1}}{\longrightarrow} & \cP_i(w_{y'})
\end{array}
\end{equation}

Set $\vec{\gamma} = \wtx + \whu \cdot \xi$ so that $T_{\vec{\gamma}} \circ F_{\whu} (\vec{a}) = T_{\wtx} \circ F_{\whu} (\vec{a} + \xi)$.

By condition \eqref{eq-dgs1} of Definition~\ref{def-approxeuclid}, and using that $\rF \leq \lambda_{\epsilon_0}$,
for all $\vec{a} \in D(\rF/2)$ we have
\begin{equation}\label{eq-linapprox}
\wtd \left (\wtexp_{\wtx ,\whu}(\vec{a} + \xi) , T_{\vec{\gamma}} \circ F_{\whu} (\vec{a}) \right) \leq \ve_0 \rF ~ , ~
\wtdp \left (\wtexp_{\wtyp ,\whvp}(\vec{a}) , T_{\wtyp} \circ F_{\whvp} (\vec{a}) \right) \leq \ve_0 \rF
\end{equation}
Set $\vec{a} = \vec{0}$ in the first estimate of \eqref{eq-linapprox}, then by Lemma~\ref{lem-approxeuclid2} we obtain
\begin{equation}\label{eq-offsetest}
\| \wty - \vec{\gamma} \|_{\whu} ~ \leq ~  \wtd (\wty , \vec{\gamma} ) + \ve_0 \rF ~ = ~ \wtd (\wtexp_{\wtx ,\whu}(\xi) , T_{\vec{\gamma}}(\vec{0}) ) + \ve_0 \rF ~ \leq ~ 2\ve_0 \rF
\end{equation}

The ``obvious'' next step is to replace the orthonormal framing $\whvp$ for $\mR^n$ for the norm $\| \cdot \|_{\whvp}$ with the new framing $\whv = \whu$, and then the claim of Proposition~\ref{prop-framed} would   follow. However, $\whu$ need not be an orthonormal framing the norm $\| \cdot \|_{\whvp}$, so it is necessary to adjust the framing $\whu$ using the Gram-Schmidt orthogonalization process. This introduces additional errors, which depend on the ``distance'' from $\whu$ to $\whvp$ in the Lie group $GL(\mR^n)$. We formulate this error as follows, using estimates derived from   the Gram-Schmidt orthogonalization process. This   derivation of the following result is straightforward, and we omit the proof.

\begin{lemma}\label{lem-GS}
Let $\mR^n$ have the standard Euclidean inner product with norm $\| \cdot \|$.
There exists $\e_n > 0$ and a monotone continuous function $\delta_n \colon [0, \e_n] \to [0, \e_n]$ with $\delta_n(0) = 0$, such that given $0 < \e \leq \e_n$ set $\delta = \delta_n(\e) > 0$, then for any basis
$\{ \vec{f}_1', \ldots, \vec{f}_n'\} \subset \mR^n$, whose vectors satisfy
\begin{enumerate}
\item $1 -\delta < \|\vec{f}_j'\| < 1+\delta$, for $1 \leq j \leq n$,
\item $| \vec{f}_i' \bullet \vec{f}_j' | < \delta$, for $1 \leq i \ne j \leq n$,
\end{enumerate}
then there exists orthonormal vectors $\{\vec{f}_1, \ldots, \vec{f}_n\}$ such that
$\|\vec{f}_j - \vec{f}_j' \| \leq \e$. \hfill $\Box$
\end{lemma}
Recall that $\whe = \{\vec{e}_1 , \ldots , \vec{e}_n\}$ is the standard orthogonal basis for $\mR^n$.
Scale these unit vectors by a factor of $d_2 = \rF/5$ so they lie in the domain of $\Psi_{x,y'}$, and set $\vec{z}_j = F_{\whu}(d_2 \vec{e}_j)$. Note that $\| \vec{z}_j \|_{\whu} = d_2$.

Recall that $\xi = (\psi^g_{x, \whu})^{-1}(y)$. Then   $d_2 \vec{e}_j + \xi \in D(\wtx , \rF/2)$ so that 
$\ds \wtz_j = \wtexp_{\vec{x},\whu}( d_2 \vec{e}_j + \xi) \in \wtD_i(\vec{x} , \rF/2)$ is well-defined.
Then by Lemma~\ref{lem-geodesicisom},
\begin{eqnarray}
| \, \wtd(\wtz_j , \wty) - d_2 \, | = | \, \wtd(\wtz_j , \wty) - \| \vec{z}_j \|_{\whu} \, | = | \, \wtd(\wtz_j , \wty) - \| \, \whu \cdot (d_2 \vec{e}_j + \xi) - \whu \cdot \xi \, \|_{\whu} \, | & \leq & \ve_0 \rF \label{eq-GSest1} \\
|\, \wtd(\wtz_j , \wtz_k) - \sqrt{2} d_2 \, | = |\, \wtd(\wtz_j , \wtz_k) - \| \whu \cdot (d_2 \vec{e}_j + \xi) - \whu \cdot (d_2 \vec{e}_k + \xi) \|_{\whu} \, |
& \leq & \ve_0 \rF \label{eq-GSest2}
\end{eqnarray}
The estimates \eqref{eq-GSest1} and \eqref{eq-GSest2} imply that the set $\{(\wtz_1 - \wty)/d_2 , \ldots , (\wtz_n - \wty)/d_2\}$ is an ``almost orthonormal'' collection for the metric $\wtd$.

By Lemma~\ref{lem-secdiv}, the map $\phi_i(x,x')$ is an $\ve_0 \rF$-isometry,
and as $\phi_i(x,x')$ is the identity map in the coordinates $\vp_i$,
we obtain   estimates corresponding to \eqref{eq-GSest1} and \eqref{eq-GSest2} for the metric $\wtdp$,
\begin{eqnarray}
| \, \wtdp(\wtz_j , \wty) - d_2 \, | ~ \leq ~ 2 \ve_0 \rF \quad & \Longrightarrow & \quad
| \, \| \, \wtz_j - \wty\, \|_{\whvp} - d_2 \, | ~ \leq ~ 3 \ve_0 \rF \label{eq-GSest3}\\
|\, \wtdp(\wtz_j , \wtz_k) - \sqrt{2} d_2 \, | ~ \leq ~ 2 \ve_0 \rF \quad & \Longrightarrow & \quad
|\, \| \, \wtz_j - \wtz_k \, \|_{\whvp} - \sqrt{2} d_2 \, | ~ \leq ~ 3 \ve_0 \rF \label{eq-GSest4}
\end{eqnarray}
which follow from estimate \eqref{eq-dgs0a} of Lemma~\ref{lem-approxeuclid2}.

Define $\vzpj = \wtexp_{\wtyp,\whvp}^{-1}(\wtz_j)$.
Then by estimate \eqref{eq-dgs1a} of Lemma~\ref{lem-approxeuclid2}, and noting that $\vzpj \in \wtD_i(\wtyp , \rF/2)$,
\begin{equation}\label{eq-approx4}
\| \, \wtz_j - T_{\wtyp} \circ F_{\whvp}(\vzpj) \, \|_{\whvp} ~ = ~
\| \, \wtexp_{\wtyp,\whvp}(\vzpj) - T_{\wtyp} \circ F_{\whvp}(\vzpj) \, \|_{\whvp} ~ \leq ~
\ve_0 \rF
\end{equation}

Then by \eqref{eq-GSest3} and \eqref{eq-GSest4}, and using that the map $T_{\wtyp} \circ F_{\whvp}$ is an isometry from the norm $\| \cdot \|$ to the norm $\| \cdot \|_{\whvp}$,
we obtain for the Euclidean norm on $\mR^n$, for $1 \leq j \ne k \leq n$,
\begin{eqnarray}
| \, \| \vzpj \| - d_2 \, | ~ & \leq & ~ 3\ve_0 \rF \label{eq-GSest5} \\
|\, \| \vzpj - \vzpk \| - \sqrt{2} d_2 \, | ~ & \leq & ~ 3\ve_0 \rF \label{eq-GSest6}
\end{eqnarray}

Set $\vec{f}_j' = \vzpj/d_2$, and observe that \eqref{eq-GSest5} implies the collection $\{ \vec{f}_1', \ldots , \vec{f}_n'\}$ satisfies hypothesis (1) of Lemma~\ref{lem-GS} for $\delta = 15 \ve_0$.

It remains to estimate $| \vec{f}_j' \bullet \vec{f}_k' |$ for $1 \leq j \ne k \leq n$.
Write $\vec{f}_k' = \vec{f}_{j,k}' + \vec{f}_{k, k}'$
where $\vec{f}_{j,k}'$ is collinear with $\vec{f}_j'$ and
$ \vec{f}_j' \bullet \vec{f}_{k,k}' = 0$.
Then $| \vec{f}_j' \bullet \vec{f}_k' | = | \vec{f}_j' \bullet \vec{f}_{j,k}' | = \| \vec{f}_j' \| \| \vec{f}_{j,k}' \|$.

Note also that
$\|\vec{f}_j'\|^2 = \|\vec{f}_{k, k}' \|^2 + \| \vec{f}_{j,k}' \|^2$ hence $ \|\vec{f}_{k,k}'\|^2 = (\|\vec{f}_j'\|^2 - \| \vec{f}_{j,k}' \|^2)$.
By \eqref{eq-GSest6} we have
\begin{equation}\label{eq-GSest7}
\sqrt{2} - 15\ve_0 ~ \leq ~ \|\vec{f}_j' - \vec{f}_k' \| = \|(\vec{f}_j' - \vec{f}_{j,k}') - \vec{f}_{k,k}' \| ~ \leq ~ \sqrt{2} + 15 \ve_0
\end{equation}
After squaring and using the orthogonality of the vectors, we obtain
$$ 2 - 30\sqrt{2} \ve_0 + 225 \ve_0^2 ~ \leq ~ \|(\vec{f}_j' - \vec{f}_{j,k}')\|^2 + \|\vec{f}_{k,k}' \|^2 ~ \leq ~ 2 + 30 \sqrt{2} \ve_0 + 225 \ve_0^2 $$
 Note that  $\vec{f}_j'$ and $\vec{f}_{j,k}'$ are collinear, hence
$ \| \vec{f}_j' - \vec{f}_{j,k}' \|^2 = \| \vec{f}_j' \|^2 - 2 \| \vec{f}_j'\| \cdot \| \vec{f}_{j,k}' \| + \| \vec{f}_{j,k}' \|^2$.
Then  using that  $ 2 - 100 \ve_0 < 2 - 30\sqrt{2} \ve_0 + 225 \ve_0^2$, we have
$$ 2 - 100 \ve_0 ~ < ~
\| \vec{f}_j' \|^2 - 2 \| \vec{f}_j'\| \cdot \| \vec{f}_{j,k}' \| + \| \vec{f}_{j,k}' \|^2 + \| \vec{f}_{k,k}' \|^2
~ < ~ 2 + 100 \ve_0 $$
From the identity $ \|\vec{f}_{k,k}' \|^2 = (\|\vec{f}_j'\|^2 - \| \vec{f}_{j,k}' \|^2)$ one derives
$\|\vec{f}_j' \bullet \vec{f}_k' \| = \| \vec{f}_{j,k}' \| < 100 \ve_0$.

Thus, the the collection $\{ \vec{f}_1', \ldots , \vec{f}_n'\}$ satisfies both hypotheses of Lemma~\ref{lem-GS} for $\delta = 100 \ve_0$. We assume that $\ve_0 < \ve_4/20$ in \eqref{eq-epsilon0} so that
$\ve_4 < \ve_5 = 2\ve_4 - 20 \ve_0$. By choice of $\ve_0$ in \eqref{eq-epsilon0}, we have $100 \ve_0 < \delta_n(\ve_5)$ so we obtain the orthonormal framing $\whf = \{\vec{f}_1, \ldots, \vec{f}_n\}$ of $\mR^n$ satisfying $\|\vec{f}_k - \vec{f}_k'\| \leq \epsilon_5$

Define $\vec{v}_j = F_{\whvp}( \vec{f}_j)$, then $\whv = \{\vec{v}_1, \ldots, \vec{v}_n\}$ is an orthonormal frame
for the norm $\| \cdot \|_{\whvp}$ so defines an orthonormal framing of $T_{y'}\F$. Note that
$F_{\whv} = F_{\whvp} \circ F_{\whf}$ and calculate:
\begin{eqnarray*}
\lefteqn{ \| T_{\vec{\gamma}} \circ F_{\whu}(d_2 \vec{e}_j ) - T_{\wtyp} \circ F_{\whv}(d_2\vec{e}_j) \|_{\whvp} \quad } \\
~ & \leq & ~ \| ( T_{\vec{\gamma}} \circ F_{\whu}(d_2 \vec{e}_j ) - \wtz_j \|_{\whvp} + \| \wtz_j - T_{\wtyp} \circ F_{\whv}(d_2\vec{e}_j) \|_{\whvp} \\
~ & = & ~ \| T_{\vec{\gamma}} \circ F_{\whu}(d_2 \vec{e}_j ) - \wtz_j \|_{\whvp} + \| \wtz_j - T_{\wtyp} \circ F_{\whvp}(d_2\vec{f}_j) \|_{\whvp} \\
~ & \leq & ~ \| T_{\vec{\gamma}} \circ F_{\whu}(d_2 \vec{e}_j ) - \wtz_j \|_{\whvp} + \| \wtz_j - T_{\wtyp} \circ F_{\whvp}(d_2\vec{f}_j') \|_{\whvp} + \|T_{\wtyp} \circ F_{\whvp}(d_2\vec{f}_j') - T_{\wtyp} \circ F_{\whvp}(d_2\vec{f}_j) \|_{\whvp} \\
~ & \leq & ~ 3\ve_0\rF + \ve_0 \rF + 2d_2 \ve_5 = 4 \ve_0 \rF + 2\ve_5 \rF/5
\end{eqnarray*}
where we use successively the definitions of the quantities involved, Lemmas~\ref{lem-approxeuclid2} and \ref{lem-geodesicisom}, the estimate \eqref{eq-linapprox}, the estimate \eqref{eq-approx4}, and Lemma~\ref{lem-GS}.
Then by the approximations \eqref{eq-linapprox} and Lemmas~\ref{lem-approxeuclid2} and \ref{lem-geodesicisom}, we have for all $\vec{a} \in D(2\rF/5)$ that
\begin{equation}\label{eq-linapprox3}
\wtd \left (\wtexp_{\wtx ,\whu}(\vec{a} + \xi) , \wtexp_{\wtyp ,\whvp}(\vec{a}) \right) \leq 7 \ve_0 \rF + 2\ve_5 \rF/5
\end{equation}
Hence by Lemma~\ref{lem-geodesicisom} and our choice $\ve_5 = 2\ve_4 - 20 \ve_0$, we obtain
\begin{equation}\label{eq-linapprox4}
\|\Psi_{x,y'}(\vec{a}) - \vec{a} \| = \|\wtexp_{\wtyp,\whvp}^{-1} \circ \wtexp_{\wtx ,\whu} \circ T_{\xi}(\vec{a}) \| \leq 8 \ve_0 \rF + 2\ve_5 \rF/5 < \ve_4 \rF
\end{equation}
completing the proof of Proposition~\ref{prop-framed}.
\endproof

\subsection{Robustness criteria}

The fine control of the affine structure of geodesic coordinates provided by Proposition~\ref{prop-framed}
is used in establishing robustness criteria for leafwise Delaunay triangulations in the next section. In preparation, we define a ``non-linear'' form of the robustness criteria in Definition~\ref {def-robust1} and Proposition~\ref{prop-robustvarest}.

Recall that ${\rm Span}(\vec{v}_0, \ldots , \vec{v}_k) \subset \mR^n$ is the \emph{affine} span
of the vectors $\{\vec{v}_0, \ldots , \vec{v}_k\}$.

\begin{defn}\label{def-robustnets}
Let $\rho > 0$ and $x \in U_i$ such that $D_{\F}(x,\rF) \subset \cP_i(x)$. Let $1 \leq m \leq n$.
A set $\{y_0, \ldots , y_m\} \subset D_{\F}(x,\rF/2)$ is \emph{$\rho$-robust} if for each $1 \leq k < m$, the following leafwise metric conditions hold:
\begin{enumerate}
\item Fix an orthonormal frame $\whu = \{\vec{u}_1, \ldots, \vec{u}_n\} \subset T_{y_k}\F$ with geodesic coordinates $\psi^g_{y_k, \whu}$;
\item for each $0 \leq j \leq k$, set $\vec{v}_j = (\psi^g_{y_k, \whu})^{-1}(y_j)$;
\item Set $\ds H(y_0 , \ldots, y_{k}; y_k) = \psi^g_{y_k, \whu} \left\{ {\rm Span}(\vec{v}_0, \ldots , \vec{v}_{k}) \cap D(\rF) \right\}$.
\end{enumerate}
Then the point $y_{k+1}$ lies at distance at least $\rho$ from the submanifold $\ds H(y_0 , \ldots, y_{k}; y_k)$
\end{defn}

We   show that the robustness condition for points in Definition~\ref{def-robustnets} implies the robustness condition Definition~\ref{def-robust1} holds for their vector coordinates in geodesic normal coordinates.
\begin{prop}\label{prop-robustnets}
Given constants
\begin{equation}
\rF/10 = d_1 < d_1 + 2\ve_0 \rF < e_1 < e_2 < d_2 - 2\ve_0 \rF < d_2 = 2\rF/10
\end{equation}
and $0 < \rho_0 < d_1$, let $x \in \fM$ and suppose $\{y_0, \ldots , y_m\} \subset D_{\F}(x,\rF/2)$ satisfy
  $\ds e_1 \leq \dF(y_j, y_k) \leq 2e_2$ for $0 \leq j \ne k \leq m$ and 
  $\{y_0, \ldots , y_m\} $ is $\rho_0$-robust.

Given an orthonormal frame $\whu$ of $T_{x}\F$, set
$\vec{w}_j = (\psi^g_{x, \whu})^{-1}(y_j) \in D(\rF/2)$ for $0 \leq j \leq m$.
Then $\{ \vec{w}_0, \ldots , \vec{w}_m\} \subset \mR^n$
is $\rho_m$-robust, where
$\rho_{\ell} = \rho_{\ell}(\rho_0, \ve_4 \rF, d_1, d_2)$ is defined by \eqref{eq-robust-norm} for $1 \leq \ell \leq m$.
\end{prop}
\proof
We proceed by induction.  By assumption,
$\dF(y_0 , y_1) \geq e_1 > d_1+ \ve_0 \rF$,
so by Lemma~\ref{lem-geodesicisom} we have $\| \vec{w}_1 - \vec{w}_0 \| \geq d_1 \geq \rho_0 - 2\ve_4 \rF = \rho_1$.

Now assume that the collection $\{\vec{w}_0 , \ldots , \vec{w}_{\ell}\}$ is $\rho_{\ell}$-robust, for each $1 \leq \ell < m$.
We show that $\{ \vec{w}_0, \ldots , \vec{w}_{\ell+1}\}$ is $\rho_{{\ell}+1}$-robust.

Let $U_i$ be a coordinate chart such that $D_{\F}(x, \rF) \subset U_i$.

Let $T_{\vec{w}_{\ell}} \colon \mR^n \to \mR^n$, $T_{\vec{w}_{\ell}}(\vec{x}) = \vec{x} + \vec{w}_{\ell}$ be translation by $\vec{w}_{\ell}$.
Define the composition
$$\Psi_{\ell} \equiv (\psi^g_{y_{\ell}, \whv})^{-1} \circ \psi^g_{x, \whu} \circ T_{\vec{w}_{\ell}} \colon D(2d_2) \to D(\rF/2)$$
for an orthonormal frame $\whv = \whv_{\ell}$ of $T_{y_{\ell}}\F$ as provided by Proposition~\ref{prop-framed} so that $\| \Psi_{\ell}(\vec{x}) - \vec{x} \| \leq \ve_4 \, \rF$, where $\ve_4$ is defined by \eqref{eq-epsilon4a} and \eqref{eq-epsilon4b}. This is possible by our choice of $\ve_0$ in \eqref{eq-epsilon0}. (In this application of Proposition~\ref{prop-framed}, we take $y' = y_{\ell} \in \cP_i(x)$ which is on the same plaque as $x$.)

For $0 \leq j \leq m$, define $\vec{z}_j = (\psi^g_{y_{\ell} , \whv})^{-1}(y_j)$, and also set $\vwpj = \vec{w}_j - \vec{w}_{\ell}$. Then $\Psi_{\ell}(\vwpj) = \vec{z}_j$.
Using that $\Psi_{\ell}$ is $\ve_4 \rF$ close to the identity, we have that each $\| \vec{z}_j - \vwpj \| \leq \ve_4 \rF$.

Note that $\{\vec{w}_0 , \ldots , \vec{w}_{\ell}\}$ is $\rho_{\ell}$-robust if and only if the collection
$\{\vec{w}_0' , \ldots , \vec{w}_{\ell}' \}$ is $\rho_{\ell}$-robust.

By the definition that $\{y_0, \ldots , y_m\}$ is $\rho_0$-robust, the point $y_{{\ell}+1}$ lies at distance at least $\rho_0$ from the submanifold $\psi^g_{y_{\ell} , \whv}( {\rm Span}(\vec{z}_0, \ldots , \vec{z}_{{\ell}}))$.
So by Lemma~\ref{lem-geodesicisom}, the vector $\vec{z}_{{\ell}+1}$ lies at distance at least $\rho_0 - \ve_0 \rF \geq \rho_{\ell}$ from the linear span ${\rm Span}(\vec{z}_0, \ldots , \vec{z}_{\ell})$.
Thus, we also have that the collection
$\{\vec{z}_0, \ldots , \vec{z}_{{\ell}+1} \}$ is $\rho_{\ell}$-robust.

It is given that $e_1 \leq \dF(y_j, y_k) \leq e_2$ for $0 \leq j \ne k \leq m$,
so by Lemma~\ref{lem-geodesicisom} we have
$$d_1 < e_1 - 2\ve_0 \rF \leq \| \vec{z}_j - \vec{z}_k\| \leq 2 e_2 + 2\ve_0 \rF < 2 d_2$$
for all $0 \leq j \ne k \leq m$. We can thus apply Proposition~\ref{prop-robustvarest} for $\ve = \ve_4 \rF$
to the collection $\{\vec{z}_0, \ldots , \vec{z}_{\ell + 1}\}$ to conclude that
${\rm Span}(\vec{w}_0, \ldots , \vec{w}_{{\ell}+1})$ is $\rho_{{\ell}+1}$-robust.
\endproof

\section{Nice stable  transversals} \label{sec-existence}

Let $\fM$ be a minimal matchbox manifold, and  $x \in \fM$. Assume there is given  a connected compact subset  $K_x \subset L_x$ such that there is    $\wtK_x \subset \wtL_x$   such that  $\ds \iota_x \colon \wtK_x \subset \wtL_x \to L_x \subset \fM$
is injective with image $K_x$. That is, $K_x$ is a proper base as in Definition~\ref{def-LBB}. As in the proof of Theorem~\ref{thm-stableapprox}, we can assume that   $L_x$ is a leaf without holonomy.
Let $\ds \whK_x$ be the extension of $K_x$ as defined by \eqref{eq-translationK0}. We show there exists a clopen  set  $V_x \subset \fT_*$ which is $\whK_x$-admissible  so that $\ds  \fN(\whK_x, V_{x})$ is well-defined, and there exists 
a nice stable  transversal $\cX$ for $\ds  \fN(\whK_x, V_{x})$. The existence of this transversal, combined with Theorem~\ref{thm-stableapprox}, then yields the proof of Theorem~\ref{thm-tessel}.

 The strategy for the construction of $\cX$ is straightforward.    The procedure is inductive, in that we assume that a collection of transversals are given but not necessarily complete, and which satisfy a set of regularity conditions, and then prove that it is possible to extend the collection by adding another transversal so the regularity conditions are again satisfied. This process must terminate, as $\whK_x$ is compact, and we then have constructed a nice stable $(d_1, d_2)$-uniform transversal for the Reeb neighborhood $\ds \wtfN(\whK_x, V_{x})$. The difficulty with this approach, is       that  the ``stable condition''  must be satisfied    at each stage of the process. That is,  each partial collection of transversals for $\ds \wtfN(\whK_x, V_{x})$
  must satisfy appropriate net and general position conditions when intersected with the set $\whK_x$. As the leaf $L_x$ is dense in $\fM$, the choice of each successive additional transversal $\fZ(\xi_{k}, i_{\xi_{k}}, V_x)$ will induce   additional net points in $\whK_x$ which must be in general position and stable with respect to the previous choices. Ensuring that these partial stability conditions are satisfied requires delicate restraints on the successive choices, and the procedure for making these choices   utilizes the constants and estimates introduced in the  previous sections. The proof that the resulting transversal is nice and stable will be given in section~\ref{sec-nicestable}. 

\subsection{Basic notations}

 We first establish some simplified notations to be used in the construction. As noted above, we can assume that we are given a point $x_0 \in \fM$ contained in a leaf    $L_0 \subset \fM$  without holonomy, so  the holonomy covering map $\Pi \colon \wtL_0 \to L_0$   is a diffeomorphism.
  Let $\cM_0$ be an $(e_1, e_2)$-net for $L_0$ as in section~\ref{sec-microbundles}, where $e_2 = \eFU/4$.   Also, let   $\wtcM_0 = \Pi^{-1}(\cM_0) \subset \wtL_0$, with   $z = \Pi(\wtz) \in \cM_0$ for  $\wtz \in \wtcM_0$.
There is given  the covering of $\fM$ by  coordinate charts $\{U_{i_z} \mid z \in \cM_0\}$, as in the proof of Lemma~\ref{lem-newcover}, where  $1 \leq i_{z} \leq \nu$, 
and  such that  $B_{\fM}(z, \eU) \subset U_{i_z}$.  Correspondingly,  
for each $\wtz \in \wtcM_0$,  there is a foliation  chart $\wtU_{\wtz} = \oU_{i_z} \times \{\wtz\}$ for $\wtfN_0$.

 Let $K_0 \subset L_0$ be a  compact connected subset  which is a union of the plaques in $L_0$.
 Without loss of generality, assume there is a subset 
 $\cM_0' \subset \cM_0$ so that 
   $K_0 = \cup \ \{\cP_{i_z}(z)  \mid z \in \cM_0' \}$. Set $\wtK_0 = \Pi^{-1}(K_0)$ and   $\wtcM_0' = \Pi^{-1}(\cM_0') \subset \wtL_0$. 
 Of course, $L_0$ is without holonomy, so $\Pi \colon \wtL_0 \to L_0$ is a diffeomorphism, so this is just notational semantics. The distinction becomes important when considering the Reeb neighborhood of $\wtL_0$.

    Let $R_0$ denote the diameter of $K_0$ in the leafwise metric, so for any  $x \in K_0$ we have   $K_0 \subset D_{\F}(x, R_0)$.
Fix  a basepoint $z_0 \in  \cM_0'$. Without loss of generality, we may assume that $i_{z_0} = 1$,  and let $w_0 \in \fT_{1}$ be the projection of $z_0$ to the transverse space for $U_{1}$. 
For $z \in \cM_0$,  let $h_{z}$ denote the holonomy along a nice path $\gamma$ from $z_0$ to $z$, which  is well defined as  $L_0$ is without holonomy.

 Recall that the constant $\eFU$ is the ``leafwise Lebesgue number'' defined in equation \eqref{eq-leafdiam}, so that
for all $y \in \fM$, $D_{\F}(y, \eFU) \subset D_{\fM}(y, \eU/2)$. The constants $\ve_0$, $\ve_1$, $\ve_2$, $\ve_3$, $\ve_4$, and $\rF$ are as chosen in section~\ref{sec-constants}. 
Also,  $r_*$ was determined by the choice of  $\rF$ in section~\ref{subsec-r*}. 

For $R_0' = R_0 + \dFU$, let $\dT_0 = \delta( r_*/2, R_0')$ be the constant defined in  Proposition~\ref{prop-domest}. Then we have
\begin{equation}\label{eq-fundsets}
D_{\fX}(w_0, \dT_0)  \subset   D(h_{z}) \quad , \quad  h_{z}(D_{\fX}(w_0, \dT_0)) \subset  D_{\fX}(h_{z}(w_0), r_*/2)  \quad    {\rm for~all} ~ z \in \cM_0'
\end{equation}
 Let $V_0 \subset D_{\fX}(w_0, \dT_0)$ be a clopen subset with $w_0 \in V_0$. It follows from \eqref{eq-fundsets} that 
    $V_0$ is $K_0$-admissible, so we can form   the Reeb neighborhood $\ds \wtfN(\wtK_0, V_0) \subset  \wtfN_0$.     
 Note that even though the map $\wtL_0 \to L_0 \subset \fM$ is a diffeomorphism, the same is not necessarily true for the map 
  $\ds \wtfN(\wtK_0, V_0) \to \fM$.  For this reason, we make our constructions in the space   $\wtfN_0$.
  Notations involving the ``tilde'' indicate that the construction is considered in $\wtfN_0$ as opposed to $\fM$.

 For   $z\in \cM_0'$,  set $V_{z} = h_{z}(V_0) \subset \fT_{i_{z}}$ and $w_z = h_z(w_0)$. 
Then by the definition of the function $\delta(\epsilon, r)$ and the choice of $V_0 \subset D_{\fX}(w_0, \dT_0)$,  we have by  \eqref{eq-fundsets} that 
$\ds  V_{z}   \subset B_{\fX}(w_{z} , r_*/2) \subset \fT_{i_{z}}$. 
Hence, for any $w \in V_{z}$ we have $V_{z} \subset B_{\fX}(w, r_*)$.
This implies that   the various estimates introduced in sections~\ref{sec-approx} and \ref{sec-constants} for the transverse leafwise metric distortions will be satisfied for plaques in $\ds \wtfN(\wtK_0, V_0)$.

Introduce  the plaque-saturated compact sets, 
\begin{equation}
\fU^{V}_{z} = \pi_{z}^{-1}(V_{z}) \subset U_{i_{z}}
\quad , \quad
\wtfU^{V}_{\wtz} = \fU^{V}_{z} \times \{\wtz\} \subset \wtU_{i_{\wtz}} \subset \wtfN_0
\end{equation}
  so that $\ds \wtfN(\wtK_0, V_0)$ is the union of   $\wtfU^{V}_{\wtz}$ for $z \in \cM_0'$, and the transversals we introduce are subsets.  
For $x \in \fU_{z}^{V_z}$, define a \textit{standard sections} by
\begin{equation}\label{eq-localsection}
\fZ(x, i_{z}, V_z) = \vp_{i_{z}}^{-1}(\lambda_{i_{z}}(x), V_{z}) \subset \fU^{V}_{z}  
\quad , \quad
\wtfZ(x, i_{\wtz}, V_z) \equiv \fZ(x, i_{z}, V_z) \times \{\wtz\} \subset \wtfN_0
\end{equation}

 Finally,  introduce a sequence of constants based on the scale $\rF$:
\begin{equation}\label{eq-errors1}
\begin{array}{ccc}
d_1 = .10 \cdot \rF, & d_1' = .11 \cdot \rF, & d_1'' = .12 \cdot \rF\\
d_2 = .20 \cdot \rF, & d_2' = .19 \cdot \rF, & d_2'' = .18 \cdot \rF
\end{array}
\end{equation}
Note that $d_2 = 2d_1$ and
$$\rF/10 = d_1 < d_1' < d_1'' < d_2'' < d_2' < d_2 = \rF/5$$

  \subsection{The induction hypotheses}
  
The construction of a nice stable transversal $\cX$ proceeds  by induction from a given ``partial net'' in $K_0$ by choosing points which complete it to a $(d_1 , d_2)$-uniform  net satisfying the stability hypotheses.  The idea is to formulate the notion of a \emph{regular partial  $V_0$-transversal} for $K_0$, which satisfies somewhat stronger hypotheses than are eventually required. 

Assume   as given a set of points   $\ds \Xi_p = \{\xi_1, \ldots, \xi_{p}\} \subset K_0$ and 
  $\ds \Lambda_p = \{z_1, \ldots , z_{p}\} \subset \cM_0'$ such that $\xi_j \in B_{\F}(z_j , \eFU/2)$.
  Label the   indices $\ds \theta_j = i_{z_j}$ for $1 \leq j \leq p$ where  $z_{\theta_j} \in \cM_0'$.

Set  $\ds \cX_j = \fZ(\xi_j, \theta_j , V_0)$ for $1 \leq j \leq p$, and define the partial $V_0$-transversal 
 $\whcX_p = \cX_1 \cup \cdots \cup \cX_{p}$. 
 
 The conditions imposed on  $\whcX_p$ are all leafwise.  For $x \in \fM$  introduce also the intersections
 \begin{equation}\label{eq-leafnets}
\cN_p(x) = \whcX_p \cap L_x \quad , \quad \cN_p = \whcX_p \cap L_0
\end{equation}

 The first hypothesis is the   separation condition  $\dF(y,z) \geq d_1'$  for all $y \ne z\in \cN_p(x)$.
 Note that $d_1' > d_1$ so a sufficiently small perturbation of the net will still satisfy the $d_1$-separation condition.

Also note that if $p=1$ and 
$\xi \ne \xi' \in \cN_1$, then $\xi, \xi' \in \cX_1$ implies they are contained   in distinct plaques of $L_0$,
hence $\dF(\xi, \xi') \geq 2\dFU > d_1'$. 

In addition, impose two additional stability conditions criteria  on the  nets $\cN_p(x)$, that they satisfy Definitions~\ref{def-ve1reg} and \ref{def-ve2robust} below, which only apply when  $p \geq 2$.

Let $\Delta'_{\F}(\whcX_p)$ denote the   \emph{$d_2'$-bounded} leafwise simplicial complex  for $\whcX_p$. By definition, a  $(k+1)$-tuple $\{y_0, \ldots , y_k\} \subset \whcX_p$ defines a $k$-simplex $\Delta(y_0 , \ldots, y_k) \in \Delta'_{\F}(\whcX_p)$ if
there exists $\omega \in L_{y_0}$ and $0 < r \leq d_2'$ such that
$B_{\F}(\omega,r) \cap \whcX_p = \emptyset$ and
$\{y_0, \ldots , y_k\} \subset S_{\F}(\omega, r) \cap \whcX_p$.

\begin{defn} \label{def-ve1reg}
The transversal $\whcX_p$  for $K_0$ 
is \emph{$\ve_1$-regular} if for all   $\Delta(y_0 , \ldots, y_n) \in \Delta'_{\F}(\cX_p)$, 
 and  for all  $\xi \in \cN_p(y_0) - \{y_0 , \ldots, y_n\}$, then 
\begin{equation}\label{eq-unifreg}
\dF(\xi , \omega(y_0 , \ldots, y_n)) ~ \geq ~ r(y_0 , \ldots, y_n) +  \ve_1 \rF
\end{equation}
\end{defn}

Given a simplex $\Delta(y_0 , \ldots, y_k) \in \Delta'_{\F}(\whcX_p)$, we say that the vertices are \emph{properly ordered} if
there exists $1 \leq i_0 < i_1 < \cdots < i_k \leq p$
and points $\xi_{i_{j}} \in \whcX_p \subset L_{y_0}$ such that
$y_{j} = \cX_{i_{j}} \cap \cP_{\theta_{i_{k}}}(y_k)$.

\begin{defn} \label{def-ve2robust}
Let $\rho = 3\ve_2\rF/2$.
A regular partial   $V_0$-transversal $\whcX_p$  for $K_0$ 
is \emph{$\rho$-robust} if for all  $\Delta(x_0 , \ldots, x_n) \in \Delta'_{\F}(\cN_p)$, 
the collection $\{x_0, \ldots , x_n\}$ is $\rho$-robust, in the sense of Definition~\ref{def-robustnets}, where we 
assume the vertices $\{x_0 , \ldots , x_n\}$ are properly ordered.
\end{defn}
 
 Finally, we give a condition for when we are done with the inductive construction.
 \begin{defn}
A regular partial $V_0$-transversal $\whcX_p$ for $K_0$ is \emph{$\delta$-complete} if
$$K_0 ~ \subset ~  \PF(\cN_p, \delta) ~ = \bigcup_{z \in \cN_p} ~ D_{\F}(z, \delta)$$
\end{defn}

Normalize the first section  $\cX_1$ by defining
\begin{equation}\label{eq-netL0}
\xi_1 = z_0 \in K_0 ~ , ~ \Xi_1 = \{\xi_1\} ~  , ~ \Lambda_1 = \{z_1\} ~ , ~ \theta_1 = i_{z_1} = 1 ~, ~ \cX_1 = \fZ(\xi_1, \theta_1 , V_0)
\end{equation}

\subsection{The inductive construction} The proof of the next result gives the inductive construction.  
\begin{prop}\label{prop-induct}
Let $\whcX_p$ for $p \geq 1$ be a regular partial $V_0$-transversal satisfies the conditions of Definitions~\ref{def-ve1reg} and \ref{def-ve2robust}. 
If $\whcX_p$ is not $d_2''$-complete for $K_0$ then there exists 
$\xi_{p+1}   \in K_0$ so that for
\begin{itemize}
\item $\ds \Xi_{p+1} = \{\xi_1, \ldots, \xi_{p+1}\} \subset L_0$ 
\item $\ds \Lambda_{p+1} = \{z_1, \ldots , z_{p}, z_{p+1}\} \subset \cM_0'$ such that $\xi_j \in B_{\F}(z_j , \eFU/2)$
\item   $\ds \theta_j = i_{z_j}$ for $1 \leq j \leq p+1$ with $z_{\theta_j} \in \cM_0'$ 
\item  $\cX_{p+1} = \fZ(\xi_{p+1}, \theta_{p+1} , V_0)$ ~ , ~ $\ds \whcX_{p+1} = \cX_1 \cup \cdots \cup \cX_p \cup \cX_{p+1}$
\end{itemize}
  then  $\ds \whcX_{p+1}$  is $d_1'$-separated, and satisfies the conditions of Definitions~\ref{def-ve1reg} and \ref{def-ve2robust}. 
\end{prop}
\proof
Set   $\ds \cN_{p} = \whcX_{p} \cap L_{0}$. As $\whcX_p$ is not $d_2''$-complete,  we have  $\ds K_0 -  \PF(\cN_p, d_2'') \ne \emptyset$.
The set $K_0$ is a union of plaques of $\F$  which are convex, so  there exists  $\xi_{p+1}'  \in K_0 $  such that 
\begin{equation} \label{eq-netp2}
B_{\F}(\xi_{p+1}' , \rF/200) \subset  K_0 \cap \left\{    \PF(\whcX_p, d_2'') -  \PF(\whcX_p, d_1'') \right\}.
\end{equation}
Then for all  $z \in \whcX_p$ we have 
  $\dF(\xi_{p+1}' , z) > d_1''$.
Choose
$z_{p+1} \in \cM_0 \cap   B_{\F}(\xi_{p+1}' , \eFU/2)$,
which is possible by the assumption that $\cM_0$ is a net which is $\eFU/2$-dense.
Set $\theta_{p+1} = i_{z_{p+1}}$. Adding the $V_0$-transversal $\fZ(\xi_{p+1}', \theta_{p+1} , V_0)$ 
will result in a   $d_1'$-separated $V_0$-transversal $\whcX_{p+1}'$.
However, such $\whcX_{p+1}'$ need not satisfy the conditions of Definitions~\ref{def-ve1reg} and \ref{def-ve2robust}. To ensure that these conditions also hold,  we  modify the choice of $\xi_{p+1}'$ to a point $\xi_{p+1} \in B_{\F}(\xi_{p+1}' , \rF/200)$.

 Consider the disk
$D_{\F}(\xi_{p+1}', 4d_2) \subset B_{\F}(\xi_{p+1}', \rF) \subset B_{\F}(\xi_{p+1}', \lF)$.
Introduce the set
\begin{equation}\label{eq-localsets}
\Omega(\xi_{p+1}') = D_{\F}(\xi_{p+1}', 4d_2) \cap \whcX_p
\end{equation}
Since the points of $\cN_p(y_0)$ are $d_1'$-separated, and $d_2 = 2 d_1$, the metric conditions
\eqref {eq-dgs1}, \eqref {eq-dgs2} and \eqref {eq-dgs3} and a standard volume estimate yields that the
cardinality of $\Omega(\xi_{p+1}')$ is at most $10^n$.

Let $\Omega^{(n)}(\xi_{p+1}') \subset \Delta^{(n)}_{\F}(\whcX_p)$ be the subset of all $n$-simplices whose vertices are contained in $\Omega(\xi_{p+1}')$.
The cardinality of the set $\Omega^{(n)}(\xi_{p+1}')$ is thus bounded above by the constant
$\ds C_n$ defined in \eqref{eq-Cn}.

For each $n$-simplex $\Delta(y_0 , \ldots, y_n) \in \Omega^{(n)}(\xi_{p+1}')$ recall that
$\omega(y_0 , \ldots, y_n) \in D_{\F}(\xi_{p+1}', 4d_2)$ denotes the center of the circumscribed sphere for its vertices, so $\{ y_0 , \ldots, y_n \} \subset S_{\F}(\omega(y_0 , \ldots, y_n), r(y_0 , \ldots, y_n))$.
For a constant $\kappa > 0$,
form the annular region
\begin{equation}\label{eq-annularregion}
A_{\F}\left( y_0 , \ldots, y_n; \kappa \right) = \PF(S_{\F}(\omega(y_0 , \ldots, y_n), r(y_0 , \ldots, y_n)), \kappa)
\end{equation}

\begin{lemma}\label{lem-volest3} Let $\kappa = 2 \ve_1 \rF$ then 
\begin{equation}\label{eq-volest2}
{\rm Vol}_{\F}~A_{\F}\left( y_0 , \ldots, y_n; \kappa \right)
\leq 200 \cdot 2^n ~ \ve_1 (\rF/5)^n
\end{equation}
\end{lemma}
\proof
Let $\Phi_n$ be the constant such that $\ds {\rm Vol}_{\whu}~D_{\mR^n}(s) = \Phi_n s^n$, where ${\rm Vol}_{\whu}$ denotes the volume with respect to the frame $\whu$.  Note that $(\sqrt{2})^{n} \leq \Phi_n \leq 2^n$.

By the condition \eqref {eq-dgs3} for $0 < s \leq \rF \leq \lambda_{\ve_0}$ we have
$$\left | \Phi_n s^n - {\rm Vol}_{\F}~D_{\F}\left( \omega(y_0 , \ldots, y_n), s \right) \right | \leq \ve_0 \cdot s^n
\leq \ve_0 \cdot (\rF)^n$$
Hence, we have
\begin{eqnarray}
\lefteqn{ {\rm Vol}_{\F}~A_{\F}\left( y_0 , \ldots, y_n ; \kappa \right) \quad } \nonumber\\
& = & ~ \left | {\rm Vol}_{\F}~D_{\F}\left( \omega(y_0 , \ldots, y_n), r(y_0 , \ldots, y_n) + \kappa \right) -
{\rm Vol}_{\F}~D_{\F}\left( \omega(y_0 , \ldots, y_n), r(y_0 , \ldots, y_n) - \kappa \right) \right | \nonumber \\
& \leq & ~ \Phi_n \cdot \left\{ (r(y_0 , \ldots, y_n)+ \kappa)^n - (r(y_0 , \ldots, y_n)- \kappa)^n \right\} + 2\ve_0 (\rF)^n \nonumber
\end{eqnarray}

Given $\kappa = 2 \ve_1 \rF$ with $20n \ve_1 < 1$, and $\rF/10 \leq r \leq \rF/5$, elementary estimates yield
\begin{eqnarray*}
\left\{ (r + \kappa)^n - (r - \kappa)^n \right\} ~ & = & ~ r^n \cdot \left\{ (1 + \kappa / r)^n - (r - \kappa / r)^n \right\} \\
~ & \leq & ~ r^n \cdot \{ (\exp (n \kappa / r) - \exp ( - n \kappa / r) \} \\
~ & \leq & ~ (\rF/5)^n \cdot \{ (\exp (20 n ~ \ve_1) - \exp ( - 20 n ~ \ve_1) \} \\
~ & \leq & ~ 100 n ~ \ve_1 (\rF/5)^n
\end{eqnarray*}
Combining these estimates and    condition 2 in \eqref{eq-epsilon0}, we obtain
\begin{eqnarray}
{\rm Vol}_{\F}~A_{\F}\left( y_0 , \ldots, y_n, \kappa \right)
& \leq & ~ \Phi_n \cdot \left\{ (r(y_0 , \ldots, y_n)+ \kappa)^n - (r(y_0 , \ldots, y_n)- \kappa)^n \right\} + 2\ve_0 (\rF)^n \nonumber \\
& \leq & ~ \{ \Phi_n \cdot 100 n ~ \ve_1 + 2 \cdot 5^n ~\ve_0 \} (\rF/5)^n \nonumber \\
& \leq & ~ \{ 2^n \cdot 100 n ~ \ve_1 + 2 \cdot 5^n ~ \ve_0 \} (\rF/5)^n \nonumber \\
& \leq & ~ 200n \cdot 2^n ~ \ve_1 (\rF/5)^n \label{eq-volestep0} \quad \quad \quad \quad \quad \quad \quad \quad \quad   \Box
\end{eqnarray}

The total volume of all such annular regions intersecting $D_{\F}(\xi_{p+1}' , \rF/100)$ is bounded above by
\begin{align*}
C_n \cdot 200n \cdot 2^n ~ \ve_1 (\rF/5)^n
\end{align*}
We also derive an estimate of the leafwise volume of the disk $\ds D_{\F}(\xi_{p+1}', \rF/200)$. Note that 
$$ | \Phi_n \cdot (\rF/200)^n - {\rm Vol}_{\F}~D_{\F}(\xi_{p+1}', \rF/200) | \leq \ve_0 \cdot (\rF/200)^n = (\ve_0/2^n) \cdot (\rF/100)^n$$
so that
\begin{equation}
{\rm Vol}_{\F}~D_{\F}(\xi_{p+1}', \rF/200) ~ \geq ~ \Phi_n \cdot (\rF/200)^n - \ve_0 \cdot (\rF/200)^n ~ \geq ~ (1/40)^n \cdot (\rF/5)^n
\end{equation}

Now, given $\ds \ve_1 = 1/(C_n \cdot 1000 n \cdot 100^n)$ by \eqref{eq-epsilon1}, it follows that
\begin{equation}
C_n \cdot 200n \cdot 2^n ~ \ve_1 \cdot (\rF/5)^n ~ \leq ~ \frac{1}{4} \cdot 1/40^n \cdot (\rF/5)^n
\end{equation}
Thus, the total volume of all annular regions intersecting
$\ds D_{\F}(\xi_{p+1}', \rF/200)$ is less than 1/4 of its volume.
Therefore, if we choose 
$\xi_{p+1} \in B_{\F}(\xi_{p+1}' , \rF/200) \subset K_0$ which lies outside of the union of these annular regions, then for all $\Delta(y_0 , \ldots, y_n) \in \Delta^{(n)}_{\F}(\cN_p)$
\begin{equation}\label{eq-unifreg2}
\dF(\xi_{p+1} , \omega(y_0 , \ldots, y_n)) ~ \geq ~ r(y_0 , \ldots, y_n) + 2\ve_1 \rF
\end{equation}
 Note that this estimate on $L_0$ is  stronger  than  \eqref{eq-unifreg}, and it will be shown in Lemma~\ref{lem-separated}   that \eqref{eq-unifreg} holds  for all $x \in \fM$.

It remains to  modify the choice of $\xi_{p+1}$ so that  the robustness condition Definition~\ref{def-robustnets} is also satisfied.
The strategy is to again use volume estimates, in this case for the sets of points for which the robustness condition fails, then chose $\xi_{p+1}$ in the complement. 
   
For $1 \leq k < n$, let $\{y_0 , \ldots, y_{k}\} \subset \Omega(\xi_{p+1}')$ be a collection of distinct points with $y_k \in \cX_{i_k}$ where $1 \leq i_0 < \cdots < i_{k} \leq p$.
Let $\whu = \{\vec{u}_1, \ldots, \vec{u}_n\} \subset T_{y_k}\F$ be an orthonormal frame, and introduce the corresponding geodesic coordinates
\begin{equation}
\psi^g_{y_k, \whu} \colon D(\rF) \to D_{\F}(y_k , \rF) \subset L_0
\end{equation}
Define $\vec{y}_j = (\psi^g_{y_k, \whu})^{-1}(y_j)$ for $0 \leq j \leq k+1$. Note that $\vec{y}_{k} = \vec{0}$.

   Let ${\rm Span}(\vec{y}_0, \ldots , \vec{y}_{k}) \subset \mR^n$ be the linear submanifold through the origin of dimension $k$ which they span. Then define a submanifold of $D_{\F}(\xi_{p+1}' , \rF/200)$,
\begin{equation}\label{eq-subplane}
H(y_0 , \ldots, y_{k}; \xi_{p+1}') = \psi^g_{y_k, \whu}\left\{ {\rm Span}(\vec{y}_0, \ldots , \vec{y}_{k}) \cap D(2d_2) \right\} \cap D_{\F}(\xi_{p+1}' , \rF/200)
\end{equation}
which has diameter at most $\rF/100$, and thus has $(n-1)$-volume bounded above by $(\rF/100)^{n-1}$.  

Form the $2\ve_2 \rF$-thickening of $H(y_0 , \ldots, y_{k}; \xi_{p+1}')$,
\begin{equation}\label{eq-slabs}
\cS(y_0 , \ldots, y_{k}; \xi_{p+1}', 2\ve_2 \rF) = \PF(H(y_0 , \ldots, y_{k}; \xi_{p+1}'), 2\ve_2 \rF) \cap D_{\F}(\xi_{p+1}', \rF/200)
\end{equation}

Then by the estimate \eqref{eq-dgs3a} and $\ve_0 \leq \ve_2$, its volume is bounded above by
\begin{equation}\label{eq-slabsest}
4 (\ve_2 \rF) \cdot (\rF/100)^{n-1} + \ve_0 (\rF/100)^n \leq 5 (\ve_2 \rF ) \cdot (\rF/100)^{n-1}
\end{equation}

The total number of such submanifolds $H(y_0 , \ldots, y_{k}; \xi_{p+1}')$ in $D_{\F}(\xi_{p+1}', 4d_2)$ is bounded above by the constant $C_n$ from \eqref{eq-Cn}, hence the total
volume of all such sets which intersect $D_{\F}(\xi_{p+1}', \rF/200)$ is thus bounded above by
\begin{equation}\label{eq-slabvol}
C_n \cdot 5( \ve_2 \rF) \cdot (\rF/100)^{n-1} =
\frac{C_n \cdot (5 \rF) \cdot (\rF/100)^{n-1}}{C_n \cdot 2000 \cdot 2^n} =
\frac{1}{4} \cdot 1/40^n \cdot (\rF/5)^n
\end{equation}
where we use the definition of $\ve_2$ in \eqref{eq-epsilon2}.

Thus, the total volume of all such slabs intersecting
$\ds D_{\F}(\xi_{p+1}', \rF/200)$ is less than 1/4 of its volume,
so we may choose $\xi_{p+1} \in B_{\F}(\xi_{p+1}', \rF/200)$ which is disjoint from the union of all annular and slab regions introduced above. 
This completes the choice of the new point $\xi_{p+1}$.

We now must check that all of the required hypotheses for a nice stable $V_0$-transversal are satisfied.
First, we note    that $\whcX_{p+1}$ is $d_1'$-separated.  Set $ \cN_{p+1}(x)  = \whcX_{p+1} \cap L_x$ and  $ \cN_{p+1} = \whcX_{p+1} \cap L_0$. 

\begin{lemma}\label{lem-separated}
For all   $y \ne z \in \whcX_{p+1}$ we have $\dF(y,z) \geq .114 \cdot \rF > d_1'$.
\end{lemma}
\proof
If $y,z$ lie on distinct leaves, there is nothing to show. Assume     $y \ne z \in \cN_{p+1}(x)$. Then by definition, there exists $\xi_i , \xi_j \in \whcX_{p+1}$ for $1 \leq i, j \leq p+1$ such that
$y \in \fZ(\xi_{i}, \theta_i , V_0) \cap L_x$ and
$z \in \fZ(\xi_{j}, \theta_j , V_0) \cap L_x$.
Without loss of generality we can assume that $i \geq j$.

If $z \not\in \cP_{\theta_{i}}(y)$ then $\dF(y,z) \geq \dFU > \rF > d_1''$.
Thus, we may assume that $z \in \cP_{\theta_{i}}(y)$, so $i > j$.

Set $y' = \fZ(\xi_{i}, \theta_{i} , V_0) \cap \cP_{\theta_{i}}(\xi_{i}) = \xi_i$ and
set $z' = \fZ(\xi_{j}, \theta_j , V_0) \cap \cP_{\theta_{i}}(\xi_{i})$.

Then $\dF(y', z') \geq d_1'' - \rF/200 = .115 \cdot \rF$ by the choice of $\xi_i'$ satisfying \eqref{eq-netp2} and the choice of $\xi_i$.

Apply Lemma~\ref{lem-secdiv} for the pairs $\{y, z\}$ and $\{y', z'\}$ and note that
$2\ve_0 < 1/1000$ to obtain
\begin{equation}\label{eq-netuniform1}
\dF(y,z) ~ \geq ~ \dF(y', z') - 2\ve_0 \cdot \rF > .115 \cdot \rF - .001 \cdot \rF = .114 \cdot \rF
\end{equation}
Thus, for all $x \in \fM$ the net $\cN_{p+1}(x)$ is $(.114 \cdot \rF)$-separated.
\endproof

A simple consequence of Lemma~\ref{lem-separated} is that if $\whcX_p$ is uniformly $d_1'$-separated, then the collection of leafwise disks
$\ds \{D_{\F}(\xi_{\ell} , d_1'/2) \mid \xi_{\ell} \in \Xi_p\}$ 
are pairwise disjoint. As the set $K_0$ was assumed to be compact, this implies the cardinality of the set $\Xi_p$ has an a priori bound. 
That is, we can repeat the construction in Proposition~\ref{prop-induct} at most a finite number of times, until we obtain 
a regular partial $V_0$-transversal $\whcX_{p_*}$ for $K_0$ which is   $d_2''$-complete, for some $p_* > 0$. 

\medskip

Set  $\cX = \whcX_{p_*}$, set $\cN = \cX \cap L_0$  and for $x \in \fM$,  let $\cN(x) = \cX \cap L_x$.

\begin{prop}\label{prop-dense}
For $x \in \fM$ the set $\cN(x)$ is a $(d_1' , d_2')$-net for $\fN(K_0, V_0) \cap L_x$.
\end{prop}
\proof
 The set $\cN(x)$ is $d_1'$-separated   by Lemma~\ref{lem-separated}.
 
 We must show that $\cN(x)$   is $d_2'$-dense in   $\ds \fN(K_0, V_0)$. For    $y \in \fN(K_0, V_0)$ we show there exists $z \in \cN(y)$ such that $\dF(y,z) \leq .181 \cdot \rF < d_2'$.

Let $y \in \fN(K_0, V_0)$. Then by the constructions in section~\ref{sec-existence}, there exists 
$z \in \cM_0'$ for which  $y = \cP_{i_z}(y) \cap \fZ(y, i_z , V_0)$. 
Choose $\zeta \in \fU^{V_0}_{z} \cap K_0$ and set $y' = \cP_{i_{z}}(\zeta) \cap \fZ(y, i_z , V_0) \in K_0$.

We are given that $K_0 \subset \PF(\cN , d_2'')$, so there exists $\xi \in \cN$ such that $\dF(y' , \xi) \leq d_2'' = .18 \cdot \rF$.

By definition of $\cN$, there exists $\xi_j \in \Xi$ for some $1 \leq j \leq p_*$ such that
$\xi = \cP_{\theta_i}(\xi) \cap \fZ(\xi_{j}, \theta_i , V_0)$.

Let $z = \cP_{\theta_i}(y) \cap \fZ(\xi_{j}, \theta_i , V_0) \in \cN(y)$, and
apply Lemma~\ref{lem-secdiv} for the pairs $\{y, z\}$ and $\{y' , \xi \}$ to obtain
\begin{equation}\label{eq-netuniform2}
\dF(y,z) ~ \leq ~ \dF(y', \xi) + 2\ve_0 \, \rF \leq d_2'' + .001 \cdot \rF = .181 \cdot \rF
\end{equation}
Thus $\cX$ is $.181 \cdot \rF$-dense in $\fN(K_0, V_0)$. 
\endproof

 \section{Stability of parametrized Delaunay triangulations} \label{sec-nicestable}
 
The transversal $\cX$ is  a $(d_1' , d_2')$-net for $\fN(K_0, V_0)$ by   Proposition~\ref{prop-dense}. Moreover,  for each of the finite set of points in  $\Xi_{p_*} = \{\xi_1, \ldots, \xi_{p_*}\} \subset L_0$, the collection of $n$-simplices $\Omega^{(n)}(\xi_k)  \subset \Delta^{(n)}_{\F}(\whcX)$
contained in the disk $\ds D_{\F}(\xi_k, 4d_2)$   satisfy the regularity and robustness conditions  of Definitions~\ref{def-ve1reg} and \ref{def-ve2robust}. It remains to show that these conditions are satisfied  for any $x \in \fM$ and   all $n$-simplices lying in $\ds D_{\F}(x, 4d_2)$. The constants of section~\ref{sec-constants} were chosen so that this will be true, although the proofs of this assertion are rather involved. This will yield the transverse stability of the leafwise Delaunay triangulations $ \ds \Delta_{\F}(\cN(x))$.

At first inspection, the stability of simplices in $\ds \Delta_{\F}(\cN(x))$ for $x \in \cX$ appears to be ``intuitively clear'', as the lengths of the edges change continuously with $x$. In fact,  this is basically correct for dimension $n \leq 2$. The difficulty is that for $n > 2$, as $x$ varies, the ``small variations'' of the points of $\cN(x)$ may result in an abrupt change in the Delaunay simplicial structure, if any face of a Voronoi cell has too small of a diameter relative to the size of the variation. In the literature for Voronoi tessellations of $\mR^n$, this difficulty appears to be formulated as a ``conditioning'' criteria. In our context, of a varying Riemannian metric, we show the nets $\cN(x)$ are ``well-conditioned'' as $x$ varies,  using 
Proposition~\ref{prop-inductrobust} along with Propositions~\ref{prop-inductsphere} and \ref{prop-inductsimplex} below.

\subsection{Robustness} We  show that $\cX$ satisfies  Definition~\ref{def-ve2robust}. It suffices to    show     this robustness condition is stable, as it holds for the simplices in $\Omega^{(n)}(\xi_k)$ for $1 \leq k \leq p_*$ by construction. The proof   uses an induction procedure which invokes Proposition~\ref{prop-robustnets} repeatedly, and invokes the constants defined by \eqref{eq-epsilon4h} and $\ve_4$ as derived from the inequalities \eqref{eq-epsilon4a} and \eqref{eq-epsilon4b}. 
For $1 \leq \ell \leq n$, recall the definitions of $\whrho_{\ell}$ and $\whrhop_{\ell}$ from 
\eqref{eq-epsilon4a}, \eqref{eq-epsilon4b}, 
and set $\wtrho_{\ell} = \whrho_{\ell} \rF/10$ and $\wtrhop_{\ell} = \whrhop_{\ell} \rF/10$.

\begin{prop}\label{prop-inductrobust} Let $\rho = 3\ve_2\rF/2$. Then for each 
$\Delta(x_0 , \ldots, x_n) \in \Delta'_{\F}(\whcX_{p+1})$, such that the vertices $\{x_0 , \ldots , x_n\}$ are properly ordered, then 
 the collection $\{x_0, \ldots , x_n\}$ is $\rho$-robust.
\end{prop}
\proof
We proceed via induction on $1 \leq m < n$. We recall some notation. Let $(y_0,\ldots,y_n) \subset L_0$
such that $\Delta(y_0 , \ldots, y_n) \in \Delta'_{\F}(\whcX_{p+1})$. By permuting the order of the vertices,
we can assume that the vertices $\{x_0 , \ldots , x_n\}$ are properly ordered. That is, 
there exists $1 \leq i_0 < i_1 < \cdots < i_n \leq p+1$
and points $\xi_{i_{k}} \in \whXi_{p+1} \subset L_0$ such that
$y_{k} = \cX_{i_{k}} \cap \cP_{\theta_{i_{n}}}(y_n)$. For notational convenience, set $\cP_{n}(z) = \cP_{\theta_{i_{n}}}(z)$. 

The subtlety of the proof lies in the fact that the robust condition in Definition~\ref{def-robustnets} is with respect to the geodesic coordinates about the last vertex in the collection of properly ordered points, but the inductive hypotheses are in terms of the geodesic coordinates about each successive vertex, not just the last one. The change of coordinates from one vertex to another introduces   an error in the robust condition. Consequently, at each stage of the induction, the robust constants $\wtrho_{\ell}$ decrease to account for this error. This fact is behind the arcane definition in formula \eqref{eq-epsilon4h}. 

The first step of the induction, $m=1$, is trivial. Note that $\wtrho_0 = 18\ve_2 \rF/10$.
 Then given $\{x_0 , x_1\} \subset \cP_{1}(\xi_{i_1})$ as above, with $x_{\ell} \in \cX_{i_{\ell}}$ and $i_0 \ne i_1$ then $\dF(x_1, x_0) \geq d_1'$ by Lemma~\ref{lem-separated}. Hence
$$d_1' \geq 2\ve_2\rF > \wtrho_0 > \wtrho_1$$ and so  $\{x_0 , x_1\}$ is $\wtrho_1$-robust.

Now assume that $1 < m < n$.
We make an inductive hypothesis which is uniform for all simplices. That is, for fixed $m < n$, assume that for all $\Delta(y_0 , \ldots, y_n) \in \Delta'_{\F}(\whcX_{p+1})$, then for all subsets of points $\{x_0, \ldots , x_m\}$ defined for $x_n \in \cX_{i_{m}}$ as above, the set $\{x_0, \ldots , x_m\}$ is $\wtrho_m$-robust.
We then show that each transverse translate $\{x_0, \ldots , x_{m+1}\}$  of $\{y_0, \ldots , y_{m+1}\}$ is $\wtrho_{m+1}$-robust.

Consider first the case $z_{m+1} = \xi_{i_{m+1}} \in \cX_{i_{m+1}}$, and set $z_k = \cX_{i_{k}} \cap \cP_{n}( \xi_{i_{m+1}})$ for $0 \leq j \leq n$.
By the inductive hypothesis, the set
$\{z_0 , \ldots , z_m\}$ is $\wtrho_m$-robust. We verify the conditions of Definition~\ref{def-robustnets}
for the vertex, $z_{m+1}$.
The point $ \xi_{i_{m+1}}$ was chosen to that it lies outside of all $2\ve_2 \rF$-neighborhoods as defined in \eqref{eq-slabs} of the images under the exponential map of affine subspaces spanned by local collections of at most $n+1$ points. It follows, in particular, that the distance from $ \xi_{i_{m+1}}$ to the submanifold
$\ds H(z_0 , \ldots, z_{m}; z_m)$
in Definition~\ref{def-robustnets}.3 is at least $2\ve_2 \rF > \wtrho_m$. Thus,
$\{z_0 , \ldots , z_{m+1}\}$ is also $\wtrho_m$-robust.

Note that $d_1' \leq \dF(z_j, z_k)$ for $0 \leq j \ne k \leq n$ by Lemma~\ref{lem-separated}.
We are given $\Delta(y_0 , \ldots, y_n) \in \Delta'_{\F}(\whcX_{p+1})$ which implies that the vertices $\{y_0 , \ldots, y_n\}$ admit a circumscribed sphere, which must have radius at most $.181 \cdot \rF$ by Proposition~\ref{prop-dense}.
Thus, $\dF(y_j , y_k) \leq .362 \cdot \rF$. The map $\phi_{y_n , z_n}$ is an $\ve_0 \rF$-isometry by Lemma~\ref{lem-secdiv}, so $\dF(z_j , z_k) \leq .362 \cdot \rF + \ve_0 \cdot \rF < .380 \cdot \rF = 2d_2'$.
It follows that the set of points $\{z_0, \ldots , z_n\}$ satisfy the hypotheses of
Proposition~\ref{prop-robustnets} for $e_1 = d_1'$, $e_2 = d_2'$, and $\rho = \wtrho_m$.

For simplicity, set $\zeta = z_{m+1}$, and choose an orthonormal frame $\whu$ of $T_{\zeta}\F$.
Then for $0 \leq j \leq m$, let $\vec{z}_j = (\psi^g_{\zeta, \whu})^{-1}(z_j)$ for the geodesic coordinates
$\psi^g_{\zeta, \whu} \colon D(\rF) \to D_{\F}(\zeta,\rF)$.
Then the collection $\{ \vec{z}_0, \ldots , \vec{z}_{m+1}\} \subset \mR^n$
is $\wtrhop_m$-robust by Proposition~\ref{prop-robustnets} and the choice of $\ve_4$ in \eqref{eq-epsilon4a} and \eqref{eq-epsilon4b}.

The robustness for the set $\{ \vec{z}_0, \ldots , \vec{z}_{m+1}\}$ is used to show it for $\{x_0, \ldots , x_{m+1}\}$.
Set $\zeta' = x_{m+1}$, then
by Proposition~\ref{prop-framed}, there exists an orthonormal framing $\whv$ of $T_{\zeta'}\F$ so that the composition
$$\Psi_{\zeta, \zeta'} \equiv (\psi^g_{ \zeta' ,\whv})^{-1} \circ \phi_i(\zeta, \zeta') \circ \psi^g_{\zeta ,\whu} \circ T_{\zeta} \colon D(\rF/2) \to \mR^n \cong T_{\zeta'}\F$$
is $\ve_4 \rF$-close to the identity.
Set $\vec{w}_j = (\psi^g_{ \zeta' ,\whv})^{-1}(x_k) = \Psi_{\zeta, \zeta'}(\vec{z}_j)$, then $\|\vec{w}_j - \vec{z}_j\| \leq \ve_4 \rF$.

The set $\{\vec{w}_1 , \ldots , \vec{w}_{m+1}\}$ satisfies the hypotheses of Proposition~\ref{prop-robustvarest} for
$e_1 = d_1$, $e_2 = d_2$, $\ve = \ve_4 \rF$ and $\rho = \wtrhop_{m}$. Therefore, $\{\vec{w}_1 , \ldots , \vec{w}_{m+1}\}$ is
$(\wtrho_{m+1}+ \ve_2 \rF/1000)$-robust.

Finally, by Lemma~\ref{lem-geodesicisom} the geodesic map $\ds \psi^g_{\zeta', \whv}$ is an $\ve_0 \rF$-isometry,
hence the distance from $x_{m+1}$ to the submanifold $\ds H(x_0 , \ldots, x_{m}; x_m)$ in Definition~\ref{def-robustnets}.3 is at least $ \wtrho_{m+1} + \ve_2 \rF/1000 - \ve_0 \rF > \wtrho_{m+1}$.

This completes the inductive step.
It remains to note that $\wtrho_n > 3\ve_2 \rF/2$ by definition \eqref{eq-epsilon4h}.
\endproof

\subsection{Circumscribed spheres} \label{subsec-circumscribed}
The next step towards showing that $\cX$ is regular and stable is to show that the circumscribed    sphere condition is stable.

Let $\Delta(y_0 , \ldots, y_n) \in \Delta'_{\F}(\cX)$.
By permuting the order of the vertices,
we can assume that there exists $1 \leq i_0 < i_1 < \cdots < i_n \leq p_*$
and points $\xi_{i_{k}} \in \Xi_{p_*} \subset L_0$ such that
$y_{k} = \cX_{i_{k}} \cap \cP_{\theta_{i_{n}}}(y_n)$.

For $x_n \in \cX_{i_n} \subset \fU_{\theta_{i_{n}}}^{V_0}$ let $ \cP_{n}(x_n) = \cP_{\theta_{i_{n}}}(x_n)$
denote the plaque containing $x_n$ in the chart $\vp_{i_n}$. Then 
set $x_k = \cX_{i_{k}} \cap \cP_{n}(x_n)$ for $0 \leq k \leq n$. We   show the existence of a  circumscribed sphere for the set  $\{ x_0,\ldots, x_n\}$, and that   $\Delta(x_0 , \ldots, x_n) \in \Delta'_{\F}(\cX)$.

\begin{prop}\label{prop-inductsphere}
For   $x_n \in \cX_{i_n}$ with $(x_0 , \ldots, x_n)$ defined as above,  there exists $r(x_0, \ldots , x_n)$ and 
$\omega(x_0, \ldots , x_n) \in \cP_{n}(x_n) $
  such that
\begin{equation}\label{eq-inductsphere1}
\{x_0, \ldots , x_n\} ~ \subset ~ S_{\F}(\omega(x_0, \ldots , x_n), r(x_0, \ldots , x_n)) ~ \cap ~ \cN(x_n)
\end{equation}
Moreover, the center satisfies, for $\ve_3$ defined by \eqref{eq-epsilon3},
\begin{equation}\label{eq-inductsphere2}
\dF(\omega(x_0, \ldots , x_n) , \omega'(y_0 , \ldots, y_n) ) ~ \leq ~ \ve_3 \rF/2
\end{equation}
where $\omega'(y_0 , \ldots, y_n) = \phi_{i_{n}}(y_n, x_n)( \omega(y_0, \ldots , y_n))$ is the translate for the center of the circumscribed sphere for the $n$-simplex $\Delta(y_0 , \ldots, y_n) \in \Delta^{(n)}_{\F}(\cX)$. In particular, this implies
\begin{equation}\label{eq-leafvar8}
| \, r(x_0, \ldots , x_n) - r(y_0, \ldots , y_n) \, | < \ve_3 \rF/2 + 2\ve_0 \rF < \ve_3 \rF
\end{equation}
\end{prop}
\proof
By rearranging the order of the vertices if necessary, we may assume that there are indices
$i_0 < i_1 < \cdots < i_n \leq p_*$ and points $\xi_{i_{k}} \in \Xi_{p_*} \subset L_0$ such that
$y_{k} = \cX_{i_{k}} \cap \cP_n(y_n)$ for $0 \leq k \leq n$.

Let $\omega = \omega(y_0, \ldots , y_n) \in \cP_{n}(y_n) $ denote the center of the circumscribed sphere for
$\{y_0 , \ldots, y_n\}$, and let $r(y_0, \ldots , y_n)$ denote its radius.
Then $d_1'/2 \leq r(y_0 , \ldots, y_n) \leq d_2'$ as $\cN(y_n)$ is $d_2'$-dense and $d_1'$-separated. Note that this implies
$\{y_0 , \ldots, y_n, \omega\} \subset D(y_n, \rF/5)$.

By Proposition~\ref{prop-inductrobust}, the set $\{x_0 , \ldots , x_n \} \subset \cP_n(x_n)$ is $\wtrho_n$-robust.

Let $\phi_{i_{n}}(y_n, x_n) \colon \cP_{n}(y_n) \to \cP_n(x_n)$ be the
transverse transport map for the chart $\vp_{i_n}$.

Let $\omega'(y_0, \ldots , y_n) = \phi_{i_{n}}(y_n, x_n)( \omega)$ denote the translation of $\omega(y_0, \ldots , y_n)$ to $ \cP_{n}(x_n)$.

For $0 \leq j \leq n$, we have the radius equalities $\dF(y_{j} , \omega ) = r(y_0, \ldots , y_n)$, hence
by Lemma~\ref{lem-secdiv},
\begin{equation}\label{eq-newcenter1}
r(y_0 , \ldots, y_n) - 2\ve_0 \rF ~ \leq ~ \dF(x_j , \omega'(y_0, \ldots , y_n) ) ~ \leq ~ r(y_0 , \ldots, y_n) + 2\ve_0 \rF
\end{equation}
Indeed,  note that $x_{\ell}$ is the transverse transport of $y_{\ell}$ for the coordinate system $\vp_{i_{\ell}}$, while $\omega'(y_0, \ldots , y_n) $ is the transport of $\omega(y_0, \ldots , y_n) $ for the coordinate system $\vp_{i_n}$ and $i_{\ell} \ne i_n$. Thus, we must use \eqref{eq-error1} in place of the sharper estimate \eqref{eq-error2}. Similarly, for $0 \leq j \ne k \leq n$, we have
\begin{equation}\label{eq-newcenter2}
d_1' ~ \leq ~ \dF(x_j ,x_k) ~ \leq ~ 2d_2' + 2\ve_0 \rF < 2d_2
\end{equation}
It follows that we also have $\{x_0 , \ldots, x_n, \omega'(y_0, \ldots , y_n) \} \subset D(x_n, \rF/5)$.

The first step is to construct a circumscribed sphere with center $\vec{\omega}(\vec{v}_0, \ldots, \vec{v}_{n})$ for the linearized problem in the tangent space $T_{x_n}\F$, and then modify the construction to obtain a circumscribed sphere with center $\omega(x_0, \ldots , x_n) \in \cP_n(x_n)$ for the leafwise metric.

Choose $\xi \in \cP_n(x_n)$ so that $\{x_0 , \ldots , x_n, \omega'(y_0, \ldots , y_n) \} \subset B_{\F}(\xi, 2d_2)$.
Let $\whu = \{\vec{u}_1, \ldots, \vec{u}_n\} \subset T_{\xi}\F$ be an orthonormal frame,
with corresponding geodesic coordinates $\psi^g_{\xi, \whu}$ about $\xi$.

Set $\vec{v}_k = (\psi^g_{\xi, \whu})^{-1}(x_k)$ for $0 \leq k \leq n$, then
$\{\vec{v}_0, \ldots , \vec{v}_{n} \} \subset \mR^n$ is $\wtrhop_{n+1}$-robust by Proposition~\ref{prop-robustnets}.

Now set $\vop(y_0, \ldots , y_n) = (\psi^g_{\xi, \whu})^{-1}(\omega'(y_0, \ldots , y_n))$. Then by
Lemma~\ref{lem-geodesicisom} and \eqref{eq-newcenter1} we have
\begin{equation}\label{eq-newcenter5}
r(y_0 , \ldots, y_n) - 3\ve_0 \rF ~ \leq ~ \| \vec{v}_j - \vop(y_0, \ldots , y_n) \| ~ \leq ~ r(y_0 , \ldots, y_n) + 3\ve_0 \rF
\end{equation}
while Lemma~\ref{lem-geodesicisom} and  \eqref{eq-newcenter2} implies, for $0 \leq j \ne k \leq n$,
\begin{equation}\label{eq-newcenter6}
d_1 < d_1' - \ve_0 \rF ~ \leq ~ \| \vec{v}_j - \vec{v}_k \| ~ \leq ~  d_2' + 3\ve_0 \rF <  d_2
\end{equation}

We can thus apply Proposition~\ref{prop-varestimate3} for $e_1 = d_1$, $e_2 = d_2 = 2d_1$,
$\rho = \wtrhop_{n+1} > 3\ve_2\rF/2$ and $C_1 = 3\ve_0 \rF$ to conclude that there exists a circumscribed sphere
$S(\omega(\vec{v}_0, \ldots , \vec{v}_n), r(\vec{v}_0, \ldots , \vec{v}_n)) \subset D(\rF)$ such that
\begin{eqnarray}
\left \| \omega(\vec{v}_0, \ldots, \vec{v}_{n}) - \vop(y_0, \ldots , y_n) \right \| ~ & \leq & ~ 3\ve_0 \cdot \left\{ n^{3/2} (2d_2)^{n-1}/ (3\ve_2 \rF/2)^{n-1} \right\} \cdot \rF \nonumber\\
~ & \leq & ~ 3\ve_0 \cdot \left\{ n^{3/2} \cdot (4/15\ve_2 )^{n-1}\right\} \cdot \rF \label{eq-newcenter7}
\end{eqnarray}
That is, the vector $ \omega(\vec{v}_0, \ldots, \vec{v}_{n}) $ is a solution of the linearized problem of finding the center of a circumscribed sphere, and \eqref{eq-newcenter7} estimates the Euclidean distance to the translated center.

The task now is to convert this approximate answer to a solution for the leafwise metric. As before,
let $\wtd$ denote the distance function on $ D(\rF)$ induced from $\dF$ by $\psi^g_{\xi, \whu} \colon D(\rF/2) \to L_{\xi}$. Then by Lemma~\ref{lem-geodesicisom},
\begin{equation}\label{eq-leafvar1}
| \wtd(\, \vec{a}, \vec{b} \, ) - \|\vec{a} - \vec{b}\| | \leq \ve_0 \, \rF \quad , \quad \text{for all} ~ \vec{a}, \vec{b} \in D(\rF/2)
\end{equation}

Introduce the equidistant submanifolds for the metric $\wtd$,
\begin{equation}\label{eq-leafvar2}
\cH( \vec{v}_{j} , \vec{v}_{k} ) = \{ \vec{z} \in D(\rF/2) \mid \wtd(\vec{z}, \vec{v}_{j}) = \wtd(\vec{z}, \vec{v}_{k}) \}
\end{equation}
and the ``thickened'' equidistant sets for the leafwise metric, for $\e > 0$,
\begin{equation}\label{eq-leafvar3}
\cH( \vec{v}_{j} , \vec{v}_{k} ; \e) = \{ \vec{z} \in D(\rF/2) \mid - \e \leq \wtd(\vec{z}, \vec{v}_{j}) - \wtd(\vec{z}, \vec{v}_{k}) \leq \e \}
\end{equation}
\begin{equation}\label{eq-leafvar4}
\cB(\vec{v}_{0} , \ldots , \vec{v}_{n} ; \e) = \cH( \vec{v}_{0} , \vec{v}_{n} ; \e) \cap \cdots \cap \cH( \vec{v}_{n-1} , \vec{v}_{n} ; \e)
\end{equation}
Then \eqref{eq-newcenter1} implies the translation
$\ds \vop(y_0, \ldots , y_n) \in \cB(\vec{v}_{0} , \ldots , \vec{v}_{n} ; 4 \ve_0 \rF)$, so
this set is not empty. The key idea is to obtain a bound for its diameter,
from which the proof of Proposition~\ref{prop-inductsphere} follows. To this end,
define the set of approximate solutions of the linearized problem by
\begin{equation}\label{eq-leafvar4a}
B(\vec{v}_{0} , \ldots , \vec{v}_{n} ; \e) = \{ \vec{z} \in D(\rF/2) \mid - \e \leq \| \vec{z} - \vec{v}_{j} \| - \| \vec{z} - \vec{v}_{n} \| \leq \e ~ , ~ 0 \leq j < n\}
\end{equation}
Note that the actual solution satisfies $\omega(\vec{v}_0, \ldots, \vec{v}_{n}) \in B(\vec{v}_{0} , \ldots , \vec{v}_{n} ; \e)$ for all $\e > 0$.

\medskip

\begin{lemma}\label{lem-boxdiam}
$\ds \cB(\vec{v}_{0} , \ldots , \vec{v}_{n} ; \e) \subset B(\vec{v}_{0} , \ldots , \vec{v}_{n} ; \e + 4 \ve_0 \rF) $
\end{lemma}
\proof
Using \eqref{eq-leafvar1} for $ \vec{z} \in D(\rF/2)$, we have that
\begin{equation}\label{eq-leafvar5}
\left| ( \| \vec{z} - \vec{v}_{j} \| - \| \vec{z} - \vec{v}_{k} \| ) \right| - 2 \ve_0 \rF ~ \leq ~
\left| \wtd(\vec{z}, \vec{v}_{j}) - \wtd(\vec{z}, \vec{v}_{k}) \right| ~ \leq ~
\left| ( \| \vec{z} - \vec{v}_{j} \| - \| \vec{z} - \vec{v}_{k} \| ) \right| + 2 \ve_0 \rF
\end{equation}
and the claim follows.
\endproof

Thus, we now have $\ds
\vop(y_0, \ldots , y_n) \in \cB(\vec{v}_{0} , \ldots , \vec{v}_{n} ; 4 \ve_0 \rF) \subset B(\vec{v}_{0} , \ldots , \vec{v}_{n} ; 8 \ve_0 \rF)$.

\begin{lemma}\label{lem-centererrorest}
Let $\vec{z} \in B(\vec{v}_{0} , \ldots , \vec{v}_{n} ; 8 \ve_0 \rF)$, then
\begin{equation}\label{eq-leafvar6}
\| \vec{z} - \omega(\vec{v}_0, \ldots , \vec{v}_n) \| ~ \leq ~ 32 \ve_0 \cdot \left\{ n^{3/2} \cdot (4/15\ve_2 )^{n-1}\right\} \cdot \rF
\end{equation}
\end{lemma}
\proof
Using the notation of Propositions~\ref{prop-varestimate2} and \ref{prop-varestimate3},
with $\vec{v}_j$ in place of $\vec{z}_j$ and $\vec{z}$ in place of $\omega$,
and $e_1 = d_1$, $e_2 = d_2 = 2d_1$, $\rho = \wtrhop_{n+1} > 3\ve_2 \rF/2$ and $C_1 = 8\ve_0 \rF$, then
$\vec{\zeta} = \vec{z} - \omega(\vec{v}_0, \ldots , \vec{v}_n)$ is a solution of the matrix inequality
\begin{equation}\label{eq-matrixineq3a}
\bV \cdot \vec{\zeta} \in B(0, 2\sqrt{n} \cdot d_2 \cdot 8 \ve_0 \rF)
\end{equation}
Then by \eqref{eq-Vest3}, we have the estimate $\|\bV^{-1} \| \leq n \cdot (2d_2)^{n-1}/d_1 ( 3\ve_2 \rF/2)^{n-1}$
which yields
\begin{eqnarray*}
\| \vec{z} - \omega(\vec{v}_0, \ldots , \vec{v}_n) \| ~ & \leq & ~ \left\{ n \cdot (2d_2)^{n-1}/d_1(3\ve_2 \rF/2)^{n-1}\right\} \cdot \left\{ 2\sqrt{n} \cdot d_2 \cdot 8 \ve_0 \rF \right\} \\
~ & \leq & ~ \ve_0 \cdot \left\{ 32 n^{3/2} \cdot (4/15\ve_2 )^{n-1}\right\} \cdot \rF
\end{eqnarray*}
where we use that $d_1 = \rF/10$ and $d_2 = 2\rF/10$ to simplify, yielding \eqref{eq-leafvar6}.
\endproof

It follows from Lemmas~\ref{lem-boxdiam} and \ref{lem-centererrorest} that
the closed set $\cB(\vec{v}_{0} , \ldots , \vec{v}_{n} ; 4 \ve_0 \rF)$ is bounded,
and is non-empty as $\ds \vop(y_0, \ldots , y_n) \in \cB(\vec{v}_{0} , \ldots , \vec{v}_{n} ; 4 \ve_0 \rF)$.
Thus, the intersection
\begin{equation}\label{eq-leafvar7}
\wtomega( \vec{v}_{0} , \ldots , \vec{v}_{n} ) = \cH( \vec{v}_{0} , \vec{v}_{n} ) \cap \cdots \cap \cH( \vec{v}_{n-1} , \vec{v}_{n} ) \subset \cB(\vec{v}_{0} , \ldots , \vec{v}_{n} ; 4 \ve_0 \rF)
\end{equation}
is non-empty by transversality of the submanifolds $\ds \cH( \vec{v}_{j} , \vec{v}_{n} )$.
Moreover, \eqref{eq-leafvar6} implies that
\begin{equation}\label{eq-leafvar6alt}
\| \wtomega( \vec{v}_{0} , \ldots , \vec{v}_{n} )- \omega(\vec{v}_0, \ldots , \vec{v}_n) \| ~ \leq ~
32 \ve_0 \cdot \left\{ n^{3/2} \cdot (4/15\ve_2 )^{n-1}\right\} \cdot \rF
\end{equation}
Combine this with the estimate \eqref{eq-newcenter7} to obtain
\begin{equation}\label{eq-leafvar6bat}
\| \wtomega( \vec{v}_{0} , \ldots , \vec{v}_{n} ) - \vop(y_0, \ldots , y_n) \| ~ \leq ~
35 \ve_0 \cdot \left\{ n^{3/2} \cdot (4/15\ve_2 )^{n-1}\right\} \cdot \rF
\end{equation}
Then set
\begin{equation}\label{eq-spherecenter}
\omega(x_{0} , \ldots , x_{n}) = \psi^g_{\xi, \whu}(\wtomega( \vec{v}_{0} , \ldots , \vec{v}_{n} ))
~ , ~ r(x_{0} , \ldots , x_{n}) = \dF( x_0 , \omega(x_{0} , \ldots , x_{n}) )
\end{equation}
so we have
$\{x_{0} , \ldots , x_{n} \} \subset S_{\F}(\omega( x_{0} , \ldots , x_{n} ) , r(x_{0} , \ldots , x_{n}))$ as desired.

Recall that $\omega'(y_0, \ldots , y_n)) = \psi^g_{\xi, \whu}(\vop(y_0, \ldots , y_n))$,
then by Lemma~\ref{lem-geodesicisom} we have
\begin{equation}\label{eq-leafvar6cat}
\dF\left( \omega(x_{0} , \ldots , x_{n}) , \omega'(y_0, \ldots , y_n) \right) ~ \leq ~
\ve_0 \cdot \left\{1 + 35 n^{3/2} \cdot (4/15\ve_2 )^{n-1}\right\} \cdot \rF < \ve_3 \rF/2
\end{equation}
where the bound by $\ve_3 \rF/2$ follows from \eqref{eq-epsilon0}.

Finally, the estimate \eqref{eq-leafvar8} follows from
\begin{eqnarray*}
| \, r(x_0, \ldots , x_n) - r(y_0, \ldots , y_n) \, | & = & | \, \dF(x_n , \omega(x_{0} , \ldots , x_{n})) - \dF(y_n , \omega(y_0, \ldots , y_n)) \, | \\
& \leq & | \, \dF(x_n , \omega(x_{0} , \ldots , x_{n})) - \dF(x_n , \omega'(y_0, \ldots , y_n)) \, | + 2\ve_0 \rF \\
& \leq & | \, \dF( \omega(x_{0} , \ldots , x_{n}) , \omega'(y_0, \ldots , y_n)) \, | + 2\ve_0 \rF \\
& \leq & \ve_3 \rF/2 + 2\ve_0 \rF < \ve_3 \rF
\end{eqnarray*}
This completes the proof of Proposition~\ref{prop-inductsphere}.
\endproof

\medskip

\subsection{Stability}
We have now established that for a simplex $\Delta(y_0 , \ldots, y_n) \in \Delta'_{\F}(\cX)$, if $\{x_0 , \ldots , x_n \}$ is a transverse translate of the set $\{y_0 , \ldots, y_n\}$,
then $\{x_0 , \ldots , x_n \}$ is $\wtrho_n > 3\ve_2/2$ robust and admits a circumscribed sphere whose radius varies according to the estimate \eqref{eq-inductsphere2}.
It remains to show that $\cX$ is stable,  that is, $\Delta(x_0 , \ldots, x_n) \in \Delta'_{\F}(\cX)$. The only ingredient left   to show is that the circumscribed sphere for the set $\{x_0 , \ldots, x_n\}$ does not contain other points of $\cX$ in its interior.

\begin{prop}\label{prop-inductsimplex}
Let $\Delta(y_0 , \ldots, y_n) \in \Delta^{(n)}_{\F}(\cX)$. Assume given $1 \leq i_0 < i_1 < \cdots < i_n \leq p_*$
and $\xi_{i_{k}} \in \Xi_{p_*} \subset L_0$ such that
$y_{k} = \cX_{i_{k}} \cap \cP_{\theta_{i_{n}}}(y_n)$.
Then for all $x_n \in \cX_{i_n}$, $\Delta(x_0 , \ldots, x_n) \in \Delta'_{\F}(\cX)$.
\end{prop}
\proof
Let $\omega(y_0, \ldots , y_n)\in \cP_n(y_n)$ be the center of the circumscribed sphere of radius $r(y_0, \ldots , y_n)$.
Then it is given that for all $\xi \in \cN(y_n) - \{y_{0} , \ldots , y_{n}\}$ we have that
$\ds \dF(\xi , \omega(y_0, \ldots , y_n)) > r(y_0, \ldots , y_n)$.

We must show that the circumscribed sphere for the set $\{x_0 , \ldots, x_n\}$ obtained in Proposition~\ref{prop-inductsphere},
with center $\omega(x_{0} , \ldots , x_{n})$ and radius $r(x_{0} , \ldots , x_{n})$, contains no points of $\cX$ in its interior. That is, we must show that
\begin{equation}\label{eq-leafvar10}
\dF(\xi' , \omega(x_{0} , \ldots , x_{n}) > r(x_{0} , \ldots , x_{n}) \quad \text{for all} ~ \xi' \in \cN(x_n) - \{x_{0} , \ldots , x_{n}\}
\end{equation}

Let $n \leq m \leq p_*$ be the largest $m$ such that the condition \eqref{eq-leafvar10} holds for all $\Delta(y_0 , \ldots, y_n) \in \Delta'_{\F}(\cX)$ with $i_n \leq m$. If $m = p_*$ then we are done, so assume that $m < p_*$ and we show this leads to a contradiction. So we assume that we are given a simplex $\Delta(y_0 , \ldots, y_n)$ with $i_n = m+1$, such that
there is some $x_n \in \cX_{i_n}$ and $\xi' \in \cN(x_n) - \{x_{0} , \ldots , x_{n}\}$ for which \eqref{eq-leafvar10} fails.

First, consider the case where there exists
$\xi' \in \cN(x_n) - \{x_{0} , \ldots , x_{n}\}$ such that
\begin{equation}\label{eq-leafvar11}
\dF(\xi' , \omega(x_{0} , \ldots , x_{n}) \leq r(x_{0} , \ldots , x_{n}) - 2\ve_3 \rF
\end{equation}
Let $1 \leq q \leq p_*$ be such that $\xi' \in \cX_{q}$ and set $\xi = \fZ(\xi', \theta_{i_q}, V_0) \cap \cP_{n}(y_n) \in \cN(y_n)$. Then
\begin{eqnarray*}
\dF(\xi , \omega(y_0, \ldots , y_n)) & < & \dF(\xi' , \omega'(y_0, \ldots , y_n)) + 2\ve_0 \rF \quad \text{by Lemma~\ref{lem-geodesicisom}} \\
& < & \dF(\xi' , \omega(x_0, \ldots , x_n)) + 2\ve_0 \rF + \ve_3 \rF/2 \quad \text{by \eqref{eq-inductsphere2}} \\
& < & r(x_{0} , \ldots , x_{n}) + 2\ve_0 \rF + \ve_3 \rF/2 - 2\ve_3 \rF \quad \text{by \eqref{eq-leafvar11}} \\
& < & r(y_{0} , \ldots , y_{n}) - 3\ve_3 \rF/2 + 2\ve_0 \rF + (2\ve_0 + \ve_3/2) \rF \quad \text{by \eqref{eq-leafvar8}} \\
& < & r(y_{0} , \ldots , y_{n}) + (2\ve_0 - \ve_3 ) \rF < r(y_{0} , \ldots , y_{n}) \quad \text{by (\ref{eq-epsilon0}.4)}
\end{eqnarray*}
which contradicts the hypothesis that $\Delta(y_0 , \ldots, y_n) \in \Delta'_{\F}(\cX)$.
Thus, we may assume that
\begin{equation}\label{eq-leafvar12}
r(x_{0} , \ldots , x_{n}) - 2 \ve_3 \rF < \dF(\xi' , \omega(x_{0} , \ldots , x_{n})) \leq r(x_{0} , \ldots , x_{n})
\end{equation}

Let $1 \leq q \leq p_*$ be the least such $q$ such that there exists $x_n \in \cX_{i_n}$ and
\eqref{eq-leafvar12} holds for $\xi' \in \cX_q$. Note that $q \ne i_{k}$ for $0 \leq k \leq n$, and $\cX_q = \fZ(\xi_q, \theta_{i_q}, V_0)$ for some
$\xi_q \in \Xi$, so that $\xi' = \fZ(\xi_q, \theta_{i_q}, V_0) \cap \cP_{n}(x_n)$.

We now use that $\xi_q$ was chosen inductively to avoid the annular $2\ve_1 \rF$-thickening of the circumscribed sphere for each $n$-simplex in $\Omega^{(n)}(\cN(\xi_q))$.

First, we assume that $q > i_n$. For this subcase, we transfer the problem to the plaque $\cP_n(\xi_q)$.
Set $z_k = \cX_{i_{k}} \cap \cP_{n}(\xi_q)$ for $0 \leq k \leq n$.
Note that that $\{z_0, \ldots , z_n\}$ admits a circumscribed sphere by Proposition~\ref{prop-inductsphere}, with center
$\omega(z_0, \ldots , z_n)$ which satisfies
\begin{equation}\label{eq-centerest10}
\dF(\omega(z_0, \ldots , z_n) , \omega'(y_0 , \ldots, y_n) ) ~ \leq ~ \ve_3 \rF/2
\end{equation}
where $\omega'(y_0 , \ldots, y_n) = \phi_{i_{n}}(y_n, z_n)(\omega(y_0, \ldots , y_n))$. Let $r(z_0, \ldots , z_n)$ denote the radius of the sphere, which by \eqref{eq-leafvar8} satisfies
\begin{equation}\label{eq-radest1}
r(y_0 , \ldots , y_n) - \ve_3 \, \rF ~ \leq ~ r(z_0 , \ldots , z_0) ~ \leq ~ r(y_0 , \ldots , y_n) + \ve_3 \, \rF
\end{equation}
We claim that
$\Delta(z_0 , \ldots, z_n) \in \Delta'_{\F}(\whcX_{i_{q-1}})$.
If not, then there exists $\eta' \in \whcX_{i_{q-1}} \cap \cP_n(z_n)$ with
$\ds \dF(\eta' , \omega(z_{0} , \ldots , z_{n})) \leq r( z_{0} , \ldots , z_{n})$. This contradicts the minimality of the choice of $q$ above.
Thus, as $\xi_q$ was chosen to satisfy the inequality \eqref{eq-unifreg2}, we have the estimate
\begin{equation}\label{eq-radest2}
\dF(\xi_q , \omega(z_{0} , \ldots , z_{n})) ~ > ~ r(z_0 , \ldots , z_0) + 2\ve_1 \, \rF
\end{equation}
On the other hand, $x_n$ was chosen so that for the circumscribed sphere with center $\omega(x_0, \ldots ,x_n)$ and radius $r(x_0, \ldots , x_n)$ we have the inequality \eqref{eq-leafvar12} above.

Apply Proposition~\ref{prop-inductsphere} to the cases $x_n \in \cX_{i_n}$ and also $z_n \in \cX_{i_n}$ to obtain the estimates
\eqref{eq-leafvar8} for both. Together, they imply
\begin{equation}\label{eq-radest3}
r(x_{0} , \ldots , x_{n}) - 2 \ve_3 \rF ~ \leq ~ r(z_{0} , \ldots , z_{n}) ~ \leq ~ r(x_{0} , \ldots , x_{n}) + 2 \ve_3 \rF
\end{equation}
Also, \eqref{eq-leafvar12} and \eqref{eq-radest3} imply
\begin{equation}\label{eq-radest4}
\dF(\xi' , \omega(x_{0} , \ldots , x_{n})) - 2 \ve_3 \rF \leq r(x_{0} , \ldots , x_{n}) - 2 \ve_3 \rF ~ \leq ~ r(z_{0} , \ldots , z_{n})
\end{equation}
which for $ \omega'(x_{0} , \ldots , x_{n}) = \phi_{i_n}(x_n , z_n)(\omega(x_{0} , \ldots , x_{n}))$ yields
\begin{equation}\label{eq-radest5}
\dF(\xi_q , \omega'(x_{0} , \ldots , x_{n})) ~\leq ~ r(z_{0} , \ldots , z_{n}) + 2 \ve_3 \rF + \ve_0 \rF
\end{equation}
and thus by \eqref{eq-centerest10} and its corresponding version for $x_n$ yields
\begin{eqnarray}
\dF(\xi_q , \omega(z_{0} , \ldots , z_{n})) ~ & \leq & ~ \dF(\xi_q , \omega'(x_{0} , \ldots , x_{n}) ) + \dF( \omega'(x_{0} , \ldots , x_{n}) , \omega(z_{0} , \ldots , z_{n})) \nonumber \\
~ & \leq & ~ ( r(z_{0} , \ldots , z_{n}) + (2 \ve_3 \rF + \ve_0 \rF ) + ( \ve_3 \rF + 2 \ve_0 \rF) \nonumber \\
~ & \leq & ~ r(z_{0} , \ldots , z_{n}) + 4\ve_3 \rF \label{eq-radest6}
\end{eqnarray}
which by the choice $\ve_3 < \ve_1/2$ in \eqref{eq-epsilon3} contradicts \eqref{eq-radest2}.
Thus, the case $q > i_n$ is not possible.

Finally, consider the case where $q < i_n$. That is,
the smallest $q$ such that there exists $x_n \in \cX_{i_n}$ and \eqref{eq-leafvar12} holds for some $\xi' \in \cX_q$ occurs for $q < i_n$.
This means that in the process of constructing $\cX$, we have chosen a point $\xi_q$ which has distance greater than $2 \ve_1 \rF$ from all previously circumscribed spheres for the net $\cN_{q-1}$, but when we add the point $\xi_{i_n}$ the Delaunay triangulation $\ds \Delta'_{\F}(\cN_{i_n})$ abruptly changes on some leaves.
The translates of $\xi_{i_n}$ are contained both inside and outside of circumscribed spheres, as the translates of $\xi_q$ also wander inside and outside. We show this is impossible, due to the choice of $\ve_1$ and of the constants $\ve_3$ and $\ve_4$ which control how much the centers of circumscribed spheres ``wander'' for transverse variation at most $r_*$.

Recall, we assume there is given $\Delta(y_0 , \ldots, y_n) \in \Delta'_{\F}(\cX)$ and $x_n \in \cX_{i_n}$ such that the transverse translate $\{x_0 , \ldots , x_n \}$ of the set $\{y_0 , \ldots, y_n\}$ is $\wtrho_n > 3\ve_2/2$ robust. Thus by Proposition~\ref{prop-inductsphere}, there is a circumscribed sphere with center $\omega(x_0, \ldots , x_n)$ and radius $r(x_0, \ldots , x_n)$ which satisfy
\begin{equation}\label{eq-inductsphere2aaa}
\dF(\omega(x_0, \ldots , x_n) , \omega'(y_0 , \ldots, y_n) ) ~ \leq ~ \ve_3 \rF/2
\end{equation}
\begin{equation}\label{eq-leafvar13}
| \, r(x_0, \ldots , x_n) - r(y_0, \ldots , y_n) \, | < \ve_3 \rF
\end{equation}
where $\omega'(y_0 , \ldots, y_n) = \phi_{i_{n}}(y_n, x_n)( \omega(y_0, \ldots , y_n))$ is the translate for the center of the circumscribed sphere for the $n$-simplex $\Delta(y_0 , \ldots, y_n) \in \Delta'_{\F}(\cX)$.
There is also given $1 \leq q < i_n$ so that the translate $\xi' = \cX_q \cap \cP_{n}(x_n) \in \cN(x_n)$ satisfies \eqref{eq-leafvar12}.
\begin{equation}\label{eq-leafvar14}
r(x_{0} , \ldots , x_{n}) - 2 \ve_3 \rF < \dF(\xi' , \omega(x_{0} , \ldots , x_{n})) \leq r(x_{0} , \ldots , x_{n})
\end{equation}

Let $\{x_0', \ldots , x_n'\} = \{x_0, \ldots , x_{n-1}, \xi'\}$ denote a reordering of the set so that $x_k' = \cX_{i_k'} \cap \cP_n(x_n)$ for $0 \leq k \leq n$
with $1 \leq i_0' < \cdots < i_n' \leq p_*$. Then these points satisfy, for $0 \leq k \leq n$,
\begin{equation}
r(x_{0} , \ldots , x_{n}) - 2 \ve_3 \rF \leq \dF(x_k' , \omega(x_{0} , \ldots , x_{n})) \leq r(x_{0} , \ldots , x_{n})
\end{equation}
The proof of Proposition~\ref{prop-inductrobust} applied to the set $\{x_0', \ldots , x_n'\} $ yields that the collection is $\wtrho_n$-robust, so admits a circumscribed sphere by Proposition~\ref{prop-inductsphere}, with center $\omega(x_0', \ldots , x_n')$ and radius $r(x_0' , \ldots , x_n' )$.
From the proof of Proposition~\ref{prop-inductsphere}, we have the estimates
\begin{equation}\label{eq-leafvar15}
\dF(\omega(x_0', \ldots , x_n') ,\omega(x_0, \ldots , x_n) ) ~ \leq ~ \ve_3 \rF/2
\end{equation}
\begin{equation}\label{eq-leafvar16}
| \, r(x_0', \ldots , x_n') - r(x_0, \ldots , x_n) \, | < \ve_3 \rF
\end{equation}
Thus, combining \eqref{eq-leafvar15} and \eqref{eq-leafvar16}, for $\zeta' = x_n$, we obtain
\begin{equation}\label{eq-leafvar17}
\dF(\zeta', \omega(x_0', \ldots , x_n') ) ~ \leq ~ r(x_0' , \ldots , x_n') + 3\ve_3 \rF/2
\end{equation}
Now let $\zeta = \xi_{i_n} \in \cX_{i_n}$. The last step is to translate the points $\{x_0', \ldots , x_n'\}$ to the plaque $\cP_n(\zeta)$, to obtain points
$z_k' = \cX_{i_k'} \cap \cP_n(\zeta)$. Then $\{z_0' , \ldots , z_n'\}$ is $\wtrhop_n$-robust by Proposition~\ref{prop-inductrobust}, and admits a circumscribed sphere with center $\omega(z_0', \dots , z_n')$ and radius $r(z_0' , \ldots , z_n')$ by Proposition~\ref{prop-inductsphere}.
Moreover, this center and radius satisfy
\begin{equation}\label{eq-leafvar18}
\dF(\omega(z_0', \ldots , z_n') , \omega'(x_0', \ldots , x_n') ) ~ \leq ~ \ve_3 \rF/2
\end{equation}
\begin{equation}\label{eq-leafvar19}
| \, r(z_0', \ldots , z_n') - r(x_0', \ldots , x_n') \, | < \ve_3 \rF
\end{equation}
Combining \eqref{eq-leafvar17}, \eqref{eq-leafvar18} and \eqref{eq-leafvar19}, we obtain
\begin{eqnarray}
\dF(\zeta, \omega(z_0', \ldots , z_n') ) ~ & \leq & ~ \dF(\zeta, \omega'(x_0', \ldots , x_n') ) + \dF(\omega'(x_0', \ldots , x_n'), \omega(z_0', \ldots , z_n') ) \nonumber \\
~ & \leq & ~ \dF(\zeta', \omega(x_0', \ldots , x_n') ) + 2 \ve_0 \rF + \ve_3 \rF/2 \nonumber \\
~ & \leq & ~ r(x_0' , \ldots , x_n') + 3\ve_3 \rF/2 + 2 \ve_0 \rF + \ve_3 \rF/2 \nonumber \\
~ & \leq & ~ r(z_0', \ldots , z_n') + \ve_3 \rF + 3\ve_3 \rF/2 + 2 \ve_0 \rF + \ve_3 \rF/2 \nonumber \\
~ & < & ~ r(z_0', \ldots , z_n') + 4\ve_3 \rF + 2\ve_0 \rF \label{eq-leafvar20}
\end{eqnarray}
By the choice of $\ve_3$ in \eqref{eq-epsilon3} we have $ 4\ve_3 + 2\ve_0 < 2\ve_1$, so that \eqref{eq-leafvar20} contradicts the choice of $\zeta = \xi_{i_n}$ to satisfy
$$ \dF(\zeta, \omega(z_0', \ldots , z_n') ) \geq r(z_0', \ldots , z_n') + 2\ve_1\rF$$
Thus, the case $q < i_n$ again leads to a contraction.
This completes the proof of Proposition~\ref{prop-inductsimplex}.
\endproof

Thus, we have shown that   $\cX$ as defined is a  nice stable  transversal $\cX$ for $\ds  \fN(\whK_x, V_{x})$.

 \section{Proofs of Theorem~\ref{thm-foliated} and Theorem~\ref{thm-tessel}} \label{sec-proofs}

\subsection{Proof of Theorem~\ref{thm-foliated}} \label{subsec-Fol}

Let $\fM$ be an equicontinuous matchbox manifold. By Theorem~\ref{thm-minimal}, the dynamics of its holonomy pseudogroup is minimal.
By Theorem~\ref{thm-invariants}, given any  $w_0 \in    \fT_*$,  say with $w_0 \in \fT_{i_0}$, there exists a descending chain of clopen subsets
$$ w_0 \in  \cdots \subset V_{\ell +1} \subset V_{\ell} \subset \cdots V_0 \subset \fT_{i_0}$$
such that for all $\ell \geq 0$, 
$w_0 \in V_{\ell}$ and ${\rm diam}_{\fX}(V_{\ell}) < \dTU/2^{\ell}$.
Here, $\dTU$ is the constant of equicontinuity defined in Proposition~\ref{prop-uniformdom}, so that each set $V_{\ell}$ is in the domain of the holonomy of any path starting at $x_0$ where $x_0 = \tau_{i_0}(w_0)$. Thus,  we can form the associated Thomas tube $\wtfN(V_{\ell})$ as defined by  \eqref{eq-thomas}.

 Choose $w_0$ corresponding to a leaf $L_0$ without holonomy in $\fM$, let $\wtL_0 \to L_0$ be the holonomy cover, which is a diffeomorphism, and lift $w_0$ to a basepoint $\wtw_0 \in \wtL_0$. Then let $\wtcM_0$ be the net in $\wtL_0$ defined in section~\ref{sec-microbundles}, with $\wtw_0 \in \wtcM_0$. 
 For $\wtz \in \wtcM_0$ there is a holonomy transport map $h_{\wtz}$ defined by choosing a path from $\wtw_0$ to $\wtz$, and we define the holonomy transport of $V_{\ell}$ by $V_{\wtz}^{\ell} = h_{\wtz}(V_{\ell}) $ as defined by \eqref{def-basicblocks}. Then the collection of clopen sets $\{ V_{\wtz}^{\ell} \mid \wtz \in \wtcM_0 \}$, which covers the transverse space $\fT_*$, is actually a finite set. Thus, for each $\ell$ there exists a compact connected subset $\wtK_{\ell} \subset \wtL_0$ so that  $V_{\ell}$ is $\wtK$-admissible, and the Reeb neighborhood, $\ds \wtfN(\wtK_{\ell}, V_{\wtz}^{\ell})$ as defined by \eqref{eq-tesselspace},  maps onto $\fM$.

Let $r_* > 0$ be defined   by   \eqref{eq-transsec2}). Apply Theorem~\ref{prop-uniformdom} for   $\e = r_*$, to conclude that for   $\ell_0$ sufficiently large, 
all holonomy translates $V_{\wtz}^{\ell_0}$ of the set $V_{\ell_0}$ have diameter less than $r_*$. We can then apply the methods of section~\ref{sec-existence}   to construct a complete regular   $V_{\ell_0}$-transversal for $\wtK_{\ell_0}$. Then by Theorem~\ref{thm-stableapprox}, 
there exists a foliated homeomorphism into, 
 $\ds  \Phi \colon  \wtK_{\ell_0} \times V_{\ell_0} \to  \wtfN(V_{\ell_0})$, 
whose image contains the Reeb neighborhood $\ds \wtfN(\wtK_{\ell}, V_{\wtz}^{\ell})$. This defines the   transverse Cantor foliation $\wtcH_{\ell_0}$ on a neighborhood of $\ds \wtfN(\wtK_{\ell}, V_{\wtz}^{\ell})$ by the methods of section~\ref{sec-proofstableappr}.
The transversal $\whcX_{p_*} \subset \wtfN(V_{\ell_0})$  is invariant, in the sense of Definition~\ref{def-standtranscover2}, and thus is the pull-back of a transversal in $\fN(V_{\ell_0}) \subset \fM$. Thus,  the transverse Cantor foliation  $\wtcH_{\ell_0}$  is $\Pi$-equivariant, and descends to a transverse Cantor foliation $\cH_{\ell_0}$ on  $\fM$.

Note that the existence of $\cH_{\ell_0}$ on  $\fM$  is the result cited in   \cite[Theorem~8.3]{ClarkHurder2011b}.
Finally, note  that \cite[Proposition~8.4]{ClarkHurder2011b}
shows that   the quotient space $M \equiv \fM/\cH$ is an $n$-dimensional topological manifold, and \cite[Proposition~8.8]{ClarkHurder2011b}
implies that the projection to the leaf space $\fM \to \fM/\cH \cong M$ is a Cantor bundle.
This complete the proof of Theorem~\ref{thm-foliated}.

\subsection{Proof of Theorem~\ref{thm-tessel}} \label{subsec-BBL}

Let $\fM$ be a matchbox manifold, $L_x \subset M$ the leaf through $x \in \fM$, and $\wtL_x$ the holonomy covering of $L_x$. We are given a proper base $K_x \subset L_x$ so that there is a connected compact subset  
  $\wtK_x \subset \wtL_x$  such that the composition
$\ds \iota_x \colon \wtK \subset \wtL_x \to L_x \subset \fM$
is injective with image $K_x$. 
Introduce the set $\whK_x$ defined by \eqref{eq-translationK0} with diameter $\whR_K$. 
Let $w_x = \pi_{i_x}(x)$ be the image of $x$ in a transversal space $\fT_{i_x}$.

Let $r_* > 0$ be defined   by   \eqref{eq-transsec2}. Apply Proposition~\ref{prop-domest} for   $\e = r_*/2$, to conclude that there exists 
$\delta_* = \delta( r_*/2 , \whR_K)$, such that if $x \in V_x \subset \fT_{i_x}$ is a clopen neighborhood with $V_x \subset B_{\fX}(w_x, \delta_*)$ then
for any path with initial point $x$ and   length at  most $\whR_K$ the holonomy translate $h_{\gamma}(V_x)$ of the set $V_x$ has diameter less than $r_*$. 
Thus, $V_x$ is $\whK_x$-admissible, in the sense of Definition~\ref{def-admissibledisjoint}. 

Since $K_x$ is compact and the map $\Pi \colon K_x \to \fM$ is injective by assumption, by restricting the diameter of the clopen neighborhood $V_x$ further, we may assume that $V_x$ is $K_x$-disjoint. In particular, the map $\Pi \colon \fN(K_x, V_{x}) \to \fM$ is injective.

Now apply the methods of section~\ref{sec-proofstableappr}. Choose a basepoint $w_0 \in V_x$ such that the leaf $L_0$ it defines is without holonomy. 
Introduce the translated set $\whK_0 \subset \wtL_0$. Then by construction, each translate $h_{\gamma}(V_x)$ has diameter less than $r_*$ so we can   apply the methods of section~\ref{sec-existence}  to construct a complete regular   $V_{x}$-transversal for $\wtK_{0}$. Then by Theorem~\ref{thm-stableapprox}, 
there exists a foliated homeomorphism into, 
 $\ds  \Phi \colon  \wtK_{0} \times V_{x} \to  \fM$, 
whose image contains the Reeb neighborhood $\ds \fN(K_x, V_{x})$. Using that $\Pi \colon K_x \to \fM$ is injective, we can chose a clopen sub-neighborhood $x \in V_x' \subset V_x$ such that the restriction  $\ds  \Phi \colon  K_x \times V_{x}' \to  \fM$ is injective. 
This defines a  transverse Cantor foliation $\cH$ on a neighborhood of $\ds \fN(K_x, V_{x})$, which 
 completes the proof of Theorem~\ref{thm-tessel}.


\end{document}